\documentclass[3p,11pt,authoryear]{elsarticle}

\usepackage{lipsum}
\makeatletter
\def\ps@pprintTitle{%
 \let\@oddhead\@empty
 \let\@evenhead\@empty
 \def\@oddfoot{}%
 \let\@evenfoot\@oddfoot}
\makeatother

\usepackage{fancyhdr}
\usepackage{amssymb}
\usepackage{multirow}
\usepackage{mathtools} 
\let\oldforall\forall
\let\forall\undefined
\DeclareMathOperator{\forall}{\oldforall}
\usepackage{amsmath}
\usepackage{blindtext}
\usepackage{nccmath}
\usepackage{amsthm}
\usepackage{multicol}
\usepackage{subcaption}
\usepackage{float}
\usepackage{lscape}
\usepackage{xcolor,soul}
\usepackage[normalem]{ulem}
\usepackage{changepage}
\usepackage{algorithm}
\usepackage{comment}
\usepackage{float}
\usepackage{setspace}
\usepackage[figuresright]{rotating}
\usepackage{graphicx}
\usepackage{lineno}
\usepackage{natbib}
\usepackage{makecell}
\usepackage{bbm}
\usepackage{accents}

\DeclareUnicodeCharacter{00A0}{~}

\usepackage{algorithm}
\usepackage[noend]{algpseudocode}
\usepackage[english]{babel}
\newtheorem{theorem}{Theorem}

\newtheorem{proposition}{Proposition}
\newtheorem{lemma}[theorem]{Lemma}

\theoremstyle{remark}

\makeatletter
\def\BState{\State\hskip-\ALG@thistlm}
\makeatother
 
\usepackage{amsthm}
\usepackage{amssymb}
\usepackage[hidelinks]{hyperref}

\begin{document}

\begin{frontmatter}





\title{Expected Optimal Distances of Random Bipartite Matching in $D$-dimensional Spaces:
Approximate Formulas and Applications to Mobility Services}


\author[label1]{Shiyu Shen}
\author[label1]{Yuhui Zhai}
\author[label1]{Yanfeng Ouyang}

\address[label1]{Department of Civil and Environmental Engineering, University of Illinois at Urbana-Champaign, Urbana, IL 61801, USA}

\begin{abstract}
The bipartite matching problem has been at the core of many theoretical and practical challenges across various domains. 
Although many well-known algorithms can solve each bipartite matching problem instance efficiently, it remains an open question how one could estimate the expected optimal matching distance for arbitrary numbers of randomly distributed vertices in $D$-dimensional spaces (referred to as a random bipartite matching problem, or RBMP). This paper proposes a comprehensive modeling framework that yields closed-form approximate formulas for estimating the expected optimal matching cost across three interrelated but increasingly complex versions of RBMPs:
(i) RBMP-I, where edge costs are independently and identically distributed (i.i.d.); 
(ii) RBMP-S, where edge costs represent distances between vertices uniformly distributed on the surface of a hyper-sphere in a $D$-dimensional Euclidean space; and
(iii) RBMP-B, where the vertices are uniformly distributed in a hyper-ball within a $D$-dimensional L$^p$ metric space.
A series of Monte-Carlo simulation experiments are conducted to verify the accuracy of the proposed formulas under varying parameter combinations (e.g., numbers of bipartite vertices, spatial dimensions, and distance metrics). These proposed distance estimates could be key for strategic performance evaluation and resource planning in a wide variety of application contexts. As an illustration, we focus on on-demand mobility services (e.g., e-hailing taxi system) where matches must be made between customers and vehicles that are randomly distributed over time and space. We show how the proposed distance formulas provide a theoretical foundation for the empirically assumed Cobb-Douglas matching function in the field, and reveal conditions under which it can work well. Our formulas can also be easily incorporated into optimization models to select on-demand mobility operation strategies (e.g., whether newly arriving customers shall be instantly matched or pooled into a batch for matching). 
Agent-based simulations are conducted to verify the predicted performance of the demand pooling strategy for two types of e-hailing taxi systems. The results not only demonstrate the accuracy of the proposed model estimates under various service conditions, but also offer valuable managerial insights into the service operators' optimal strategies.

\end{abstract}




\end{frontmatter}

\section{Introduction}
\noindent 
\subsection{Background}
The bipartite matching problem is a fundamental problem in the field of applied mathematics and combinatorial optimization. 
Its variations have been applied to a variety of theoretical or practical contexts. 
In the field of physics, they can be used to capture important properties of various disordered complex systems, such as identifying the patterns and energy configurations of atomic magnets in spin glass systems \citep{mezard_replicas_1985}. 
In the field of biology, they can be used to describe interactions between species in an ecosystem \citep{simmons_motifs_2019}, or to analyze pairwise protein-protein interactions \citep{tanay_revealing_2004}. 
In the field of computer science, they are formulated to enhance graph-based pattern recognition by mapping the underlying data structures of images/signals to their features/labels \citep{yu_learning_2020}. In the field of social media and e-commerce, they are used to capture user/information interactions among distinct socioeconomic groups \citep{zhou_bipartite_2007, wu_graph_2022}.

In particular, bipartite matching problems are widely applicable to many transportation/mobility systems, 
where one fundamental challenge is to find matches between travel demand and resource supply over time (e.g., a planning horizon) and space (e.g., a service region).
For instance, they can be used to address the operations of multiple elevators in a tall building, where customers arrive randomly at different floors and are matched to one of the elevators --- All the vertices of interest would be distributed in a 1-dimensional space. In 2- or 3-dimensional spaces, they can be used to describe how surface courier vehicles and idle taxis, or freight drones and passenger aerial vehicles, are matched to their customers in a city. 
In particular, the latter examples have received a lot of attention, as ride-hailing services for passengers (e.g., as those offered by Transportation Network Companies, TNCs) and freight (e.g., as those offered by courier service companies) are booming all over the world.
In their typical operations, customers and vehicles evolve in the system as random vertices in a service region, and the service platform periodically (e.g., every a few seconds) makes vehicle dispatch and allocation decisions to best serve the customers. 
In each decision epoch, the system captures a snapshot of its current state by gathering information on both idle vehicles (e.g., locations) and new customers (e.g., origin and destination locations, and the elapsed waiting time).
A bipartite graph is constructed, where one subset of vertices include all idle vehicles, and the other subset all new customer origins. Costs of the edges could be based on distance (or travel cost, time), customer preference, and customer priority. Matches are then optimized by the platform based on a predefined objective, such as minimizing the total matching distances for pickups (between the vehicles and the customer origins).

This bipartite matching based mobility service scheme stands out for its ease of computation and implementation at the operational level. For each decision epoch, the associated problem instance can be solved quickly using linear programming methods or algorithms; e.g., the Hungarian algorithm \citep{kuhn_hungarian_1955}, Jonker–Volgenant algorithm \citep{jonker_shortest_1987}, and their variations can generate optimal solutions in polynomial time. In recent years, advanced machine-learning approaches, such as those reviewed in \citep{zhang_survey_2023}, can effectively tackle larger problem instances within a short time.

For strategic planning, however, the operators usually need to estimate the service efficiency under a large number of possible realizations of supply and demand scenarios (e.g., different vehicle and customer distributions), rather than finding exact vehicle-customer matches for one particular problem instance. 
The average matching distance between the vehicles and the customer origins, also commonly referred to as the average ``deadheading" distance, stands as a key indicator of service efficiency. It indicates the ``unproductive" efforts made by both customers (i.e., waiting for pickup) and vehicles (i.e., running empty) within the mobility system. 
Understanding the relationship between the average matching distance and vehicle/customer distribution can help improve service efficiency in many ways.
Operators, for example, often need to set operation standards such as the maximum time a passenger may have to wait for a vehicle assignment, or the time between consecutive decision epochs (i.e., pooling interval). A longer pooling interval may lead to more customers/vehicles appearing in one matching problem instance, which may potentially reduce the resulting matching distance. However, it also implies that the customers need to wait longer to find a match. Knowledge on the average matching distance directly helps find a balance between these conflicting objectives and identify the optimal operational standard. 
Moreover, operators often deploy new tactical-level strategies to further enhance their service efficiency, such as swapping vehicles and customers across already-matched pairs  
whenever new options become feasible \citep{ouyang_measurement_2023, shen_dynamic_2023}. Knowledge on the average matching distances before and after a random swap directly helps analyze the effectiveness and worthiness of these strategies. 

\subsection{Random Bipartite Matching Problem}
\begin{figure}[ht!]
    \centering
    \begin{subfigure}[b]{0.28\textwidth}
        \centering
        \includegraphics[height=1.8in]{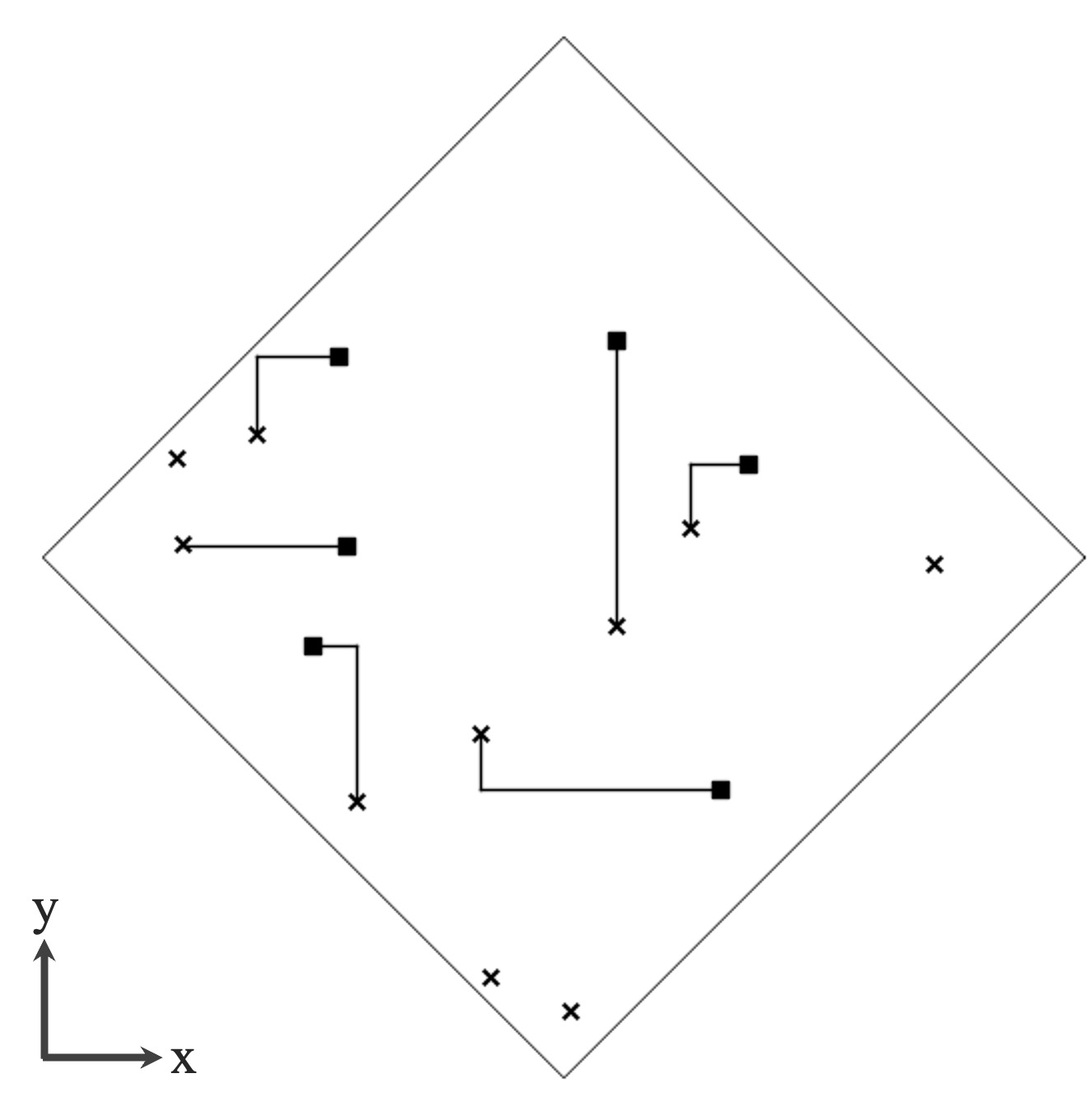}
        \caption{$D=2, p=1$.}
    \end{subfigure}%
    \hspace{0.5cm} 
    \begin{subfigure}[b]{0.28\textwidth}
        \centering
        \includegraphics[height=1.8in]{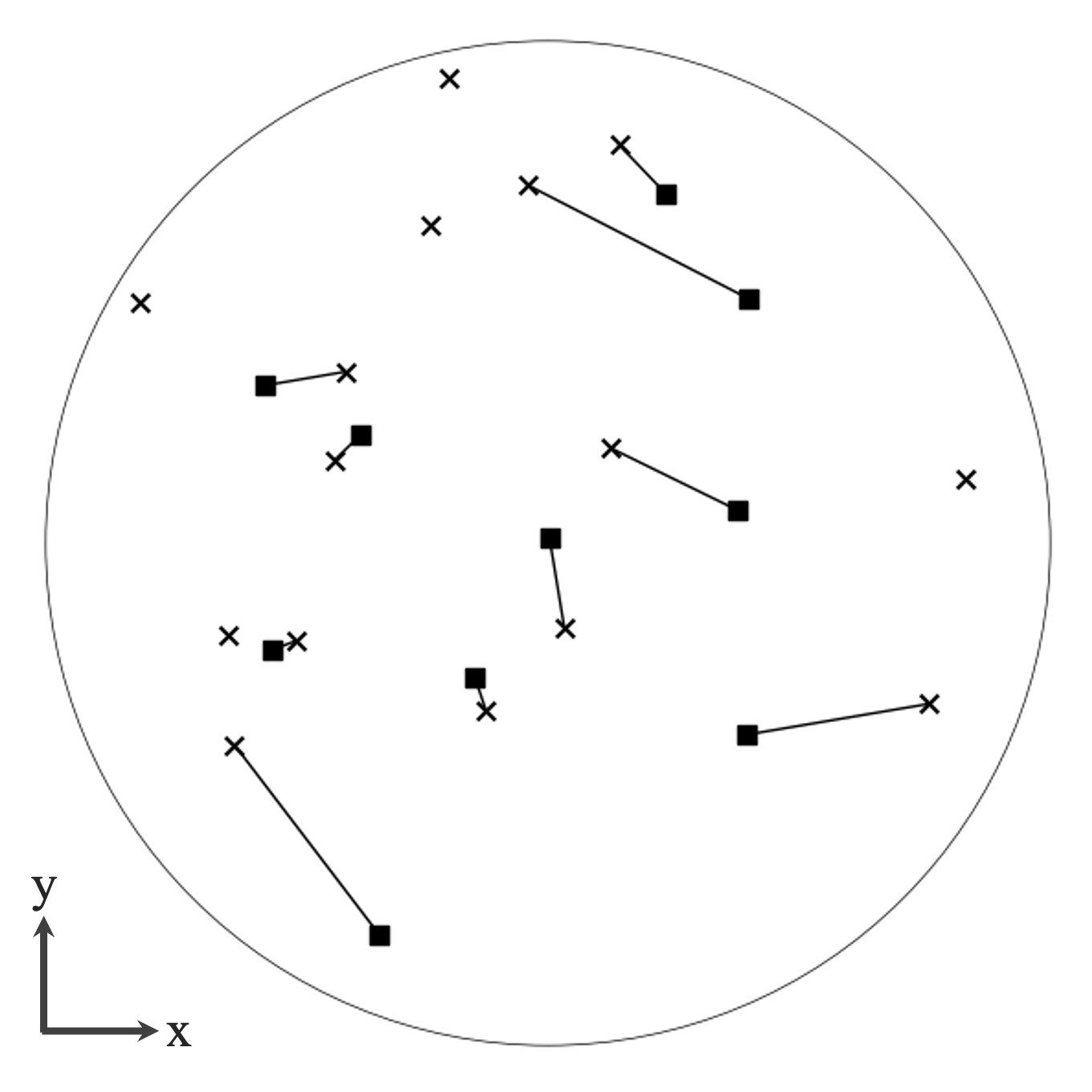}
        \caption{$D=2, p=2$.}
    \end{subfigure}
    \hspace{0.5cm} 
    \begin{subfigure}[b]{0.28\textwidth}
        \centering
        \includegraphics[height=1.8in]{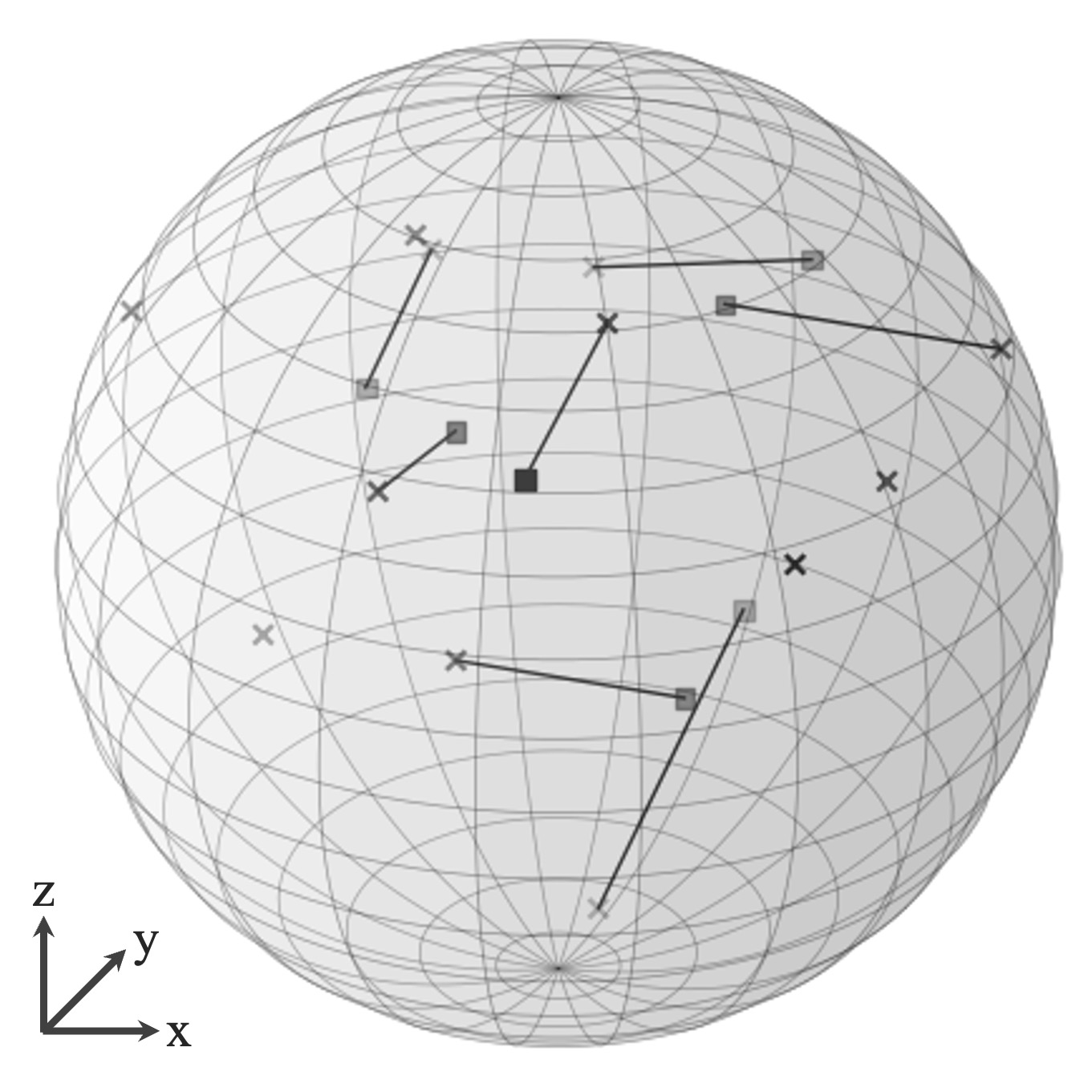}
        \caption{$D=3, p=2$.}
    \end{subfigure}
    \caption{Examples of RBMP instances in a $D$-dimensional L$^p$ metric space.}
    \label{fig:RBMP_examples}
\end{figure}

These needs give rise to a stochastic version of the bipartite matching problem where 
(i) both subsets of vertices are randomly distributed within a given domain in a $D$-dimensional L$^p$ metric space;
(ii) the weight (cost) on each edge is determined by the distance between the respective vertices; 
(iii) for each realization of random vertices (i.e., one instance), an optimal matching solution is found to minimize the average matching cost per matched pair. 
The ``Random Bipartite Matching Problem" (RBMP) is interested in estimating the average ``optimal matching cost" per vertex pair across all possible realizations. Figure \ref{fig:RBMP_examples} illustrates several representative examples of RBMP instances when $D=2, 3$ and $p=1, 2$. 

To the best of our knowledge, estimating the expected ``optimal matching distance" for these RBMPs remains a challenging task. 
While one could use the state-of-the-art techniques (as those employed by the TNCs) to solve a sufficiently large number of (simulated or historical) problem instances and conduct statistical analysis, analytical models are favored for two reasons. First, they can provide more insights, such as those developed in \citet{daganzo_analysis_2020} and \citet{ouyang_performance_2021}, by directly connecting the key performance metrics (e.g., the expected vehicle distance traveled) with mobility service strategies. 
Moreover, this type of analytical models can be easily incorporated into more comprehensive optimization/equilibrium models, helping the operators or regulators in optimizing their service offerings to achieve higher service efficiency or social welfare \citep{zha_economic_2016, ouyang_performance_2021, liu_mobility_2021, liu_planning_2023}.
However, as we will discuss in the literature review, existing analytical models and formulas for estimating matching distances are limited only to special scenarios such as: (i) the densities of vertex distributions are (nearly) equal; or (ii) the number of dimensions is limited to only one (i.e., $D=1$); or (iii) the distance is measured by only Euclidean or squared Euclidean metric in an unbounded domain. 
More importantly, many of these studies primarily focus on identifying asymptotic approximations where the number of vertices approaches infinity.

\begin{figure}[ht!]
    \centering
    \includegraphics[height=2in]{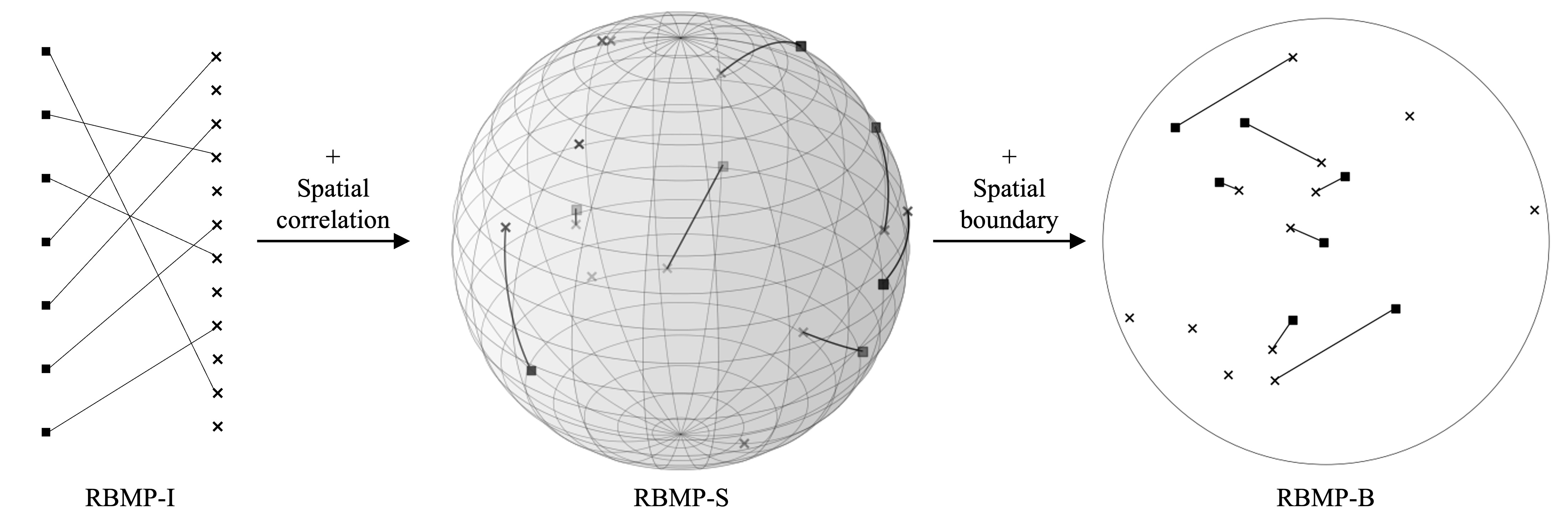}
    \caption{Three versions of RBMPs studied in this paper.}
    \label{fig:RBMP_model}
\end{figure}

In response to these challenges, this paper proposes a comprehensive modeling framework that yields closed-form approximate formulas for the expected optimal matching cost across three interrelated but increasingly complex versions of RBMP, as illustrated in Figure \ref{fig:RBMP_model}.
We first develop an analytical model to estimate the distribution and moments of the optimal matching cost for RBMP-I, where edge costs are independently and identically distributed (i.i.d.).
The model is then extended to RBMP-S, in which edge costs represent distances between vertices uniformly distributed over an unbounded domain in a $D$-dimensional space; e.g., the surface of a sphere.
In so doing, a correction method is proposed to capture correlations among vertex distributions and matched pairs in the spatial domain. 
Finally, we extend the model to RBMP-B when the vertices are distributed in a bounded domain of the $D$-dimensional space; e.g., a circular area on a plane. 
Additional correction must be added to account for the impacts of domain boundary on vertex distributions and matching likelihoods.

\subsection{Contribution}

The proposed models, for the first time, identify a desirable convergence property in RBMP-I when the edge costs follow a class of power-law distributions.
It reveals why the greedy matching cost rapidly approaches the optimal cost as one subset of vertices becomes dominant, and leads to closed-form approximate formulas to estimate the moments of the optimal matching cost in RBMP-I. 
Building on these results, we further obtain the first set of closed-form approximate formulas for the expected matching distances in RBMP-S and RBMP-B, in consideration of spatial correlations among the vertex distributions. 
They are applicable to arbitrary vertex set sizes, spatial dimensions, and a broad class of distance metrics. 
The proposed distance formulas also provide a theoretical explanation of the conditions under which existing RBMP-B models (e.g., a widely-used empirical Cobb-Douglas meeting function from the field of transportation science) is applicable, and how some of its parameter values should be set.

Beyond the theoretical significance, the proposed distance estimate formulas are practically important because they can be efficiently used for strategic performance evaluation and resource planning of many real-world problems.
For example, by incorporating the proposed formulas into an optimization framework, we reveal the effectiveness of demand pooling strategies for ride-hailing services across a range of service conditions (e.g., varying demand rates and fleet sizes). 
The findings show that instant customer matching may already be the optimal strategy in a system with a sufficiently large fleet, while demand pooling can be beneficial only when the system is stuck in an inefficient equilibrium state. In addition, the models can be extended to provide managerial insights into various other operational strategies, such as fleet deployment, pricing, vehicle repositioning, etc.

The accuracy of the proposed distance formulas is verified by a set of Monte-Carlo simulations conducted over a wide range of problem settings. The results indicate that our proposed formulas can provide highly accurate distance estimations. 
In addition, a series of agent-based simulations are conducted to verify how our formulas can effectively help enhance the design of optimal mobility service strategies. The results also show that the model estimates match quite well with the simulation measurements under all considered settings. 

The remainder of this paper is organized as follows. Section \ref{sec: literature} provides a review of the related literature. Section \ref{sec:RBMP-I} presents models and formulas for RBMP-I as a building block. Sections \ref{sec:RBMP-S} and \ref{sec:RBMP-B} present models and formulas for RBMP-S and RBMP-B, respectively. 
Section \ref{sec: experiment} presents numerical experiments. 
Section \ref{sec: application} shows how the proposed formulas can be applied, as an example, to improve on-demand mobility systems. 
Finally, Section \ref{sec: conclusion} offers concluding remarks and suggestions for future research. 

\section{Literature Review} \label{sec: literature}

Statistical physicists and mathematicians seem to be among the first to estimate the ``average optimal cost" of RBMPs.
\citet{mezard_replicas_1985} used the ``replica method'' from the mean field theory to derive analytical formulas for the expected asymptotic optimal cost in a balanced RBMP-I, where the two vertex subsets are of equal and sufficiently large size, and the edge costs are i.i.d. from a uniform distribution. 
The accuracy of their proposed formulas was verified by extensive numerical experiments in \citet{brunetti_extensive_1991}, where the errors were shown to be less than 1\% from Monte-Carlo simulation results. 
In the following years, formulas for RBMP-I variations with exponentially distributed and i.i.d. edge costs were conjectured; e.g., for complete matching of balanced bipartite graphs \citep{parisi_conjecture_1998}, or partial matching of unbalanced bipartite graphs \citep{coppersmith_constructive_1998, alm_exact_2002}. 
Based on the memory-less property of exponential distribution, these conjectures were later rigorously proven by \citet{linusson_proof_2004} and \citet{nair_proofs_2005}.

These earlier studies paved the way for finding the average optimal distance for balanced RBMP-S in a Euclidean space. 
When edge costs are measured in a metric space, correlations exist among these costs due to implied spatial properties (such as triangular inequalities) among vertex distributions. Hence, RBMP-S is more difficult to analyze than RBMP-I with simple i.i.d. cost distributions. 
\cite{mezard_euclidean_1988} formally defined the Euclidean matching problems by connecting  RBMP-I, where the edge costs follow a power-law distribution, to RBMP-S. They proposed asymptotic approximations to estimate the average optimal costs when the number of vertices approaches infinity, in consideration of first-order triangular correlations. 
Subsequent works extended these asymptotic analyses to problems where the edge cost is measured by the distance between vertices raised to a power.
For example, \citet{caracciolo_scaling_2014, caracciolo_quadratic_2015, caracciolo_scaling_2015} proposed formulas for problems measured by squared Euclidean distance, also known as the ``quadratic" cost. 
However, since all these asymptotic approximations were derived under the strong assumption that the number of bipartite vertices approaches infinity, they could only serve as bounds rather than exact estimations for finite-size problems. 
The only exception is the 1-dimensional problem, for which \citet{caracciolo_one-dimensional_2014, boniolo_correlation_2014, caracciolo_selberg_2019} proposed exact formulas when the cost exponent is greater than one.
These powered-cost problems have received particular attention because the optimal matching has a special property: after sorting both vertex sets, the $i$-th vertex from one set must be matched with the $i$-th vertex from the other set. 

When the edge cost is measured by the exact distance (as in RBMP-S), the sorted matching property no longer holds, making the problem much more challenging. 
Most existing results focus on analyzing the asymptotic scaling behavior of these problems. 
For example, based on results from the quadratic-cost problem, \citet{caracciolo_scaling_2014} formulated several simple parametric models under the key assumption that the ``leading terms'' of the exact optimal cost for finite-size balanced RBMP-S could be approximated by the product of the asymptotic estimate and a scaling coefficient. However, this scaling behavior is only conjectured to apply to the exact-cost problems when $D\ge 3$, and the coefficients in the formulas must be estimated from statistical regression on simulated data. 
In addition, \citet{kanoria_dynamic_2022} showed that the ``scaling behavior" from these static problems can be leveraged, via a hierarchical greedy algorithm, to estimate the scale of the lower and upper bounds on the minimum achievable expected distance in an unbalanced and dynamic problem. However, they did not provide any explicit formulas for exact distance estimation. 
More recently, \citet{zhai_average_2024} addressed the finite-size balanced RBMP-S when $D=1$, by relating the total optimal matching distance to the area underneath a random walk. However, this special property can not be generalized to higher-dimensional problems.

In the field of transportation science, RBMP-B is of particular interest, as real-world transportation services typically operate within bounded spatial domains. 
However, these problems are even more challenging to analyze than RBMP-S, as the presence of spatial boundaries introduces stronger correlations among edge costs. 
Earlier mathematical studies have identified the asymptotic scaling factors for balanced RBMP-B when $D = 2$ \citep{ajtai_optimal_1984} and $D \ge 3$ \citep{talagrand_matching_1992}. 
In addition, the statistical physics studies on 1-dimensional powered-cost RBMP-S also extended their models to RBMP-B \citep{caracciolo_one-dimensional_2014, boniolo_correlation_2014, caracciolo_selberg_2019}. 
Meanwhile, the transportation literature has primarily focused on the exact-cost RBMP-B when $D = 1$ and $D = 2$.
For instance, \citet{daganzo_bounds_2004} explored a related random network transportation problem, which estimates the expected shipment costs needed to clear normally distributed supply and demand quantities at uniformly distributed vertices. 
An exact formula was proposed for balanced 1-dimensional cases, but for higher-dimensional problems, the formula involves parameters that are estimated through statistical regression (similar to those found in the field of statistical physics). 
Moreover, their results cannot be directly applied to RBMP-B, where the demand or supply quantities are either positive or negative one. 
Very recently, \citet{zhai_average_2024} proposed distance formulas for RBMP-B only when $D=1$ and extended the result to a regular network. 

Other studies in the transportation field specifically explored RBMP-B in a 2-dimensional space. 
The most basic model in \citet{daganzo_approximate_1978} estimates the probability distribution and expectation of 
a special case of RBMP-B where the number of vertices in one subset is exactly one. 
However, the model ignores boundary effects and only holds when the number of vertices in the other subset is sufficiently large. 
For more general matching problems with arbitrary numbers of vertices, a widely used model by \citet{yang_equilibria_2010} --- called the ``matching function'' --- assumes that the rate of successful matches between taxi customers and vehicles (as two subsets of vertices) takes the form of a Cobb-Douglas formula. Little's formula can then be used to derive the expected matching distance in a power function form. This function has been used extensively to analyze and plan mobility services  \citep{yang_equilibrium_2011, zha_economic_2016, zhang_efficiency_2019}. 
The parameters in the Cobb-Douglas formula are typically ``guessed" or estimated via statistical regression on simulated or real-world data. 
Section \ref{sec: cobb} will provide more detailed discussion on this model.
More recently, \citet{lei_average_2024} developed an analytical model to estimate the expected tour length for visiting a random subset of uniformly distributed vertices within a compact region.
Although this study did not specifically address RBMP-B, it offered insights in identifying the expected distance between random subsets of vertices. 

In summary, to the best of our knowledge, none of the existing studies have derived formulas that can estimate the expected optimal matching distance for general unbalanced or balanced RBMP-S and RBMP-B in arbitrary dimensions and arbitrary distance metrics, which are of our primary interest. 
As summarized in Table \ref{tab:literature}, existing results are limited to special cases, 
and some of them rely heavily on statistical curve-fitting to estimate coefficients. 
As one step towards filling such a gap, this paper aims to develop a general modeling framework built upon cost estimators for finite-size RBMP-I under a class of power-law distributions, and extend it to accommodate RBMP-S and RBMP-B with arbitrary vertex sizes, arbitrary number of spatial dimensions, and a broader class of distance metrics (e.g., great-circle and L$^p$).
Our goal is to derive practical formulas, where the leading terms are supported by sufficient theoretical justification, that can directly be used to support engineering and operational decisions.

\begin{table}[ht!]
\centering
\caption{A summary of the some relevant literature.}
\label{tab:literature}
\resizebox{\columnwidth}{!}{%
\begin{tabular}{|c|cc|ccc|}
\hline
Problem   Type &
  \multicolumn{2}{c|}{Cost   distribution} &
  \multicolumn{1}{c|}{Number of   vertices} &
  \multicolumn{1}{c|}{Formula type} &
  Literature \\ \hline
\multirow{4}{*}{RBMP-I} 
 &
  \multicolumn{2}{c|}{Uniform} &
  \multicolumn{1}{c|}{Balanced} &
  \multicolumn{1}{c|}{Asymptotic} &
  \citet{mezard_replicas_1985} \\ \cline{2-6} 
 &
  \multicolumn{2}{c|}{Power-law} &
  \multicolumn{1}{c|}{Balanced} &
  \multicolumn{1}{c|}{Asymptotic} &
  \citet{mezard_euclidean_1988} \\ \cline{2-6} 
 &
  \multicolumn{2}{c|}{Exponential} &
  \multicolumn{1}{c|}{All} &
  \multicolumn{1}{c|}{Exact} &
  \citet{linusson_proof_2004,nair_proofs_2005} \\ \cline{2-6} 
 &
  \multicolumn{2}{c|}{\textbf{Power-law}} &
  \multicolumn{1}{c|}{\textbf{All}} &
  \multicolumn{1}{c|}{\textbf{Approximate}} &
  \textbf{This study} \\ \hline
 &
  \multicolumn{1}{c|}{Spatial dimension} &
  Distance metric &
  \multicolumn{3}{c|}{} \\ \hline
\multirow{4}{*}{RBMP-S} 
 &
  \multicolumn{1}{c|}{All} &
  Euclidean &
  \multicolumn{1}{c|}{Balanced} &
  \multicolumn{1}{c|}{Asymptotic} &
  \citet{mezard_euclidean_1988} \\ \cline{2-6} 
 &
  \multicolumn{1}{c|}{$\ge 3$} &
  Euclidean &
  \multicolumn{1}{c|}{Balanced} &
  \multicolumn{1}{c|}{Parametric} &
  \citet{caracciolo_scaling_2014} \\ \cline{2-6} 
 &
  \multicolumn{1}{c|}{1} &
  Euclidean &
  \multicolumn{1}{c|}{Balanced} &
  \multicolumn{1}{c|}{Approximate} &
  \citet{zhai_average_2024} \\ \cline{2-6} 
 &
  \multicolumn{1}{c|}{\textbf{All}} &
  \textbf{Great-circle} &
  \multicolumn{1}{c|}{\textbf{All}} &
  \multicolumn{1}{c|}{\textbf{Approximate}} &
  \textbf{This study} \\ \hline
\multirow{6}{*}{RBMP-B} 
 &
  \multicolumn{1}{c|}{1} &
  Euclidean &
  \multicolumn{1}{c|}{Balanced} &
  \multicolumn{1}{c|}{Approximate} &
  \citet{daganzo_bounds_2004} \\ \cline{2-6} 
 &
  \multicolumn{1}{c|}{1} &
  Euclidean &
  \multicolumn{1}{c|}{All} &
  \multicolumn{1}{c|}{Approximate} &
  \citet{zhai_average_2024} \\ \cline{2-6} 
 &
  \multicolumn{1}{c|}{2} &
  Euclidean &
  \multicolumn{1}{c|}{Balanced} &
  \multicolumn{1}{c|}{Parametric} &
  \citet{daganzo_bounds_2004} \\ \cline{2-6} 
 &
  \multicolumn{1}{c|}{2} &
  Euclidean &
  \multicolumn{1}{c|}{One} &
  \multicolumn{1}{c|}{Approximate} &
  \citet{daganzo_approximate_1978} \\ \cline{2-6} 
 &
  \multicolumn{1}{c|}{2} &
  All &
  \multicolumn{1}{c|}{All} &
  \multicolumn{1}{c|}{Empirical} &
  \citet{yang_equilibria_2010} \\ \cline{2-6} 
 &
  \multicolumn{1}{c|}{\textbf{All}} &
  \textbf{L$^p$} &
  \multicolumn{1}{c|}{\textbf{All}} &
  \multicolumn{1}{c|}{\textbf{Approximate}} &
  \textbf{This study} \\ \hline
\end{tabular}%
}
\end{table}

\section{RBMP-I}
\label{sec:RBMP-I}
This section develops approximate formulas for estimating the optimal matching cost in RBMP-I. The conditional optimal matching cost is derived explicitly using order statistics, whereas the optimal matching probability is approximated via a two-step heuristic process. 
These derivations lead to two main results: one for RBMP-I under an arbitrary cost distribution, and the other for a special case where the costs follow a power-law distribution.
The results lay a foundation for deriving approximate formulas for RBMP-S and RBMP-B.
\subsection{Problem Definition}
An RBMP-I is formally defined as follows. 
Consider a bipartite graph $G = (U \cup V; E)$, where $U = \{1,\ldots, m\}$ and $V = \{1,\ldots,n\}$ represent the two disjoint subsets of vertices with given cardinalities $m, n \in \mathbb{Z}^+$, respectively.
Without loss of generality, we assume $m \le n$, and refer to the vertices in $U$ and $V$ as the ``demand" and ``supply" vertices, respectively.
The edges in set $E$ connect every pair of vertices across $u \in U$ and $v \in V$, and the associated costs are drawn independently and identically from a given probability distribution with cumulative distribution function (CDF) $F_{C}(\cdot)$.
For each realized problem instance, the cost of edge $(u, v) \in E$ is denoted as $c_{u,v}$. 
Since $m \le n$, each demand vertex $u \in U$ will find exactly one match $v(u) \in V$, and a set of optimal matches, denoted by $\{\accentset{\ast}{v}(u), \forall u\in U\}$, can be found for this instance to minimize the total cost between the matched vertex pairs: 
\begin{align}
\min_{\{v(u), \forall u\in U\}} \quad &\sum_{u \in U} c_{u,v(u)}, \label{eq:RBMP_obj}\\
\text{s.t.} \quad 
& v(u) \ne \emptyset, \forall u\in U, \label{eq:RBMP_con1}\\
& v(u) \ne v(u'), \forall u \ne u'\in U. \label{eq:RBMP_con2}
\end{align}
RBMP-I seeks the expected optimal matching distance across all possible realizations of matching problem instances, which shall be 
dependent on only the vertex set sizes $m$ and $n$, and the edge cost distribution $F_{C}(\cdot)$. 
This section derives approximate formulas for the distribution and moments of the optimal matching cost of a randomly sampled $u\in U$ from an instance to its optimal match $\accentset{\ast}{v}(u)$, denoted by $\accentset{\ast}{C}(m,n)$. For notation simplicity, we may drop the arguments and simply use $\accentset{\ast}{C}$ (unless necessary otherwise).

In an instance, the optimal match of a realized $u \in U$ 
may be a supply vertex that has the $k$-th smallest realized costs, where $k \in \{1,2,\ldots,n\}$. 
For convenience, we refer to this supply vertex as the $k$-th nearest ``neighbor" of $u$. 
Estimating the distribution of $\accentset{\ast}{C}$ involves two key steps:
\begin{itemize}
    \item[(i)] Estimate the distribution of the cost of any randomly sampled $u \in U$ to its $k$-th nearest neighbor, denoted by the conditional random variable $C \mid k$.
    \item[(ii)] Determine the probability that any randomly sampled $u \in U$ is matched to its $k$-th nearest neighbor in the optimal solution, denoted by $\accentset{\ast}{\mathbb{P}}(k)$.
\end{itemize}

We begin by examining $C \mid k$ in step (i), whose distribution can be easily derived using the theory of order statistics. By definition, $C \mid k$ corresponds to the $k$-th order statistic among $n$ i.i.d. samples drawn from distribution $F_{C}(\cdot)$. It is well known \citep{arnold_first_2008} that the CDF and the $M$-th moment of $C \mid k$, denoted by $F_{C \mid k}(x)$ and $\mathbb{E}[C^M \mid k]$, respectively, are given by:
\begin{align}
    F_{C \mid k}(x)
    &= I_{F_{C}(x)}(k,n-k+1), \label{eq:F_X_k}
    \\
    \mathbb{E}[C^M \mid k]
    &= \int_{-\infty}^{+\infty}
    \frac{x^M}{\text{B}(k,n-k+1)}\left[1-F_{C}(x)\right]^{n-k}\left[F_{C}(x)\right]^{k-1}\text{ d}F_{C}(x). 
    \label{eq:E_X_k}
\end{align}
Here $I_{z}(a,b) = \frac{\text{B}(z;a,b)}{\text{B}(a,b)}$ is the regularized beta function, $\text{B}(z;a,b) = \int_0^z t^{a-1}(1-t)^{b-1} \text{ d}t$ is the incomplete beta function, $\text{B}(a,b)=\int_0^1 t^{a-1}(1-t)^{b-1} \text{ d}t = \frac{\Gamma(a)\Gamma(b)}{\Gamma(a+b)}$ is the beta function, and $\Gamma(z) = \int_0^{\infty}t^{z-1}e^{-t} \text{ d}t$ is the gamma function.

Since the objective of the matching problem in Equation~\eqref{eq:RBMP_obj} is to minimize the total cost across all matches for $u \in U$, the worst-case match for any $u$ is its $m$-th nearest neighbor. Therefore, $k$ can only take values from $\{1,\ldots,m\}$. The CDF and the $M$-th moment of $\accentset{\ast}{C}$ can then be derived using the law of total probability, as stated in the following proposition.
\begin{proposition} \label{prop:X_optimal}
For an RBMP-I with given $n, m \in \mathbb{Z}^+$, and $F_{C}(\cdot)$, the CDF and 
the $M$-th moment of the optimal matching cost $\accentset{\ast}{C}$ are respectively given by:
\begin{align}
    F_{\accentset{\ast}{C}}(x) 
    &= 
    \sum_{k=1}^m I_{F_{C}(x)}(k,n-k+1)\cdot \accentset{\ast}{\mathbb{P}}(k)
    ,\\
    \mathbb{E}[\accentset{\ast}{C}^M] 
    &= 
    \sum_{k=1}^m  \mathbb{E}[C^M \mid k] \cdot \accentset{\ast}{\mathbb{P}}(k),
\end{align}
where $\mathbb{E}[C^M \mid k]$ is given by Equation \eqref{eq:E_X_k}.
\end{proposition}

The more challenging task is to estimate the optimal matching probabilities $\{\accentset{\ast}{\mathbb{P}}(k) \mid 1\le k\le m \}$ in step (ii). The next section presents an approximation to these probabilities based on a two-step heuristic process.

\subsection{Approximate Matching Probabilities}
\label{sec:P_k}
The idea of the two-step heuristic process is as follows. 
In the initial step, a greedy matching process is used to sequentially match each demand vertex $u \in U$ to a supply vertex.
The probability for $u$ to be matched to its $k$-th nearest neighbor in this greedy solution is denoted by $\bar{\mathbb{P}}(k)$. 
From a probabilistic point of view, the greedy solution provides an upper bound on the matching cost, though each demand vertex may still have the opportunity to be matched to a closer neighbor in the better solution. 
This potential improvement is approximated by a local search process, where we consider the probability of swapping 
a demand vertex $u$'s initial match (say, its $k'$-th nearest neighbor) in the greedy matching solution, to a closer match  
(say, its $k$-th nearest neighbor, where $ k < k'$) in a refined solution, with the corresponding matching probability $\hat{\mathbb{P}}(k)$. 

The resulting matching probabilities $\bar{\mathbb{P}}(k)$ and $\hat{\mathbb{P}}(k)$ can be directly applied, as in Proposition \ref{prop:X_optimal}, to estimate the expected near-optimum matching costs. 
In addition, as discussed later in Section \ref{sec:power-law}, when $F_C(\cdot)$ has a power-law form, we identify a desirable  convergence property: the expected greedy matching cost rapidly converges to the expected optimal matching cost when one subset of vertices becomes dominant (i.e., $n\gg m$). This suggests that, in such a special case, the greedy matching probability $\bar{\mathbb{P}}(k)$ can be used as an effective approximation for the optimal matching probability $\accentset{\ast}{\mathbb{P}}(k)$. For more general values of $m$ and $n$, the refined matching probability $\hat{\mathbb{P}}(k)$ provides a closer approximation to $\accentset{\ast}{\mathbb{P}}(k)$. We next describe the processes for deriving $\bar{\mathbb{P}}(k)$ and $\hat{\mathbb{P}}(k)$. 

\subsubsection{Greedy matching process.}
\label{sec:greedy}
For any given RBMP-I instance, consider the following greedy matching process: 
\begin{itemize}
    \item[(i)] For each vertex $u\in U$, sort all the vertices in $V$ into a sequence based on their distance to $u$ in an ascending order, $K_u[v]: V \rightarrow \{1, \cdots, n\}$.
    \item[(ii)] Randomly select a demand vertex $u \in U$ and match it the nearest available supply vertex that has not been previously matched with another demand vertex, and denote this match by $\bar{v}(u)$.
    \item[(iii)] Repeat step (ii) with another unmatched vertex in $U$, until all vertices in $U$ have been matched.
\end{itemize}
The following lemma presents the probabilities $\{\bar{\mathbb{P}}(k) \mid 1\le k\le m \}$ derived from this greedy process.

\begin{lemma}
\label{lemma:P_k_greedy}
For an RBMP-I with given $n, m \in \mathbb{Z}^+$, the greedy matching probability $\bar{\mathbb{P}}(k)$ is given by: 
\begin{align} \label{eq:P_k_greedy}
    \bar{\mathbb{P}}(k) 
    = \frac{1}{m} \cdot \sum_{i=k}^m \frac{\binom{n-k}{i-k}}{\binom{n}{i-1}},
    \quad k \in \{1, \cdots, m\}.
\end{align}
\end{lemma} 
\noindent{\bf Proof}. 
Since all demand vertices in $U$ are selected randomly, each of them has an equal probability $\frac{1}{m}$ of being the $i$-th vertex in the sorted sequence. 
Let $\bar{\mathbb{P}}(k \mid i)$ denote the conditional probability that the $i$-th vertex is matched to its $k$-th nearest neighbor, where $1\le k\le i\le m$. By the law of total probability, 
\begin{align} \label{eq:P_k_greedy_proof}
\bar{\mathbb{P}}(k) = \sum_{i=k}^m \frac{1}{m} \bar{\mathbb{P}}(k \mid i),
\quad k \in \{1, \cdots, m\}.
\end{align}

We now derive $\bar{\mathbb{P}}(k \mid i)$. From the perspective of the $i$-th demand vertex, it will be matched to its $k$-th nearest neighbor if and only if all its $k - 1$ nearest neighbors have already been matched to the $i - 1$ prior candidates while its $k$-th nearest neighbor is still available.
This can be regarded as a sampling process without replacement from a finite population of $n$ elements (i.e., all vertices in $V$), where $i - 1$ of them are labeled as ``successes" (already matched), and the rest as ``failures" (still available). 
From the negative hyper-geometric distribution, denoted by $\text{NHG}_{n, i-1, 1}(k-1)$, we have:
\begin{align} \label{eq:P_k_i}
\bar{\mathbb{P}}(k \mid i)
= \frac{\binom{i - 1}{k - 1}}{\binom{n}{k - 1}} \cdot \frac{n - (i - 1)}{n - (k - 1)}
= \frac{\binom{n - k}{i - k}}{\binom{n}{i - 1}}, \quad k \in \{1, \cdots, i\}.
\end{align}
Here the term $\frac{\binom{i - 1}{k - 1}}{\binom{n}{k - 1}}$ represents the probability that the first $k - 1$ trials are all successful. 
The term $\frac{n - (i - 1)}{n - (k - 1)}$ represents the probability that the first failure occurs at the $k$-th trial. 
Substituting Equation \eqref{eq:P_k_i} into Equation \eqref{eq:P_k_greedy_proof} leads to Equation \eqref{eq:P_k_greedy}.
\hfill \qedsymbol\\

In addition, we propose an approximation for the greedy matching probability $\bar{\mathbb{P}}(k)$, which serves as a key component in deriving approximate formulas for the expected greedy matching cost in Section \ref{sec:power-law}.
In deriving $\bar{\mathbb{P}}(k \mid i)$ in the above proof, we interpret the event that the $k-1$ nearest neighbors have already been matched to the $i-1$ prior candidates as drawing $k-1$ successes without replacement from a finite population of $n$ elements. 
When $n\gg 1$, \citet{ahlgren_probability_2014} has shown that the negative hyper-geometric distribution (sampling without replacement) with exactly one failure, i.e., $\text{NHG}_{n, \cdot, 1}(\cdot)$, converges to a geometric distribution (sampling with replacement). 
Intuitively, this indicates the difference between sampling with and without replacement in this setting becomes negligible, and the probability of having a success remains approximately the same (i.e., $\frac{i-1}{n}$) across all $k-1$ trials. In addition, when $k < i$, the probability that the $k$-th nearest neighbor is still available is approximately $1 - \frac{i-1}{n}$. While when $k = i$, the greedy matching rule guarantees that the $k$-th nearest neighbor is always available with probability 1.
Accordingly, $\bar{\mathbb{P}}(k \mid i)$ in Equation \eqref{eq:P_k_i} can be approximated by the following.
\begin{align}
    \bar{\mathbb{P}}(k\mid i) \xrightarrow{n\gg 1}
    \begin{cases}
            \left(\frac{i-1}{n}\right)^{k-1}\left(1-\frac{i-1}{n}\right), 
            k\in \{1,\cdots,i-1\}, \\
            \left(\frac{i-1}{n}\right)^{i-1} 
            ,  k = i.
    \end{cases}
\end{align}
Substituting this into Equation \eqref{eq:P_k_greedy_proof} yields the following approximation of $\bar{\mathbb{P}}(k)$:
\begin{align} \label{eq:P_k_geometric}
\begin{split}
    \bar{\mathbb{P}}(k) \xrightarrow{n\gg 1} \frac{1}{m} \left(\frac{k-1}{n}\right)^{k-1} + \frac{1}{m} \sum_{i=k+1}^{m}  \left(\frac{i-1}{n}\right)^{k-1}\left(1-\frac{i-1}{n}\right),\quad k \in \{1, \cdots, m\}.
\end{split}
\end{align}

\subsubsection{Vertex-swapping refinement process.}
\begin{figure}[ht!]

\begin{subfigure}{\textwidth}
    \centering
    \includegraphics[width=0.8\textwidth]{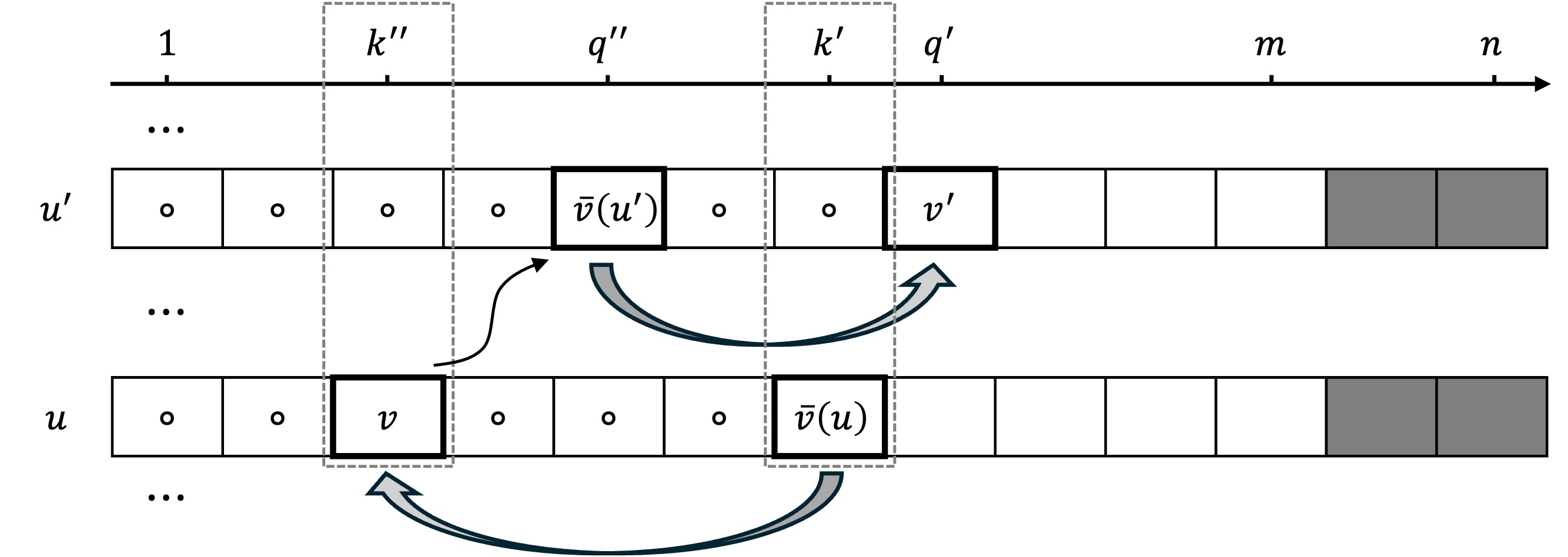}
    \caption{A potential swap between two demand vertices $u$ and $u'$.}
    \label{fig:swap_one}
\end{subfigure}

\vspace{1cm}

\begin{subfigure}{\textwidth}
    \centering
    \includegraphics[width=0.8\textwidth]{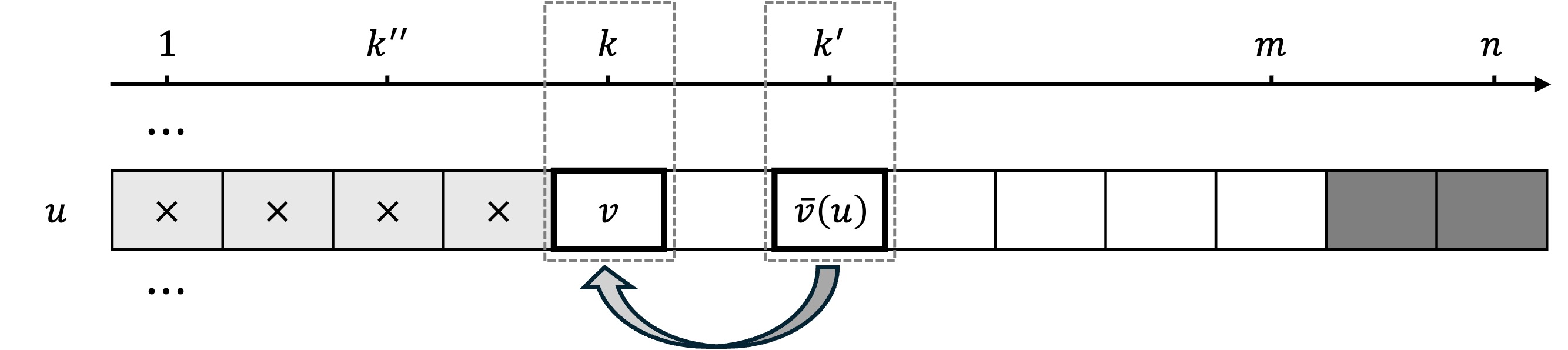}
    \caption{Vertex $u$ successfully swapped to match with its $k$-th nearest neighbor.}
    \label{fig:swap_success}
\end{subfigure}

\caption{Examples of swaps.}
\end{figure}

Next, we focus on a specific vertex $u \in U$, which is the $i$-th in the sorted sequence from the greedy matching process and is initially matched to its $k'$-th nearest neighbor; i.e., $K_u[\bar{v}(u)] = k'$. 
A successful vertex-swapping refinement for $u$ satisfies the rule below. 
\begin{itemize}
    \item[(i)] Select the supply vertices that are closer to $u$ than its current match $\bar{v}(u)$: $\hat{V}_u=\{v \in V \mid K_u[v] = 1, 2, \ldots, K_u[\bar{v}(u)]-1\}$.
    \item[(ii)] Examine each $v \in \hat{V}_u$ sequentially in ascending order of distance to see whether it is a ``feasible" swapping candidate. 
    Here, a swap from $\bar{v}(u)$ to $v$ 
    is considered feasible if it reduces total matching cost for $u$ to take the current match of another demand vertex $u'$, and let $u'$ take a new match vertex, while not affecting the matching of other demand vertices.\footnote{ 
Generally, there are many types of swaps. A one-vertex swap is feasible when a demand vertex $u$ swaps its current match to another currently ``available" and nearer supply vertex, while not affecting the matching of any other demand vertices. However, from the greedy matching process, all supply vertices in $\hat{V}_u$ must have already been matched to a previously selected demand vertex, such that no one-vertex swap is feasible in our context. More complex swaps involving more than two demand vertices are also possible, but they tend to yield diminishing marginal improvements, so we will not consider them in this paper.
}
Figure \ref{fig:swap_one} illustrates an example of a potential two-vertex swap between $u$ and $u'$, where their sorted supply vertices are represented in the vector tables. The supply vertices involved in the swap, $\bar{v}(u), v$ for $u$ and $\bar{v}(u'),v'$ for $u'$, are highlighted with bold black boxes. 
Specifically, the feasibility conditions are: 
\begin{itemize}
    \item[(a)] $v$ must already be matched to $u'$ where $\bar{v}(u')$ can only take a value from $1$ to $i-1$; i.e.,
$$1 \le K_{u'}[\bar{v}(u')] \le i - 1;$$ 
    \item[(b)] $v'$ must be farther away from $u'$ than $\bar{v}(u')$, or otherwise, additional demand vertices need to be involved; i.e., $$K_{u'}[\bar{v}(u')] < K_{u'}[v'] \le m;$$
    \item[(c)] $v'$ must either be matched to $u$ or remain unmatched, 
and all supply vertices between $\bar{v}(u')$ and $v'$ must have already been matched to other demand vertices; otherwise, a better alternative than $v'$ would exist. 
    \item[(d)] The cost reduction achieved by $u$ from swapping $\bar{v}(u)$ to $v$ must be no smaller than the cost increase incurred by $u'$ for swapping from $\bar{v}(u')$ to $v'$. 
\end{itemize}
    \item[(iii)] The first feasible $v$ found in step (ii) is selected, and $u$ swaps its current match $\bar{v}(u)$ to $v$. Figure \ref{fig:swap_success} illustrates an example for $u$ to successfully swaps its initial match $\bar{v}(u)$ to $v$, with all supply vertices closer than $v$ being infeasible. 
    If no feasible swap is found after checking all $v \in \hat{V}_u$
    , then $u$ keeps its initial match $\bar{v}(u)$. 
\end{itemize}


We next derive the probabilities $\{\hat{\mathbb{P}}(k) \mid 1 \le k \le m\}$ for the refinement process in the following two subsections. We first compute the probability that the $i$-th candidate $u$ swaps its initial match from the $k'$-th nearest neighbor to the $k''$-th nearest neighbor, based on the feasible swapping conditions described above.
Since not every $k''$-th nearest neighbor of $u$ is a valid swap candidate, we then derive the probability that $u$ successfully swaps from its $k'$-th neighbor to its $k$-th neighbor, following the sequential checking rule in refinement step (ii).

\subsubsection*{Probability for a feasible swap.}
For a swap from $\bar{v}(u)$ to $v$ to be feasible for $u$, where $K_u[\bar{v}(u)]=k', K_u[v]=k''$, all the aforementioned conditions (a)-(d) must be satisfied simultaneously. 
Note that for a randomly sampled $u$, the candidate $u'$ and the associated $\bar{v}(u')$ and $v'$ for a potential swap may vary across different problem instances. 
We first focus on condition (d) and derive the probability that it holds, conditional on the given positions of $\bar{v}(u')$ and $v'$: $K_{u'}[\bar{v}(u')] = q''$ and $K_{u'}[v'] = q'$. 
Then, we derive the probability for conditions (a)–(d) to be simultaneously satisfied, considering all possible combinations of the positions of $\bar{v}(u')$ and $v'$.
The corresponding results are presented in the following two lemmas. 

\begin{lemma} \label{lemma:P_swap_one} 
Given $K_u[\bar{v}(u)]=k', K_u[v]=k'', K_{u'}[\bar{v}(u')]=q'', K_{u'}[v']=q'$, where $1\le q''<q'\le n$ and $1\le k''<k'\le n$, the probability for condition (d) to hold is given by: 
\begin{align} \label{eq:P_swap_one}
\begin{split}
    \Pr\{\Delta_{q'',q'}\le\Delta_{k'',k'}\} = 
    \int_{0}^{\infty} \left[ \int_0^y f_{\Delta_{q'',q'}}(x) \text{ d}x \right] f_{\Delta_{k'',k'}}(y) \text{ d}y, 
\end{split}
\end{align}
where 
\begin{align} \label{eq:f_spacing}
\begin{split}
    f_{\Delta_{i,j}}(y) = 
    \frac{n!}{(i-1)!(j-i-1)!(n-j)!} 
    \int_{-\infty}^{\infty} 
    & \left[F_{C}(x)\right]^{i-1}\left[F_{C}(x+y)-F_{C}(x)\right]^{j-i-1} \\
    & \cdot \left[1-F_{C}(x+y)\right]^{n-j}f_{C}(x)f_{C}(x+y)
    \text{ d}x, 
    \quad y\ge 0.
\end{split}
\end{align}
\end{lemma}
\noindent{\bf Proof}.
Let random variable $\Delta_{i,j}\ge 0$ denote the difference between the cost to one's $j$-th nearest neighbor versus the $i$-th nearest neighbor, where $1 \le i < j \le n$. 
By definition, $\Delta_{i,j}$ represents the ``spacing" between the $j$-th and $i$-th order statistics for a sample of size $n$ drawn from distribution $F_{C}(\cdot)$, and its probability density function (p.d.f.), per \citet{arnold_first_2008}, 
is given by Equation \eqref{eq:f_spacing}. 
Since the probability for condition (d) to hold is $\Pr\{\Delta_{q'',q'}\le\Delta_{k'',k'}\}$, we compute the CDF of the difference $\Delta_{q'',q'} - \Delta_{k'',k'}$ from the convolution of their respective p.d.f.s, which is Equation \eqref{eq:P_swap_one}.
\hfill \qedsymbol\\

\begin{lemma} \label{lemma:P_swap_feasible}
Given $K_u[\bar{v}(u)]=k', K_u[v]=k''$, where $1 \le k'' < k' \le i$, the probability for conditions (a)–(d) to hold simultaneously for the $i$-th sorted demand vertex $u$ to possibly swap its match from $\bar{v}(u)$ to $v$ is:
\begin{align} \label{eq:P_swap_kp_kpp}
s(k', k'' \mid i) = 
\frac{1}{i-1} \sum_{q'' = 1}^{i - 1} \sum_{q' = q'' + 1}^{m} \frac{\binom{n - q'}{m - q'}}{\binom{n - q''}{m - q'' - 1}} 
\cdot \Pr\{\Delta_{q'', q'} \le \Delta_{k'', k'}\}, 
\quad 1 \le k'' < k' \le i.
\end{align}
where $\Pr\{\Delta_{q'', q'} \le \Delta_{k'', k'}\}$ is given by Equation \eqref{eq:P_swap_one}.
\end{lemma}
\noindent{\bf Proof}.
Given Lemma~\ref{lemma:P_swap_one}, the key to this proof is to derive the probability that $\bar{v}(u')$ and $v'$ take specific positions across all problem instances; i.e., $\Pr\{K_{u'}[\bar{v}(u')]=q'', K_{u'}[v'] = q' \}$, based on conditions (a)-(c). 

From condition (a), and given all the supply vertices are probabilistically identical, $\bar{v}(u')$ has an equal probability to take a position from $\{1,\ldots,i-1\}$, i.e., 
\begin{align} \label{eq:P_qpp}
\Pr\{K_{u'}[\bar{v}(u')]=q''\} = \frac{1}{i-1}, \quad q''\in\{1,\ldots,i-1\}.
\end{align}
Then, based on conditions (b)-(c), the conditional probability for $v'$ to take a position from $\{q''+1,\ldots,m\}$, i.e., $\Pr\{K_{u'}[v'] = q' \mid K_{u'}[\bar{v}(u')] = q''\}$, can be derived as follows. 
From the perspective of vertex $u'$, its first $q''$ nearest neighbors have already been identified, implying that $v'$ must be selected from the remaining $n - q''$ vertices. 
Moreover, these $q''$ neighbors have already been matched to $q''$ demand vertices (not including $u$), leaving $m - q''-1$ demand vertices available for matching with the supply vertices between $\bar{v}(u')$ and $v'$.
As such, the event $K_{u'}[v'] = q'$ occurs if and only if the $q'$-th nearest neighbor of $u'$ is either matched to $u$ or remains unmatched (considered as a ``failure''), while each of the $(q''+1)$-th through $(q'-1)$-th nearest neighbors is matched to one of the $m - q'' - 1$ demand vertices (considered as a ``success'').
Similarly to the derivation of Equation~\eqref{eq:P_k_i}, based on negative hyper-geometric distribution $\text{NHG}_{n - q'', m - q'' - 1, 1}(q' - q'' - 1)$, we have:
\begin{align} \label{eq:P_qp_qpp}
\begin{split}
    \Pr\{K_{u'}[v']=q' \mid K_{u'}[\bar{v}(u')]=q''\}
    = \frac{\binom{m-q''-1}{q'-q''-1}}{\binom{n-q''}{q'-q''-1}} \cdot \frac{n-q''-(m-q''-1)}{n-q''-(q'-q''-1)}&= \frac{\binom{n - q'}{m-q'}}{\binom{n - q''}{m-q''-1}}, \\
    & \quad q'\in\{q''+1,\ldots,m\}.
\end{split}
\end{align}

Then, based on the law of total probability, $s(k', k'' \mid i)$ can be computed by summing the following product 
$$
\Pr\{K_{u'}[\bar{v}(u')]=q''\} \cdot \Pr\{K_{u'}[v']=q' \mid K_{u'}[\bar{v}(u')]=q''\} \cdot \Pr\{\Delta_{q'', q'} \le \Delta_{k'', k'}\}
$$ 
across all possible values of $q''$ and $q'$, which leads to Equation \eqref{eq:P_swap_kp_kpp}. 
\hfill \qedsymbol\\

\subsubsection*{Probability for a successful swap.}

Lemma \ref{lemma:P_swap_feasible} provides the probability that the $i$-th candidate $u$ can possibly swap its initial greedy match (the $k'$-th nearest neighbor) to a better alternative (the $k''$-th nearest neighbor). By sequentially checking the feasibility of each of $u$'s $k''$-th nearest neighbor $v\in \hat{V}_u$, where $k''\in\{1,\ldots,k'-1\}$, we can derive the probability that $u$ successfully swaps to its $k$-th nearest neighbor. Then, by considering all possible values of $i \in \{1\ldots,m\}$, we can obtain an estimate of $\hat{\mathbb{P}}(k)$ as presented in the following lemma. 

\begin{lemma}
For an RBMP-I given $n, m \in \mathbb{Z}^+$, the refined matching probability $\hat{\mathbb{P}}(k)$ is given by: 
\begin{align} \label{eq:P_k_swap}
\hat{\mathbb{P}}(k) = 
\frac{1}{m} \cdot 
\sum_{i=k}^m 
\sum_{k'=k}^i 
\frac{\binom{n - k'}{i - k'}}{\binom{n}{i - 1}} \cdot s(k', k \mid i) \prod_{k''=1}^{k-1} \left[1 - s(k', k'' \mid i)\right]
, \quad k \in \{1, \cdots, m\},
\end{align}
where $s(k', k'' \mid i), \forall k'' < k'$ is given by Equation \eqref{eq:P_swap_kp_kpp}, and $s(k', k''=k' \mid i) = 1$.

\end{lemma}

\noindent{\bf Proof}.
From the perspective of the $i$-th candidate $u$ initially matched to its $k'$-th nearest neighbor $\bar{v}(u)$, given any $v\in \hat{V}_u$, the potential swapping candidate $u'$ and the associated $\bar{v}(u')$ and $v'$ vary. 
Since all vertices are probabilistically identical, we assume that the feasibility of swaps between $\bar{v}(u)$ and each $v\in \hat{V}_u$ are independent.\footnote{Strictly speaking, correlation could exist among feasible swaps but it tends to be weak. We choose to simplify the formula by making the independence assumption. Numerical experiments show later that such an assumption is reasonable.} 
The rule to perform a successful swap to its $k$-th neighbor can be interpreted as conducting a sequence of Bernoulli trials. Specifically, for each $k'' \in \{1, \cdots, k-1\}$, the trial succeeds with probability $s(k', k'' \mid i)$ and fails with probability $1 - s(k', k'' \mid i)$. The process continues until the first successful swap is identified or all $k-1$ trials fail. 
As such, the probability that $u$ successfully swap its match from the $k'$-th to the $k$-th nearest neighbor is given by:
\begin{align*} 
s(k', k \mid i) \prod_{k''=1}^{k-1} \left[1 - s(k', k'' \mid i)\right], \quad k' \in \{k+1, \cdots, i\}.
\end{align*}
Meanwhile, the probability that $u$ keeps its initial match, the $k$-th nearest neighbor, (i.e., no feasible swap occurs) is given by:
\begin{align*}
\prod_{k''=1}^{k-1} \left[1 - s(k', k'' \mid i)\right], \quad k'=k,
\end{align*}
If we define $s(k', k' \mid i) = 1$, the above two equations can be jointly written as: 
\begin{align}
\label{eq:P_swap_kp_k}
s(k', k \mid i) \prod_{k''=1}^{k-1} \left[1 - s(k', k'' \mid i)\right], \quad k' \in \{k, \cdots, i\}.
\end{align}

Following a similar logic as in the proof of Lemma \ref{lemma:P_k_greedy}, we first derive the conditional probability that the $i$-th candidate $u$ is ultimately matched to its $k$-th nearest neighbor after applying the swapping refinement, denoted by $\hat{\mathbb{P}}(k \mid i)$. 
For the trivial case when $i=1$, $u$ is initially matched to its nearest neighbor (i.e. $k=1$) and no refinement can occur, such that $\hat{\mathbb{P}}(k=1 \mid i=1)=\bar{\mathbb{P}}(k=1 \mid i=1)=1$. 
While when $i>1$, 
$\hat{\mathbb{P}}(k \mid i>1)$ can be computed as follows by summing over the probabilities that $u$ is initially matched to its $k'$-th (where $k\le k'\le i$) nearest neighbor, $\bar{\mathbb{P}}(k' \mid i)$, and then either successfully swaps to its $k$-th nearest neighbor or remains matched to its $k'$-th (where $k'=k$) nearest neighbor, as given by Equation \eqref{eq:P_swap_kp_k}. The resulting formula is as follows. 
\begin{align} \label{eq:P_k_swap_i}
\hat{\mathbb{P}}(k \mid i) =
\sum_{k'=k}^{i} \bar{\mathbb{P}}(k' \mid i) \cdot s(k', k \mid i) \prod_{k''=1}^{k-1} \left[1 - s(k', k'' \mid i)\right], \quad i\in\{1, \cdots,m\}, k \in \{1, \cdots, i\},
\end{align}
where $\bar{\mathbb{P}}(k' \mid i)$ is given by Equation \eqref{eq:P_k_i}. 
Then, by the law of total probability, the unconditional probability $\hat{\mathbb{P}}(k) = 
\sum_{i=k}^m \frac{1}{m} 
\hat{\mathbb{P}}(k\mid i)$ is given by Equation \eqref{eq:P_k_swap}.
\hfill \qedsymbol\\

\subsection{Main Results} \label{sec: formula}
We are now ready to present our main results. 
By putting all findings from the previous subsections together, we first obtain the distribution of the matching costs under the proposed matching processes for RBMP-I under arbitrary edge cost distribution $F_{C}(\cdot)$. Then, we focus on a special case when the edge costs follow a class of commonly seen power-law distributions. 
This class of distributions is of particular interest because they are related to the distribution of ``distances" between uniformly distributed vertices in a $D$-dimensional space, which will be closely related to RBMP-S and RBMP-B later. For this special case, we propose approximate closed-form formulas for both the matching probabilities and the expected matching costs. Based on these approximations, we identify a desirable convergence property: the expected greedy matching cost tend to converge to the expected optimal matching cost as $n \gg m$. This property reveals the conditions under which the proposed matching probabilities and cost formulas are applicable, which not only applies to RBMP-I, but also is useful for estimating distances in RBMP-S and RBMP-B. 

\subsubsection{Arbitrary cost distribution.}

The following results hold for RBMP-I under an arbitrary edge cost distribution $F_{C}(\cdot)$. 
First, by replacing the optimal matching probability $\accentset{\ast}{\mathbb{P}}(k)$ with $\bar{\mathbb{P}}(k)$ and $\hat{\mathbb{P}}(k)$ in Proposition \ref{prop:X_optimal}, we can respectively estimate the CDFs and the $M$-th moments of the matching costs under the greedy matching process, denoted by $\bar{C}(m,n)$ or $\bar{C}$, and after the vertex-swapping refinement, denoted by $\hat{C}(m,n)$ or $\hat{C}$, as follows.
\begin{proposition} \label{prop:X_greedy_swap}
For an RBMP-I with given $n, m \in \mathbb{Z}^+$, and $F_{C}(\cdot)$, the CDFs and the $M$-th moments of the greedy and refined matching costs are respectively:
\begin{align} 
    F_{\bar{C}}(x) 
    = 
    \sum_{k=1}^m I_{F_{C}(x)}(k,n-k+1)\cdot \bar{\mathbb{P}}(k)
    , \quad 
    \mathbb{E}[\bar{C}^M] 
    &= 
    \sum_{k=1}^m \mathbb{E}[C^M \mid k] \cdot \bar{\mathbb{P}}(k), \label{eq:E_X_greedy_prop}\\
    F_{\hat{C}}(x) 
    = 
    \sum_{k=1}^m I_{F_{C}(x)}(k,n-k+1)\cdot \hat{\mathbb{P}}(k)
    , \quad 
    \mathbb{E}[\hat{C}^M] 
    &= 
    \sum_{k=1}^m \mathbb{E}[C^M \mid k] \cdot \hat{\mathbb{P}}(k), \label{eq:E_X_swap_prop}
\end{align}
where $\mathbb{E}[C^M \mid k]$ is given by Equation \eqref{eq:E_X_k}, $\bar{\mathbb{P}}(k)$ and $\hat{\mathbb{P}}(k)$ are given by Equations \eqref{eq:P_k_greedy} and \eqref{eq:P_k_swap}.
\end{proposition}

According to Propositions \ref{prop:X_optimal} and \ref{prop:X_greedy_swap}, the expected optimal, greedy and refined matching costs can be obtained by simply setting $M = 1$ in the corresponding equations. 
Since both the expected greedy and refined costs, $\mathbb{E}[\bar{C}]$ and $\mathbb{E}[\hat{C}]$, are derived from processes that generate feasible matching solutions, they serve as upper bounds on the expected optimal cost $\mathbb{E}[\accentset{\ast}{C}]$.
In addition, $\mathbb{E}[\bar{C}]$ must always be greater than or equal to $\mathbb{E}[\hat{C}]$, as the vertex-swapping refinement process is designed to reduce (or at least not increase) the total matching cost. 
In the special case where the number of demand vertices is $m = 1$, the greedy matching process already yields the optimal solution (matched with its nearest neighbor), and no refinements are needed. As such, all three expected matching costs, $\mathbb{E}[\bar{C}(1,n)]$, $\mathbb{E}[\hat{C}(1,n)]$ and $\mathbb{E}[\accentset{\ast}{C}(1,n)]$, are equal. 
Furthermore, for any given value of $n$, the expected matching cost in the case $m = 1$ must be strictly smaller than those with any $m>1$, since this ``nearest" matching cost represents the best achievable cost for each demand vertex. These relationships among the expected matching costs are summarized in the following proposition.
\begin{proposition} \label{prop:X_inequality}
For an RBMP-I with given $n, m \in \mathbb{Z}^+$, the expected greedy, refined, optimal and nearest matching costs must satisfy the following:
\begin{align} \label{eq:E_Xs_compare}
    \mathbb{E}[\bar{C}(m,n)] \ge \mathbb{E}[\hat{C}(m,n)] \ge \mathbb{E}[\accentset{\ast}{C}(m,n)] \ge \mathbb{E}[\bar{C}(1,n)].
\end{align}
\end{proposition}

\subsubsection{Special case: $F_C(\cdot)$ following a power-law distribution.}
\label{sec:power-law}
Next, we present a set of approximate closed-form formulas for estimating the matching probabilities and costs of RBMP-I under a class of power-law distributions: 
$$F_{C}(x)=\left(\frac{x}{R}\right)^D, \quad x\in [0,R],$$ where $R>0$ and $D\in \mathbb{Z^+}$ are given parameters. 
Note that for an arbitrary $F_{C}(\cdot)$, closed-form formulas may not be available, particularly because Equations \eqref{eq:E_X_k} and \eqref{eq:P_swap_one} involve integrals that may not be readily simplified into elementary functions.
Yet, simpler formulas can be derived through approximation for the special power-law distribution case. 
These approximations allow us to reveal a convergence property between the greedy and nearest matching costs, which further leads to the approximations for the optimal matching cost.

\subsubsection*{Matching probability approximation.}
\label{sec:P_approx}
First, we propose a normal approximation in estimating $\Pr\{\Delta_{q'',q'}\le\Delta_{k'',k'}\}$ in Equation \eqref{eq:P_swap_one}.
When $D=1$, the edge cost distribution becomes uniform, and it is well known \citep{arnold_first_2008} that the scaled spacing between two costs, $\frac{\Delta_{i,j}}{R}$, follows a Beta distribution: Beta$(j-i, n-j+i+1)$. In this case, the expectation and variance of $\Delta_{i,j}$ are respectively given by:
\begin{align}
    \mathbb{E}[\Delta_{i,j}] = \frac{j-i}{n+1}\cdot R, \quad
    \mathbb{V}[\Delta_{i,j}] = \frac{(j-i)(n - j + i + 1)}{(n+1)^2(n+2)}\cdot R^2. \label{eq:mu_sigma}
\end{align}
According to \citet{wise_normalizing_1960}, under certain conditions (such as $j-i \gg 1$ and $n\gg j-i$), $\Delta_{i,j}$ is approximately normal. In these cases, $\Delta_{q'',q'}$ and $\Delta_{k'',k'}$ can be treated as two independent normal variables, such that the probability $\Pr\{\Delta_{q'',q'}\le\Delta_{k'',k'}\}$ can be approximated by the following.
\begin{align}
    \Pr\{\Delta_{q'',q'}\le\Delta_{k'',k'}\} \approx \Phi\left(\frac{\mathbb{E}[\Delta_{k'',k'}] - \mathbb{E}[\Delta_{q'',q'}]}{\sqrt{\mathbb{V}[\Delta_{q'',q'}] + \mathbb{V}[\Delta_{k'',k'}]}}\right) \label{eq:P_normal},
\end{align}
where $\Phi(\cdot)$ represents the standard normal CDF. 

When $D\neq 1$, the same normal approximation applies while the following expectation and variance of $\Delta_{i,j}$ shall be used in Equation \eqref{eq:P_normal}.
\begin{align}
\label{eq:E_Delta}
    \mathbb{E}[\Delta_{i,j}] = \int_0^{\infty} yf_{\Delta_{i,j}}(y)\text{ d}y, 
    \quad 
    \mathbb{V}[\Delta_{i,j}] = \int_0^{\infty} y^2f_{\Delta_{i,j}}(y)\text{ d}y - \mathbb{E}^2[\Delta_{i,j}],
\end{align}
where $f_{\Delta_{i,j}}(y)$ is given by Equation \eqref{eq:f_spacing}.
For computational simplicity, even for $D > 1$, we propose to use $\mathbb{E}[\Delta_{i,j}]$ and $\mathbb{V}[\Delta_{i,j}]$ from Equation \eqref{eq:mu_sigma} instead to estimate $\Pr\{\Delta_{q'',q'}\le\Delta_{k'',k'}\}$. 
As discussed in the numerical experiment section, these approximations yield reasonably good results in practice.

\subsubsection*{Greedy matching cost.} 
\label{sec:greedy_cost}
We next present an approximation for the $M$-th moment of the conditional matching cost in Equation \eqref{eq:E_X_k}.
This approximation gives a set of simplified formulas for the expected greedy matching cost and ultimately leads us to identify the convergence property.

Under the given power-law distribution, it is known \citep{arnold_first_2008} that Equation \eqref{eq:E_X_k} can be simplified into the following form:
\begin{align}
    &\mathbb{E}[C^M\mid k]
    = \frac{\Gamma(n+1)\Gamma(k+\frac{M}{D})}{\Gamma(n+1+\frac{M}{D})\Gamma(k)} \cdot R^M, \quad k\in\{1, \ldots, n\}. \label{eq:E_X_k_power}
\end{align}

Note that Equation \eqref{eq:E_X_k_power} can be substituted into Equation \eqref{eq:E_X_greedy_prop} in Proposition \ref{prop:X_greedy_swap} to obtain the exact formula for the expected greedy matching cost, $\mathbb{E}[\bar{C}]$.
Consider the case when $n \gg 1$, we can use the approximate matching probability $\bar{\mathbb{P}}(k)$ from Equation \eqref{eq:P_k_geometric} to have an approximate formula for $\mathbb{E}[\bar{C}]$, as follows.
\begin{align} \label{eq:E_X_greedy_geo}
\begin{split}
\mathbb{E}[\bar{C}] 
    \xrightarrow{n\gg 1} \frac{R\cdot\Gamma(n+1)}{m\cdot\Gamma(n+1+\frac{1}{D})} \sum_{i=1 }^{m} 
        & 
        \left[ \sum_{k=1}^{i} \left( \frac{i-1}{n} \right)^{k-1} \left( 1- \frac{i-1}{n} \right) \frac{\Gamma(k+\frac{1}{D}) }{\Gamma(k)}  
        +  
        \left( \frac{i-1}{n} \right)^{i}  \frac{\Gamma(i+\frac{1}{D}) }{\Gamma(i)}  \right].
\end{split}
\end{align}
Then, consider the ratio of two Gamma functions can be simplified using Stirling's approximation for sufficiently large $z$, as follows:
\begin{align} \label{eq:gamma_ratio}
\frac{\Gamma(z+\frac{1}{D})}{\Gamma(z)} \xrightarrow{z \gg 1} z^{\frac{1}{D}}.
\end{align}
We introduce a parameter $\kappa \in \{0, \ldots, m\}$ to specify the minimum value of $z$ for which the gamma approximation in Equation \eqref{eq:gamma_ratio} is applied, and obtain a further approximation of $\mathbb{E}[\bar{C}]$, denoted by $\mathbb{E}_\kappa[\bar{C}]$. 
\begin{align} \label{eq:E_X_greedy_kappa}
    \mathbb{E} [\bar{C}] \lessapprox
    \mathbb{E}_\kappa [\bar{C}] 
        & = \frac{R}{mn^{\frac{1}{D}}} \sum_{i=1 }^{m} 
        \left[ \sum_{k=1}^{i} \left( \frac{i-1}{n} \right)^{k-1} \left( 1- \frac{i-1}{n} \right) 
\left( \mathbf{1}_{k \le \kappa} \cdot \frac{\Gamma(k+\frac{1}{D}) }{\Gamma(k)} + \mathbf{1}_{k > \kappa} \cdot k^{\frac{1}{D}}\right) \right. \nonumber\\
        & +  \left. \left( \frac{i-1}{n} \right)^{i} \left( \mathbf{1}_{i \le \kappa} \cdot \frac{\Gamma(i+\frac{1}{D}) }{\Gamma(i)} + \mathbf{1}_{i > \kappa} \cdot i^{\frac{1}{D}}\right) \right]. %
\end{align}
The error bounds associated with this ``$\kappa$-approximation" to $\mathbb{E}[\bar{C}]$ are detailed in Appendix \ref{app:kappa_approx}, which implies the approximation is reasonably good. 
In addition, the value of $\kappa$ should be chosen to find a good trade-off between accuracy and computation efficiency.
Note here setting $\kappa = m$ yields the most accurate approximation, however requires the most computational effort.
To the extreme, we may set $\kappa = 0$; i.e., use Equation \eqref{eq:gamma_ratio} to approximate all gamma function ratios. In our numerical experiments, we find that setting $\kappa = 0$ provides a reasonably accurate approximation. 

Based on the approximations developed above, a desirable convergence property of the matching costs under the power-law distribution can be identified, where the greedy, refined, optimal, and nearest matching costs in Proposition \ref{prop:X_inequality} tend to converge when $n \gg m$. 
We first derive an error bound between the greedy matching cost $\mathbb{E}[\bar{C}]$ and the nearest matching cost $\mathbb{E}[\bar{C}(1,n)]$, as stated in the following lemma. The detailed proof can be found in Appendix \ref{app:error_upper_lower}.
\begin{lemma} \label{lemma:error_upper_lower}
For an RBMP-I given $n, m \in \mathbb{Z}^+$, and $F_{C}(x) = \left(\frac{x}{R}\right)^D$, where $R>0$ and $D\in \mathbb{Z}^+$,
\begin{align}\label{eq:error_upper_lower}
    \frac{\mathbb{E}[\bar{C}(m,n)] - \mathbb{E}[\bar{C}(1,n)]}{\mathbb{E}[\bar{C}(1,n)]} \le \frac{\left( -\frac{n}{m} \right) \ln \left(1-\frac{m}{n}\right)-1}{\Gamma(1+\frac{1}{D})} .
\end{align}
\end{lemma}

Note that when $n \gg m$, 
$\ln \left(1-\frac{m}{n}\right) \approx -\frac{m}{n}$, 
the error bound on the right-hand side of Equation \eqref{eq:error_upper_lower} quickly diminishes to zero. This implies that, when $n \gg m$, the upper and lower bounds in Proposition \ref{prop:X_inequality} converge, and by the squeeze (sandwich) theorem, the intermediate quantities (i.e., the refined and optimal matching costs) also converge, as shown in the following proposition. 

\begin{proposition} \label{prop:X_greedy_opt}
For an RBMP-I given $n, m \in \mathbb{Z}^+$, and $F_{C}(x) = \left(\frac{x}{R}\right)^D$, where $R>0$ and $D\in \mathbb{Z}^+$, when $n \gg m$, the expected greedy, refined, optimal and nearest matching costs converge.
\begin{align}
\label{eq:convergence}
    \mathbb{E}[\bar{C}(m,n)] \xrightarrow{n\gg m} \mathbb{E}[\hat{C}(m,n)] \xrightarrow{n\gg m} \mathbb{E}[\accentset{\ast}{C}(m,n)] \xrightarrow{n\gg m} \mathbb{E}[\bar{C}(1,n)] 
    = R\cdot \Gamma(1+\frac{1}{D})\cdot n^{-\frac{1}{D}}.
\end{align}
\end{proposition}

\subsubsection*{Optimal matching cost.}
The above convergence property, together with Propositions \ref{prop:X_optimal} and \ref{prop:X_greedy_swap}, suggests that the greedy and refined matching probabilities, $\bar{\mathbb{P}}(k)$ and $\hat{\mathbb{P}}(k)$, can serve as effective approximations to the optimal matching probability $\accentset{\ast}{\mathbb{P}}(k)$ across different scenarios.

When $n\gg m$, $\bar{\mathbb{P}}(k)$ can be used as an effective approximation for $\accentset{\ast}{\mathbb{P}}(k)$; i.e.,
\begin{align} \label{eq:P_greedy_approx}
    \accentset{\ast}{\mathbb{P}}(k) \xrightarrow{n\gg m} \bar{\mathbb{P}}(k)
    \approx
    \frac{1}{m} \left(\frac{k-1}{n}\right)^{k-1} + \frac{1}{m} \sum_{i=k+1}^{m}  \left(\frac{i-1}{n}\right)^{k-1}\left(1-\frac{i-1}{n}\right), \quad k \in \{1, \cdots, m\}.
\end{align}
For arbitrary values of $m$ and $n$, the refined matching cost $\mathbb{E}[\hat{C}]$ can be expected to be near-optimum, making $\hat{\mathbb{P}}(k)$ a practical approximation to $\accentset{\ast}{\mathbb{P}}(k)$; i.e.,
\begin{align} \label{eq:P_swap_approx}
    \accentset{\ast}{\mathbb{P}}(k) \approx \hat{\mathbb{P}}(k) \approx \frac{1}{m} \cdot \sum_{i=k}^m \sum_{k'=k}^i \frac{\binom{n - k'}{i - k'}}{\binom{n}{i - 1}} \cdot s(k', k \mid i) \prod_{k''=1}^{k-1} \left[1 - s(k', k'' \mid i)\right], \quad k \in \{1, \cdots, m\}.
\end{align}
Here, when computing $s(k', k'' \mid i)$, we may approximate $\Pr\{\Delta_{q'',q'} \le \Delta_{k'',k'}\}$ using the normal approximation in Equations \eqref{eq:mu_sigma} and \eqref{eq:P_normal}, and obtain the following:
\begin{align}
\label{eq:s_normal}
    s(k', k'' \mid i) \approx 
\begin{cases}
    \frac{1}{i-1} \sum_{q'' = 1}^{i - 1} \sum_{q' = q'' + 1}^{m} \frac{\binom{n - q'}{m - q'}}{\binom{n - q''}{m - q'' - 1}} 
    \cdot \Phi\left(\frac{(k'-k''-q'+q'')\sqrt{n+2}}{\sqrt{(q'-q'')(n-q'+q''+1)+(k'-k'')(n-k'+k''+1)}}\right), \\
    \qquad\qquad\qquad\qquad\qquad\qquad\qquad\qquad\qquad\qquad
    k'' \in \{1, \cdots, k'-1\},\\
    1, \quad k''=k'.
\end{cases}
\end{align}

In summary, the $M$-th moment of the optimal matching cost for RBMP-I under a power-law distribution can be approximately estimated by the formulas below. When $n\gg m$, 
\begin{align}
    \mathbb{E}[\accentset{\ast}{C}^M ]
    \xrightarrow{n\gg m} \frac{R^M\cdot\Gamma(n+1)}{\Gamma(n+1+\frac{M}{D})}\sum_{k=1}^{m} \frac{
    \Gamma(k+\frac{M}{D})}{
    \Gamma(k)}\cdot \bar{\mathbb{P}}(k),
    \label{eq:E_X_opt_special}
\end{align}
where $\bar{\mathbb{P}}(k)$ is given by Equation \eqref{eq:P_greedy_approx}. 
While for arbitrary values of $m$ and $n$, 
\begin{align}
    \mathbb{E}[\accentset{\ast}{C}^M]
    \approx \frac{R^M\cdot\Gamma(n+1)}{\Gamma(n+1+\frac{M}{D})}\sum_{k=1}^{m} \frac{
    \Gamma(k+\frac{M}{D})}{
    \Gamma(k)}\cdot \hat{\mathbb{P}}(k),
    \label{eq:E_X_opt_general}
\end{align}
where $\hat{\mathbb{P}}(k)$ is given by Equation \eqref{eq:P_swap_approx}. 

As a final remark, although the above properties are established for RBMP-I under the power-law distributions, the results are expected to hold (at least approximately) for other distance distributions that have similar behaviors, such as those in RBMP-S and RBMP-B.

\section{RBMP-S}
\label{sec:RBMP-S}
This section presents approximate formulas for estimating the optimal matching distance in RBMP-S. We use the results from RBMP-I as a starting point, but 
further account for the effects of spatial correlation among vertex distributions and matching. 

\subsection{Problem Definition}
We begin by formally defining an RBMP-S. The key feature of this variant is that the vertices in $U$ and $V$ represent independently and uniformly distributed ``physical points" on a unit-area $D$-dimensional hyper-sphere $\mathcal{S}_D$ embedded in a $(D+1)$-dimensional Euclidean space $\mathbb{R}^{D+1}$, where $D\in \mathbb{Z}^+$. 
For example, $\mathcal{S}_1$ represents a circular perimeter of a round disk; $\mathcal{S}_2$ represents the surface of a spherical ball (e.g., see Figure \ref{fig:RBMP_model}). 
Given any $D\in \mathbb{Z}^+$, the radius of the hyper-sphere $\mathcal{S}_D$, denoted $R_{\text{S}}$, is
\citep{olver_nist_2010}:
\begin{align}
    R_{\text{S}} = 
    \pi\left(\frac{\Gamma\left(\frac{D+1}{2}\right)}{2\pi^{\frac{D+1}{2}}}\right)^{\frac{1}{D}}.
\end{align}
For each realization of $U$ and $V$, the edge cost between any two vertices $u\in U$ and $v\in V$ is computed as their great-circle distance, which is related to the convex central angle between them, $\theta_{u,v}$, as follows: 
\begin{align} \label{eq:c_RBMP-S}
    c_{u,v} = \frac{R_{\text{S}}}{\pi}\cdot\theta_{u,v}, \quad \theta_{u,v} \in[0,\pi].
\end{align}
Then, the same matching problem, as defined in Equations~\eqref{eq:RBMP_obj}–\eqref{eq:RBMP_con2}, is solved for each realization, and the average of the optimal matching distance is is taken across all realizations. 

As such, unlike in RBMP-I, where the edge cost distribution $F_{C}(\cdot)$ is directly given, an RBMP-S is now characterized by the vertex set sizes $m$ and $n$, and the spatial dimension $D$, where the CDF of the edge cost, denoted by $F_{C_{\text{S}}}(\cdot)$, depends on these three parameters. 
Our primary interest is in estimating the expectation of the optimal matching cost (or distance), denoted by $\accentset{\ast}{C_{\text{S}}}(m,n,D)$, or simply $\accentset{\ast}{C_{\text{S}}}$. 

Following the main idea for RBMP-I in Section \ref{sec:RBMP-I}, the estimation of $\mathbb{E}[\accentset{\ast}{C_{\text{S}}}]$ includes: 
\begin{enumerate}
    \item[(i)] Estimating the probability for a randomly sampled demand vertex being matched to its $k$-th nearest neighbor. It is important to note that in RBMP-S, although the vertices on $\mathcal{S}_D$ are i.i.d., the edge costs (distances) are inherently correlated due to the constraints imposed by a geometric metric space (e.g., triangle inequality). To address this, we first assume i.i.d. edge costs as a first approximation, such that the matching probability can be given by $\accentset{\ast}{\mathbb{P}}(k)$ from its RBMP-I counterpart. 
    \item[(ii)] Estimating the expected distance from a randomly sampled demand vertex on $\mathcal{S}_D$ to its $k$-th nearest neighbor, $\mathbb{E}[C_{\text{S}} \mid k]$. 
    This can be derived from Equation \eqref{eq:E_X_k} once $F_{C_{\text{S}}}(x)$ is identified, as will be discussed later in Section \ref{sec:dist_k_RBMP-S}. 
    \item[(iii)] As an additional step, we relax the i.i.d. cost assumption and introduce a correction factor, denoted by $\delta_{\text{S}}(m,n,D)$ (or simply $\delta_{\text{S}}$), to account for the impacts of spatial correlation in RBMP-S across different parameter combinations, as will be discussed later in Section \ref{sec:correction_RBMP-S}.
\end{enumerate}

\subsection{Matching Distance Formula}
\label{sec:dist_RBMP-S}
Based on the results from the aforementioned three steps, $\mathbb{E}[\accentset{\ast}{C}_{\text{S}}]$ can be estimated on top of Proposition \ref{prop:X_optimal} as follows:
\begin{align} \label{eq:E_Y_optimal}
\mathbb{E}[\accentset{\ast}{C}_{\text{S}}] = 
\left( 1 + \delta_{\text{S}} \right)\sum_{k=1}^m \mathbb{E}[C_{\text{S}}\mid k]\cdot \accentset{\ast}{\mathbb{P}}(k).
\end{align} 
We next propose a set of approximate formulas applicable to different supply-demand scenarios. 

First, consider the scenario when $n\gg m$. 
In the trivial case with only $m = 1$ demand vertex, no correlation arises for the single vertex match, and hence $\delta_{\text{S}} = 0$. 
When $n \gg m$, the impact of correlation is also negligible, because each demand vertex is likely to be surrounded by sufficient supply vertices nearby, so that its matching distance would be minimally dependent on other demand vertices.
As such, we expect that 
$$\delta_{\text{S}} \xrightarrow{n\gg m} 0.$$
In addition, consider when $n$ or $m$ is sufficiently large, the distances between nearby vertices become very small, and the great-circle distance can be well approximated by the Euclidean distance. Under this approximation, the distribution $F_{C_{\text{S}}}(x)$ converges to a power-law form; i.e.,
$$F_{C_{\text{S}}}(x)
\xrightarrow{x\to 0}
\left( \frac{x}{R_{\text{S}}} \right)^D, \quad x \in [0, R_{\text{S}}].$$ 
In this case, Proposition \ref{prop:X_greedy_opt} applies, indicating that the expected greedy, optimal and nearest matching distances converge; i.e.,
\begin{align}
\label{eq:convergence_S}
\mathbb{E}[\bar{C_{\text{S}}}(m,n,D)] \xrightarrow{n\gg m} \mathbb{E}[\accentset{\ast}{C_{\text{S}}}(m,n,D)] \xrightarrow{n\gg m} \mathbb{E}[\bar{C_{\text{S}}}(1,n,D)].
\end{align}

As discussed in the final remark of Section \ref{sec:power-law}, although $F_{C_{\text{S}}}(x)$ in RBMP-S does not strictly follow a power-law form for general cases with arbitrary values of $m$ and $n$, it nevertheless captures the distance-related distribution between vertices and is expected to have a similar behavior. 
Accordingly, we propose to use the greedy matching probability $\bar{\mathbb{P}}(k)$ and the refined matching probability $\hat{\mathbb{P}}(k)$ derived for RBMP-I to approximate the optimal matching probability $\accentset{\ast}{\mathbb{P}}(k)$ in RBMP-S at different scenarios. Similarly to the results for RBMP-I in Equations \eqref{eq:E_X_opt_special} and \eqref{eq:E_X_opt_general}, the expecte optimal matching distance for RBMP-S can be approximately estimated by the formulas below. 

When $n\gg m$, 
\begin{align} \label{eq:E_Y_greedy}
\mathbb{E}[\accentset{\ast}{C}_{\text{S}}]
\xrightarrow{n\gg m}
\mathbb{E}[\bar{C}_{\text{S}}]
&=\sum_{k=1}^m \mathbb{E}[C_{\text{S}}\mid k]\cdot \bar{\mathbb{P}}(k),
\end{align}
where $\bar{\mathbb{P}}(k)$ is given by Equation \eqref{eq:P_greedy_approx}.
While for arbitrary values of $m$ and $n$, 
\begin{align} \label{eq:E_Y_swap}
\mathbb{E}[\accentset{\ast}{C}_{\text{S}}] \approx \left( 1 + \delta_{\text{S}} \right) \sum_{k=1}^m \mathbb{E}[C_{\text{S}}\mid k] \cdot \hat{\mathbb{P}}(k).
\end{align}
where $\hat{\mathbb{P}}(k)$ is given by Equation \eqref{eq:P_swap_approx}. The next two sections provide details for estimating the $k$-th nearest-neighbor distance $\mathbb{E}[C_{\text{S}}\mid k]$ and the correction factor $\delta_{\text{S}}$, respectively.

\subsection{$k$-th Nearest Neighbor Distance}
\label{sec:dist_k_RBMP-S}
Note $C_{\text{S}}$ now corresponds to the distance between two randomly sampled vertices from a total of $m + n$ vertices uniformly distributed on the unit sphere $\mathcal{S}_D$.
According to \citet{cai_distributions_2013}, for sufficiently large $m$ or $n$, the central angle $\theta_{u,v}$ between any two vertices on $\mathcal{S}_D$ becomes approximately i.i.d. with p.d.f.:
\begin{align} \label{eq:f_theta}
f_{\theta_{u,v}}(x) = \frac{\Gamma(\frac{D+1}{2})}{\sqrt{\pi} \Gamma(\frac{D}{2})}(\sin^{D-1}{x}),
\quad x\in[0,\pi].
\end{align}
The CDF of $C_{\text{S}}$ can be obtained by integrating the above p.d.f. and substituting $\theta_{u,v}$ with $c_{u,v}$ using the relationship defined in Equation \eqref{eq:c_RBMP-S}. As such, we have:  
\begin{align} \label{eq:F_Y}
    F_{C_{\text{S}}}(x)
    = I_{\sin^2\left(\pi 
    x/2R_{\text{S}}\right)}\left(\frac{D}{2}, \frac{D}{2}\right), 
    \quad x\in[0,R_{\text{S}}].
\end{align}
By setting $F_{C}(x) = F_{C_{\text{S}}}(x)$ in an RBMP-I as given by the above equation, and applying Equation \eqref{eq:E_X_k}, $\mathbb{E}[C_{\text{S}}\mid k]$ can be estimated as follows:
\begin{align} \label{eq:E_Y_k_general}
\begin{split}
    \mathbb{E}[C_{\text{S}}\mid k]
    = \int_{0}^{R_{\text{S}}}
    &\frac{x^M}{\text{B}(k,n-k+1)}\left[1-I_{\sin^2\left(\pi 
    x/2R_{\text{S}}\right)}\left(\frac{D}{2}, \frac{D}{2}\right)\right]^{n-k}\\
    &\cdot\left[I_{\sin^2\left(\pi 
    x/2R_{\text{S}}\right)}\left(\frac{D}{2}, \frac{D}{2}\right)\right]^{k-1}\text{ d}I_{\sin^2\left(\pi
    x/2R_{\text{S}}\right)}\left(\frac{D}{2}, \frac{D}{2}\right). 
\end{split}
\end{align}

In real-world transportation problems, when $D = 1$ and $D = 2$, 
the regularized beta functions in the CDFs in Equation \eqref{eq:F_Y} can be simplified into the following:
\begin{align} \label{eq:E_Y_k_special}
    I_{\sin^2\left(\pi 
    x/2R_{\text{S}}\right)}\left(\frac{1}{2}, \frac{1}{2}\right) = \frac{x}{R_{\text{S}}}, \quad
    I_{\sin^2\left(\pi 
    x/2R_{\text{S}}\right)}\left(1, 1\right) = \frac{1}{2}[1-\cos(\pi 
    x/R_{\text{S}})].
\end{align}
Specifically, when $D = 1$, it can be seen from the above equation that $C_{\text{S}}$ follows a uniform distribution, and $\mathbb{E}[C_{\text{S}}^M(m,n,1)\mid k]$ can be simplified into a closed-form by substituting $R$ with $R_{\text{S}}$ in Equation \eqref{eq:E_X_k_power}, as follows.
\begin{align} \label{eq:E_X_k_1D_RBMP-S}
    &\mathbb{E}[C_{\text{S}}(m,n,1)\mid k]=\frac{k}{n+1}\cdot R_{\text{S}}
    , \quad k\in\{1, \ldots, n\}. 
\end{align}

\subsection{Correction Factor Estimation}
\label{sec:correction_RBMP-S}
This section propose a specific functional form for the correction factor $\delta_{\text{S}}(m,n,D)$ based on several key results from the literature.

To the best of our knowledge, the only known analytical result that is tangentially related to RBMP-S is for the special case of asymptotic random mono-partite matching problem (RMMP), as in \citet{mezard_euclidean_1988}, where matching is performed within a single set of $n\rightarrow \infty$ vertices (rather than between two disjoint subsets).
In their study, a correction factor is also proposed to quantify the relative difference in the expected optimal matching costs between RMMP-I, where edge costs are i.i.d. and follow a power-law distribution $F_{C}(x) \xrightarrow{x \to 0} (x/R)^D$, and RMMP-S, where edge costs are defined over $\mathcal{S}_D$ and converge to the same distribution as $n \to \infty$. They showed that:
(i) the optimal matching costs of RMMP-I and RMMP-S scale similarly (as $n \to \infty$) under any fixed dimension $D$, and the ratio between them converges to a non-negative constant; 
and (ii) the correction factor decreases as $D \to \infty$ and remains relatively small even in low dimensions; for example, it is less than 10\% for $D \geq 2$.
These results offer valuable insights in estimating $\delta_{\text{S}}(m,n,D)$ for RBMPs, as mono-partite and bipartite matching problems share many similarities (at least regarding the impacts of spatial correlation). 

\citet{houdayer_comparing_1998} conducted extensive numerical simulations to quantify the relative differences between balanced RBMP-S and RBMP-I with very large (nearly infinite) numbers of vertices; i.e., $\delta_{\text{S}}(\infty,\infty,D)$. 
They had several key observations.
First, the scaling behavior of matching costs as $n \to \infty$ is the same in both RBMP-S and RBMP-I for $D \ge 3$, and the estimated $\delta_{\text{S}}(\infty,\infty,D)$ for $3 \le D \le 10$ are constants (as summarized in Table \ref{tab:deltas}). However, the scaling behavior is different for $D = 1$ and $D = 2$, and in these cases, $\delta_{\text{S}}(\infty,\infty,D)$ is no longer meaningful.  
Second, the correction factor $\delta_{\text{S}}(\infty,\infty,D)$ decreases as $D \rightarrow \infty$, and approximately scales with $\frac{1}{D^2}$. 
In addition, while no prior results exist for problems where $m \le n < \infty$, it is reasonable to expect that $\delta_{\text{S}}(m,n,D)$ decreases as the ratio $\frac{m}{n}$ decreases, since we know that  $\delta_{\text{S}}(m,n,D) \xrightarrow{n\gg m} 0$.

\begin{table}[ht!]
\centering
\caption{Values of $\delta_{\text{S}}(\infty,\infty,D)$ for $3\le D\le 10$ (Source: \citet{houdayer_comparing_1998}).}
\label{tab:deltas}
\resizebox{0.8\textwidth}{!}{%
\begin{tabular}{|c|c|c|c|c|c|c|c|c|}
\hline
$D$ & 3 & 4 & 5 & 6 & 7 & 8 & 9 & 10 \\ \hline
$n\cdot \mathbb{E}[\accentset{\ast}{C}(\infty,\infty)]$                        & 0.6537 & 0.6865 & 0.7244 & 0.7628 & 0.8006 & 0.8371 & 0.8724 & 0.9064 \\ \hline
$n\cdot \mathbb{E}[\accentset{\ast}{C_{\text{S}}}(\infty,\infty,D)]$              & 0.7080 & 0.7081 & 0.7349 & 0.7688 & 0.8039 & 0.8391 & 0.8736 & 0.9076 \\ \hline
$\delta_{\text{S}}(\infty,\infty,D)$ & 8.31\% & 3.15\% & 1.46\% & 0.78\% & 0.42\% & 0.24\% & 0.14\% & 0.13\% \\ \hline
\end{tabular}%
}
\end{table}

Based on all these observations, we propose the following functional form for $\delta_{\text{S}}(m,n,D)$:
\begin{align} \label{eq:delta}
    \delta_{\text{S}}(m,n,D) = \beta_{\text{S}}(n,D) \left(\frac{m}{n}\right)^{d} \frac{1}{D^2}.
\end{align}
Here, the term $\frac{1}{D^2}$ with respect to spatial dimension $D$ is from the second observation of \citet{houdayer_comparing_1998}.
The term $\beta_{\text{S}}(n,D)$, a scaling factor that depends on both $D$ and $n$, accounts for the spatial correlations arising from balanced matching problems.
Then, the term $\left(\frac{m}{n}\right)^d$ accounts for the impact of imbalance between the two vertex sets, where $d$ is a parameter.

We next estimate $\beta_{\text{S}}(n,D)$ based on best known results of $\delta_{\text{S}}(m,n,D)$ for balanced problems ($m = n$). We approach this for three scenarios: $D = 1$, $D = 2$, and $D \ge 3$. 
First, for a line with $D = 1$, we adopt the closed-form approximate formula for $\mathbb{E}[\accentset{\ast}{C_{\text{S}}}(n,n,1)]$ provided in \citet{zhai_average_2024} as a baseline estimation:
\begin{align} \label{eq:E_Y_1d}
\mathbb{E}[\accentset{\ast}{C_{\text{S}}}(n,n,1)]
\approx
\frac{1}{4}\sqrt{\frac{\pi}{2}} \cdot n^{-\frac{1}{2}}. 
\end{align} 
Combining Equations \eqref{eq:E_Y_swap} and \eqref{eq:delta}, we see that $\beta_{\text{S}}(n,D)$ equals to $\delta_{\text{S}}(n,n,1)$, which is defined as the relative difference between $\mathbb{E}[\accentset{\ast}{C_{\text{S}}}(n,n,1)]$ and its RBMP-I counterpart; i.e.,
\begin{align} \label{eq:beta_delta_RBMP-S}
\beta_{\text{S}}(n,1) = \delta_{\text{S}}(n,n,1) \approx
\frac{\mathbb{E}[\accentset{\ast}{C_{\text{S}}}(n,n,1)]}{\sum_{k=1}^n\mathbb{E}[{C_{\text{S}}}(n,n,1)\mid k]\cdot \hat{\mathbb{P}}(k)} - 1,
\end{align}
where $\mathbb{E}[{C_{\text{S}}}(n,n,1)\mid k]$ is given by Equation \eqref{eq:E_X_k_1D_RBMP-S}, and $\hat{\mathbb{P}}(k)$ is given by Equation \eqref{eq:P_swap_approx}.  
Then, substituting Equations \eqref{eq:E_X_k_1D_RBMP-S} and \eqref{eq:E_Y_1d} into Equation \eqref{eq:beta_delta_RBMP-S} gives:
\begin{align}
    \beta_{\text{S}}(n,1) \approx \frac{\frac{1}{4}\sqrt{\frac{\pi}{2}}\cdot n^{-\frac{1}{2}}(n+1)}{R_{\text{S}}\sum_{k=1}^n k \cdot\hat{\mathbb{P}}(k)} - 1.
\end{align}

For $D \ge 3$, no formulas are available for $\mathbb{E}[\accentset{\ast}{C_{\text{S}}}(n,n,D\ge 3)]$. Nevertheless, the asymptotic values of $\delta_{\text{S}}(\infty,\infty,D)$ reported in Table \ref{tab:deltas} provide benchmarks for their relative magnitudes. Since these values are already very small for all $D \ge 3$, we adopt them directly for finite-size problems as well.
Accordingly, the scaling factor can be approximated as
\begin{align} \label{eq:C_3d}
    \beta_{\text{S}}(n,D)
    \approx  
    \delta_{\text{S}}(\infty,\infty,D)\cdot D^2, \, \, \forall  D \ge 3.
\end{align}

For $D = 2$, neither formulas nor asymptotic estimations are currently available. However, it is reasonable to expect that its scaling behavior lies between that of $D = 1$ and $D = 3$. Accordingly, we approximate it as the average of the values obtained for $D = 1$ and $D = 3$:
\begin{align} \label{eq:C_2d}
    \beta_{\text{S}}(n,2) \approx 
    \frac{1}{2}\left[\beta_{\text{S}}(n,1)+\beta_{\text{S}}(n,3)\right].
\end{align}

Finally, the exponent $d$ in the imbalance discount factor $\left(\frac{m}{n}\right)^d$ can be estimated by fitting the logarithm of \eqref{eq:delta} with the data from simulations via linear regression. For simplicity, here we only consider $d\in \mathbb{Z}$. 
It is found that setting 
$$d = 3$$ provides a good approximation across all available data. 

\section{RBMP-B}
\label{sec:RBMP-B}
We next study the RBMP-B, where the results from RBMP-I are once again employed as a first approximation. We propose an additional boundary correction factor and new estimators to explicitly account for the impacts of the boundary on the $k$-th nearest-neighbor matching distance and probability.

\subsection{Problem Definition}
The key distinct feature in RBMP-B is that the vertices in $U$ and $V$ are independently and uniformly distributed in a unit-volume hyper-ball embedded in a $D$-dimensional L$^p$ metric space, denoted by $\mathcal{B}_{D,p}\subseteq \mathbb{R}^{D}$, where $D\in \mathbb{Z}^+$. 
The subplots in Figure \ref{fig:RBMP_examples} show examples for $\mathcal{B}_{2,1}$, $\mathcal{B}_{2,2}$ and $\mathcal{B}_{3,2}$, respectively. 
The radius of $\mathcal{B}_{D,p}$, denoted as $R_{\text{B}}$, is given by \citep{olver_nist_2010}:
\begin{align}
\label{eq:R_B}
    R_{\text{B}} = \frac{\left(\Gamma(\frac{D}{p}+1)\right)^{\frac{1}{D}}}{2\Gamma(\frac{1}{p}+1)}.
\end{align}
Let $\mathbf{x}_i = (x_{i,1}, \cdots, x_{i,D})$ denote the coordinates of vertex $i \in U \cup V$ in $\mathcal{B}_{D,p}$. 
The cost $c_{u,v}$ on edge $(u,v)$ is defined as the L$^p$ distance between the realized vertices; i.e.,
\begin{align}
c_{u,v} = \left( \sum_{j=1}^D |x_{u,j} - x_{v,j}|^p \right)^{1/p}.
\end{align}
An RBMP-B is now characterized by the vertex set sizes $m$ and $n$, spatial dimension $D$ and distance metric $p$. The CDF of the edge cost, denoted by $F_{C_{\text{B}}}(\cdot)$, also depends on these four parameters.
Our primary interest is in estimating the expectation of the optimal matching cost (or distance), denoted by $\accentset{\ast}{C_{\text{B}}}(m,n,D,p)$, or simply $\accentset{\ast}{C_{\text{B}}}$. 

The estimation of $\mathbb{E}[\accentset{\ast}{C}_{\text{B}}]$ involves the following steps:
\begin{itemize}
    \item[(i)] Estimating the probability for a demand vertex to be matched to its $k$-th nearest neighbor. Similarly to Section \ref{sec:RBMP-S}, we assume i.i.d. edge costs as a first approximation, and use $\accentset{\ast}{\mathbb{P}}(k)$ from the corresponding RBMP-I counterpart.
    \item[(ii)] Estimating the conditionally expected distance of a demand vertex to its $k$-th nearest neighbor, $\mathbb{E}[C_{\text{B}} \mid k]$. It is important to note that the presence of boundaries in RBMP-B introduces even stronger correlations among vertices than in RBMP-S. In RBMP-S, the distance distribution is independent of a vertex's location, whereas in RBMP-B, location may play a non-negligible role; for example, vertices near the center of $\mathcal{B}_{D,p}$ may have a different distance distribution than those near the boundary. To account for this, we will need a new step to estimate the distance conditional on a vertex's position within $\mathcal{B}_{D,p}$, as explained in Section \ref{sec:dist_k_RBMP-B}.
    \item[(iii)] Two levels of correction are required here to relax the i.i.d. assumption in RBMP-I. The first is to apply the spatial correlation correction factor $\delta_{\text{S}}$, as provided in Section \ref{sec:correction_RBMP-S}, which accounts for the correction from RBMP-I to RBMP-S. The second is to apply an additional boundary correction factor, denoted by $\delta_{\text{B}}(m,n,D)$ or $\delta_{\text{B}}$, which accounts for the correction from RBMP-S to RBMP-B, as explained in Section \ref{sec:correction_RBMP-B}.
\end{itemize}

\subsection{Matching Distance Formula}

We next present a set of distance estimation formulas for RBMP-B applicable to various scenarios.
Again, building on Proposition \ref{prop:X_optimal}, we assume $\mathbb{E}[\accentset{\ast}{C_{\text{B}}}]$ can be estimated as follows:
\begin{align} \label{eq:E_Z_opt}
\mathbb{E}[\accentset{\ast}{C}_{\text{B}}] = 
\left( 1 + \delta_{\text{B}} \right)\left( 1 + \delta_{\text{S}} \right)\sum_{k=1}^m\mathbb{E}[C_{\text{B}}\mid k]\cdot \accentset{\ast}{\mathbb{P}}(k).
\end{align}

For the scenario when $n \gg m$, Section \ref{sec:dist_RBMP-S} has already shown that the spatial correlation correction factor $\delta_{\text{S}} \xrightarrow{n\gg m} 0.$ Similarly, when each demand vertex is surrounded by many supply vertices, the effect of vertex location is minimal, and the boundary effects are expected to diminish. As such, 
$$\delta_{\text{B}} \xrightarrow{n\gg m} 0.$$
Moreover, in this case, the distance distribution of vertices located near the center of the hyperball provides a good approximation for that of vertices at other locations.
Thus, $F_{{C}_{\text{B}}}(x)$ converges to a power-law form; i.e.,
$$F_{C_{\text{B}}}(x)
\xrightarrow{n\gg m}
\left( \frac{x}{R_{\text{B}}} \right)^D, \quad x \in [0, R_{\text{B}}].$$   

Given the power-law distribution, and similarly to the analysis in Section \ref{sec:dist_RBMP-S}, we know Proposition \ref{prop:X_greedy_opt} applies. Accordingly, when $n\gg m$, the $M$-th moment of the optimal matching distance, $\mathbb{E}[\accentset{\ast}{C}_{\text{B}}^M]$ can be approximated by the greedy matching distance:
\begin{align} 
    \mathbb{E}[\accentset{\ast}{C}_{\text{B}}^M] 
    &\xrightarrow{n\gg m} 
    \mathbb{E}[\bar{C}_{\text{B}}^M]
    =\frac{R_{\text{B}}^M\cdot\Gamma(n+1)}{\Gamma(n+1+\frac{M}{D})}\sum_{i=1}^m \frac{
    \Gamma(k+\frac{M}{D})}{
    \Gamma(k)}\cdot \bar{\mathbb{P}}(k),
\label{eq:E_Z_greedy}
\end{align}
where $\bar{\mathbb{P}}(k)$ is given by Equation \eqref{eq:P_greedy_approx}.
Specifically, the expected optimal matching distance $\mathbb{E}[\accentset{\ast}{C}_{\text{B}}]$ can be estimated using the $\kappa$-approximation similarly as in Equation \eqref{eq:E_X_greedy_kappa}:
\begin{align} \label{eq:E_Z_greedy_kappa}
\mathbb{E}[\accentset{\ast}{C_{\text{B}}}]\xrightarrow{n\gg m}
\mathbb{E}[\bar{C}_{\text{B}}]
&\approx 
\mathbb{E}_\kappa [\bar{C}_{\text{B}}]\nonumber\\
& = \frac{R_{\text{B}}}{mn^{\frac{1}{D}}} \sum_{i=1 }^{m} 
        \left[ \sum_{k=1}^{i} \left( \frac{i-1}{n} \right)^{k-1} \left( 1- \frac{i-1}{n} \right) 
\left( \mathbf{1}_{k \le \kappa} \cdot \frac{\Gamma(k+\frac{1}{D}) }{\Gamma(k)} + \mathbf{1}_{k > \kappa} \cdot k^{\frac{1}{D}}\right) \right. \nonumber\\
        & +  \left. \left( \frac{i-1}{n} \right)^{i} \left( \mathbf{1}_{i \le \kappa} \cdot \frac{\Gamma(i+\frac{1}{D}) }{\Gamma(i)} + \mathbf{1}_{i > \kappa} \cdot i^{\frac{1}{D}}\right) \right]. 
\end{align}
In addition, when $n\ggg m$, it further converges to the nearest matching distance, $\mathbb{E}[\bar{C}_{\text{B}}(1,n,D,p)]$:
\begin{align}
    \mathbb{E}[\accentset{\ast}{C}_{\text{B}}]
    \xrightarrow{n\ggg m} \mathbb{E}[\bar{C}_{\text{B}}(1,n,D,p)] 
    = R_{\text{B}}\cdot \Gamma(1+\frac{1}{D})\cdot n^{-\frac{1}{D}}, \label{eq:E_Z_nearest}
\end{align}

For general cases with arbitrary values of $m$ and $n$, we propose to use the refined matching probability $\hat{\mathbb{P}}(k)$ to approximate the optimal matching probability $\accentset{\ast}{\mathbb{P}}(k)$, and have:
\begin{align} \label{eq:E_Z_swap}
    \mathbb{E}[\accentset{\ast}{C}_{\text{B}}] 
    \approx \left( 1 + \delta_{\text{B}} \right)\left( 1 + \delta_{\text{S}} \right)
    \sum_{k=1}^m\mathbb{E}[C_{\text{B}}\mid k]\cdot \hat{\mathbb{P}}(k),
\end{align}
where $\hat{\mathbb{P}}(k)$ and $\delta_{\text{S}}$ are given by Equations \eqref{eq:P_swap_approx} and \eqref{eq:delta}, respectively. 
In the following sections, we  provide details for estimating the $\mathbb{E}[C_{\text{B}}\mid k]$ and $\delta_{\text{B}}$, respectively.

\subsection{$k$-th Nearest-neighbor Distance}
\label{sec:dist_k_RBMP-B}
\begin{figure}[ht!]
    \centering
    \begin{subfigure}[b]{0.5\textwidth}
        \centering
        \includegraphics[height=2.2in]{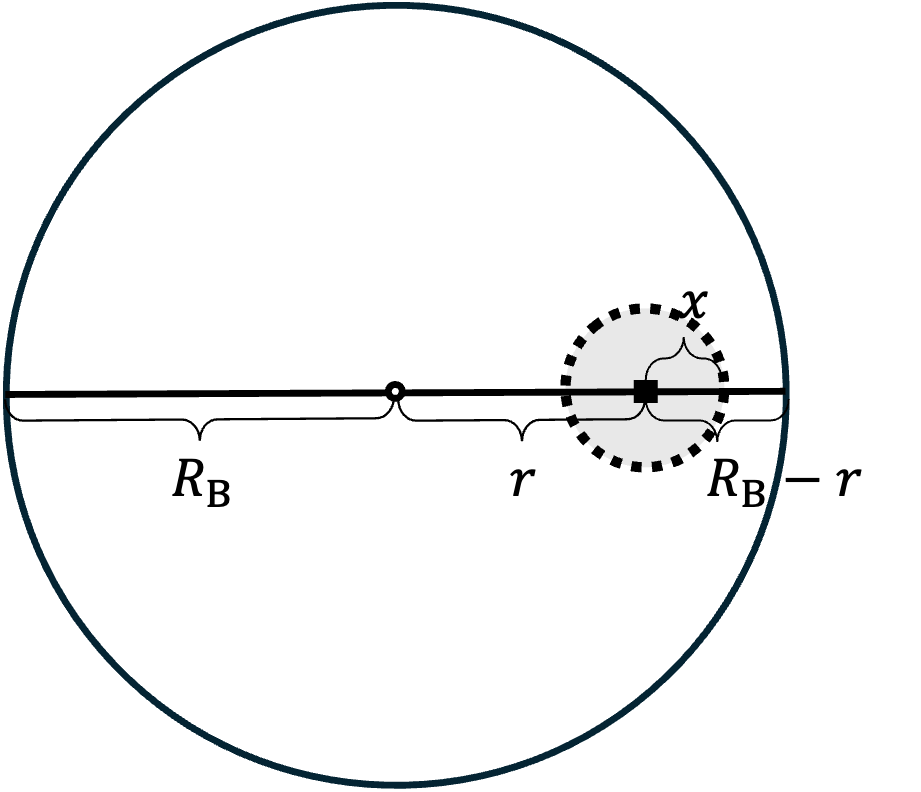}
        \caption{$0\le x \le R_{\text{B}}-r$.}
        \label{fig:intersection_inner}
    \end{subfigure}%
    ~ 
    \begin{subfigure}[b]{0.5\textwidth}
        \centering
        \includegraphics[height=2.2in]{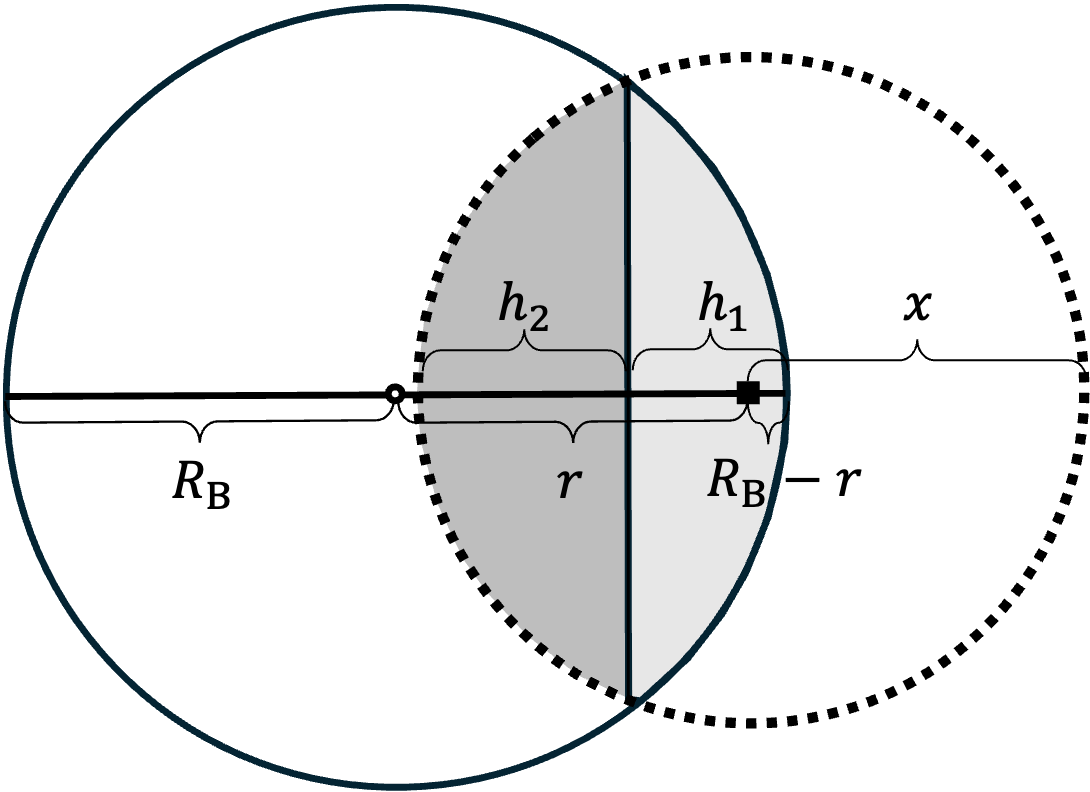}
        \caption{$R_{\text{B}}-r<x\le R_{\text{B}}+r$.}
        \label{fig:intersection_outer}
    \end{subfigure}
    \caption{Intersection between two hyper-balls in $2$-dimensional L$^2$ metric space.}
    \label{fig:intersection_volume}
\end{figure}

The key to deriving $\mathbb{E}[C_{\text{B}} \mid k]$ is to account for the location of the random demand vertex, indicated by its distance to the center of $\mathcal{B}_{D,p}$, denoted by a random variable $r$. Figure \ref{fig:intersection_volume} illustrates an example of a randomly sampled demand vertex in $\mathcal{B}_{2,2}$, indicated by the square marker. The derivation contains three steps:
\begin{itemize}
\item[(i)] Conditional on the demand vertex's location $r$, we first derive the CDF of its distance from a randomly selected supply vertex, ${C}_{\text{B}}(r)$. 
\item[(ii)] The distribution of the distance from that demand vertex to its $k$-th nearest neighbor, ${C}_{\text{B}}(r) \mid k$, can then be obtained via order statistics.
\item[(iii)] Apply the law of total probability over all possible demand vertex locations (or values of $r\in [0, R_{\text{B}}]$) to compute $\mathbb{E}[C_{\text{B}} \mid k]$.
\end{itemize}

We first look at step (i) and focus on deriving the CDF of ${{C}_{\text{B}}}(r)$. 
For a demand vertex located at distance $r$ from the center, and given that the supply vertices are uniformly distributed, the probability that ${C}_{\text{B}}(r)$ is less than a threshold $x$, i.e., $F_{{C}_{\text{B}}(r)}(x)$, equals the volume of the intersection between $\mathcal{B}_{D,p}$ (the solid-line circles in Figure \ref{fig:intersection_volume}) and another hyper-ball of radius $x$ centered at that demand vertex (the dotted-line circles in Figure \ref{fig:intersection_volume}).
This volume depends on the values of $r$ and $x$, and two cases may arise. 
When $0\le x \le R_{\text{B}}-r$, the hyper-ball of radius $x$ lies entirely within $\mathcal{B}_{D,p}$, as illustrated by the light-gray shaded region in Figur ~\ref{fig:intersection_inner}. In this case, the intersection volume is simply $(x / R_{\text{B}})^D$.
When $R_{\text{B}}-r < x \le R_{\text{B}}+r$, the 
intersection consists of two hyper-spherical caps (one from each hyper-ball), as illustrated by the light-gray and dark-gray shaded regions in Figure \ref{fig:intersection_outer}, and its volume depends on $D$ and $p$. 
%
The heights of these two caps (as shown in Figure \ref{fig:intersection_outer}) can be derived from simple geometry:
\begin{align}
    h_1 = \frac{x^2-(R_{\text{B}}-r)^2}{2r}, \quad h_2 = x+R_{\text{B}}-r-h_1.
\end{align}
Let $\mathcal{V}(R, h, D, p)$ denote the volume of a hyper-spherical cap in a $D$-dimensional L$^p$ metric space, given radius $R$ and cap height $h$. 
When $p = 2$, \citet{li_concise_2010} provides a general formula to compute $\mathcal{V}(R, h, D, p)$ for an arbitrary dimension $D$: 
\begin{align} \label{eq:V_general}
\mathcal{V}(R, h, D, p) =
\begin{cases}
\frac{\left(2R\cdot\Gamma\left(\frac{1}{p}+1\right)\right)^D}{2\Gamma\left(\frac{D}{p}+1\right)} I_{\frac{2Rh - h^2}{R^2}}\left(\frac{D+1}{2},\frac{1}{2}\right), & \text{if } h \le R, \\
\frac{\left(2R\cdot\Gamma\left(\frac{1}{p}+1\right)\right)^D}{2\Gamma\left(\frac{D}{p}+1\right)} \left[ 2 - I_{\frac{2R(2R - h) - (2R - h)^2}{R^2}}\left(\frac{D+1}{2},\frac{1}{2}\right) \right], & \text{if } h > R.
\end{cases}
\end{align}
For general cases with arbitrary values of $p$, however, no closed-form formula is available for $\mathcal{V}(R, h, D, p)$. 
We propose to use formula \eqref{eq:V_general} as an approximation. 
The CDF of ${C}_{\text{B}}(r)$ in step (i) can be expressed as follows:
\begin{align}
F_{{C}_{\text{B}}(r)}(x) =
\begin{cases}
    \left( \frac{x}{R_{\text{B}}} \right)^D,\quad 0\le x\le R_{\text{B}}-r,\\
    \mathcal{V}(R_{\text{B}},h_1, D, p)+\mathcal{V}(x,h_2, D, p), \quad R_{\text{B}}-r < x\le R_{\text{B}}+r.
\end{cases}
\end{align}
The CDF of ${C}_{\text{B}}(r)\mid k$ can be computed via order statistics in step (ii), using Equation \eqref{eq:F_X_k}:
\begin{align}
F_{{C}_{\text{B}}(r)\mid k}(x) =
\begin{cases}
    I_{\left(\frac{x}{R_{\text{B}}}\right)^D}(k,n-k+1), \quad 0\le x\le R_{\text{B}}-r,\\
    I_{\mathcal{V}(R_{\text{B}},h_1, D, p)+\mathcal{V}(x,h_2, D, p)}(k,n-k+1), \quad R_{\text{B}}-r < x\le R_{\text{B}}+r.
\end{cases}
\end{align}

Finally, we are ready to derive $\mathbb{E}[C_{\text{B}}\mid k]$ in step (iii). 
Since demand vertices are uniformly distributed within $\mathcal{B}_{D,p}$, the distance $r$ from a randomly selected demand vertex to the center follows a power-law distribution with p.d.f.: 
$$f_r(r)=\frac{D}{R_{\text{B}}^D}r^{D-1}, \quad r \in [0,R_{\text{B}}].$$ 
By the law of total probability, $\mathbb{E}[C_{\text{B}}\mid k]$ is then obtained by integrating the conditional expectation 
$$\mathbb{E}[{C}_{\text{B}}(r)\mid k]=\int_0^{R_{\text{B}}-r}x\text{ d}F_{{C}_{\text{B}}(r)\mid k}(x)+\int_{R_{\text{B}}-r}^{R_{\text{B}}+r}x\text{ d}F_{{C}_{\text{B}}(r)\mid k}(x)$$
over all possible values of $r$, as follows:
\begin{align} \label{eq:E_Z_k_general}
\begin{split}
    \mathbb{E}[C_{\text{B}}\mid k]
    &= \int_0^{R_{\text{B}}}f_r(r)\cdot\mathbb{E}[{C}_{\text{B}}(r)\mid k]  \text{ d}r\\
    &= \frac{D}{R_{\text{B}}^D}\int_0^{R_{\text{B}}} r^{D-1}
    \left[ 
    \int_0^{R_{\text{B}}-r}x\text{ d}I_{\left(\frac{x}{R_{\text{B}}}\right)^D}(k,n-k+1)
    \right.\\ 
    & \qquad\qquad\qquad\qquad
    +\left.\int_{R_{\text{B}}-r}^{R_{\text{B}}+r}x\text{ d}I_{\mathcal{V}(R_{\text{B}},h_1, D, p)+\mathcal{V}(x,h_2, D, p)}(k,n-k+1) 
    \right] \text{ d}r.
\end{split}
\end{align}

\subsection{Correction Factor Estimation}
\label{sec:correction_RBMP-B}
In light of Section \ref{sec:correction_RBMP-S}, we propose that $\delta_{\text{B}}(m,n,D)$ should have a functional form similar to that of Equation \eqref{eq:delta} for $\delta_{\text{S}}(m,n,D)$:
\begin{align}
    \delta_{\text{B}}(m,n,D) \approx \beta_{\text{B}}(n,D) \left(\frac{m}{n}\right)^{3} \frac{1}{D^2}.
\end{align}
The key difference is that we now require an updated scaling factor $\beta_{\text{B}}(n,D)$ for balanced problems, which should satisfy the following properties. In the special case where $D = 1$, $$\beta_{\text{B}}(n,1) = \delta_{\text{B}}(n,n,1),$$ which can be directly estimated by comparing the expected matching costs of RBMP-B and its RBMP-S counterpart.
Again, as a baseline, a closed-form formula for estimating the expected optimal matching cost
in a balanced 1-dimensional RBMP-B, $\mathbb{E}[\accentset{\ast}{C}_{\text{B}}(n,n,1,p)]$, is provided in \citet{zhai_average_2024} for $n\gg 1$. 
\begin{align}
    \mathbb{E}[\accentset{\ast}{C}_{\text{B}}(n,n,1,p)]
    \xrightarrow{n\gg 1} 
    \frac{\sqrt{\pi}}{4}\cdot n^{-\frac{1}{2}}.
\end{align}
Then, given the expected optimal matching cost
in a balanced 1-dimensional RBMP-S in Equation \eqref{eq:E_Y_1d}, $\beta_{\text{B}}(n,1)$ can be approximately estimated as:
\begin{align} \label{eq:H_1d}
\beta_{\text{B}}(n,1) \approx \frac{\mathbb{E}[\accentset{\ast}{C}_{\text{B}}(n,n,1,p)]}{\mathbb{E}[\accentset{\ast}{C}_{\text{S}}(n,n,1)]} - 1 = \sqrt{2} - 1.
\end{align}
For $D \ge 3$, we adopt the same scaling factor used for RBMP-S; see Equation \eqref{eq:C_3d}:
\begin{align} \label{eq:H_3d}
\beta_{\text{B}}(n,D\ge 3) \approx \beta_{\text{S}}(n,D\ge 3) \approx \delta_{\text{S}}(\infty,\infty,D \ge 3) \cdot D^2.
\end{align}
For $D = 2$, we estimate the scaling factor using the same approach as for Equation \eqref{eq:C_2d}, by averaging the values for $D = 1$ and $D = 3$:
\begin{align} \label{eq:H_2d}
\beta_{\text{B}}(n,2) \approx \frac{1}{2}[\beta_{\text{B}}(n,1) + \beta_{\text{B}}(n,3)].
\end{align}

\section{Numerical Results} \label{sec: experiment}

In this section, we present simulation results to verify the accuracy of the proposed formulas for the three problem variants: RBMP-I, RBMP-S, and RBMP-B. 
A set of Monte-Carlo simulations are conducted across various parameter combinations, including the vertex numbers $n$ and $m$, spatial dimension $D$, and distance metric $p$. 
In each simulation run, we generate $10^3$ RBMP instances. 
Each realized instance is solved by the state-of-art modified Jonker-Volgenant algorithm \citep{modified_JV}. For each parameter combination, we record the optimal matching cost (or distance) for each demand vertex, as well as the sample mean and standard deviation over all vertices across the $10^3$ instances. 

\subsection{Verification of RBMP-I}
\label{sec: sim_dist}

We begin with RBMP-I and and examine the accuracy of the proposed formulas for $\mathbb{E}[\accentset{\ast}{C}^M(m,n)]$: 
(i) Equation \eqref{eq:E_X_opt_general}, which approximates the expected optimal cost for general $m$ and $n$ values; 
(ii) Equation \eqref{eq:E_X_opt_special}, which serves as an upper bound of that expected optimal cost; and 
(iii) Equation \eqref{eq:convergence}, which serves as a lower bound. Both bounds converge to the expected optimal cost in the special case when $n \gg m$, as identified in Proposition \ref{prop:X_greedy_opt}. 
Here we set $D \in \{1,2,3\}$, 
where the edge cost distribution is given by $F_{C}(x)=\left(\tfrac{x}{R}\right)^D$, following the power-law form introduced in Section \ref{sec:power-law}. 
The maximum radius $R$ is set to $R_{\text{B}}$ as defined in Equation \eqref{eq:R_B}, with $p=2$. 
We further set $m \in \{10, 100\}$, and $n$ ranging from $m$ to a notably larger number, e.g. $3m$. 

\begin{figure}[ht!] 
\centering 
\includegraphics[width=\textwidth]{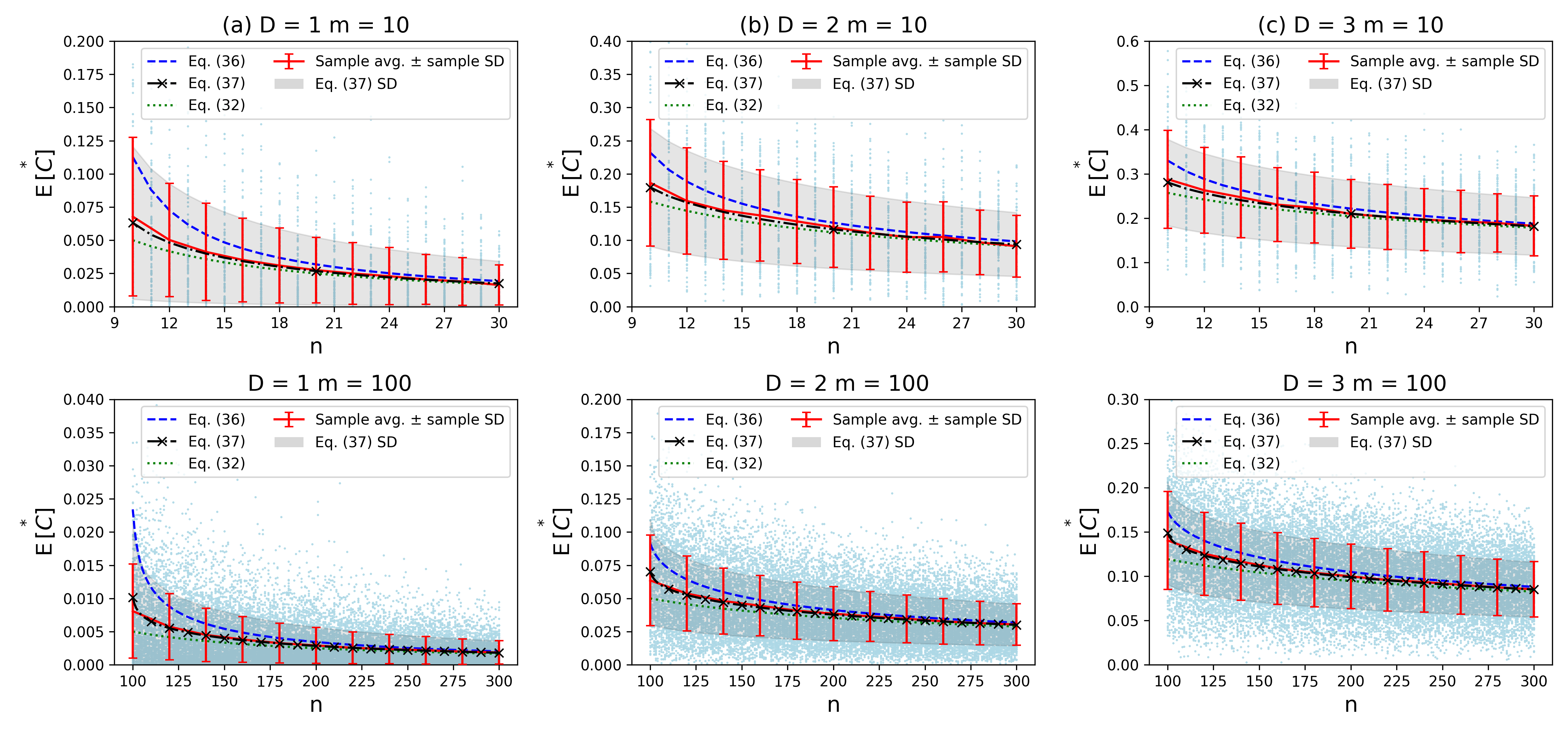}
\caption{
RBMP-I estimators vs. simulations under different dimensions.}
\label{fig: RBMP-I}
\end{figure}

Figure \ref{fig: RBMP-I} shows a comparison of the estimated expected costs (i.e., by setting $M=1$ in the formulas) with the corresponding simulation results.  
The optimal matching cost for each demand vertex from each simulated RBMP-I instance is represented by a light-blue dot. The corresponding sample mean across all vertices is represented by the red solid line. 
The estimates  from Equations \eqref{eq:E_X_opt_special}, \eqref{eq:E_X_opt_general}, and \eqref{eq:convergence} are represented by the blue dashed lines, the black dash-dotted lines with cross markers, and the green dotted lines, respectively. 

It can be observed that the estimations from Equation \eqref{eq:E_X_opt_general} closely match with the simulated sample means across all parameter combinations $(m,n,D)$. 
When $m = 10$, Equation \eqref{eq:E_X_opt_general} has an average relative error of 3.25\%, 2.13\%, and 1.41\% for $D=1,2,3$, respectively, across all $n$ values. When $m = 100$, Equation \eqref{eq:E_X_opt_general} has an average relative error of 2.01\%, 1.54\%, and 1.00\% for $D=1,2,3$, respectively. 
This indicates that our proposed formula Equation \eqref{eq:E_X_opt_general} 
is a reasonably good estimator for RBMP-I in general parameter settings.
Specifically, when $n \gg m$, it can be seen from Figure \ref{fig: RBMP-I} that Equations \eqref{eq:convergence}, \eqref{eq:E_X_opt_general} and \eqref{eq:E_X_opt_special} start to converge. 
In these cases, the simpler Equations \eqref{eq:convergence} and \eqref{eq:E_X_opt_special} can be used to estimate the optimal matching cost. For example, when $D=2,m=100$, Equations \eqref{eq:convergence} and \eqref{eq:E_X_opt_special} have average relative errors of 6.04\% and 5.85\% across all $n\geq 2m$, respectively. 
This also indicates that ``greedy matching" (as discussed in Section \ref{sec:greedy}) can already yield near-optimal solutions once $n\ge 2m$, with most demand vertices being matched to their nearest neighbors. 

In addition to the expected cost, 
Figure \ref{fig: RBMP-I} also compares the estimated standard deviation (SD) between the formulas and simulations. 
For each parameter combination, the formula-based SD is be obtained from the second moment (i.e., $M=2$) in Equation \eqref{eq:E_X_opt_general}; i.e., SD = $\sqrt{\mathbb{E}[\accentset{\ast}{C}^2] - \mathbb{E}[\accentset{\ast}{C}]^2}$.
The formula estimates of the interval $[\mathbb{E}[\accentset{\ast}{C}] -\text{SD}, \mathbb{E}[\accentset{\ast}{C}] + \text{SD}]$ are represented by the gray shaded areas, 
while the corresponding sample estimates from the simulations are shown by the red error bars. 
It is observed that the formulas can provide very good estimates of the SD as well. 
When $m = 10$, Equation \eqref{eq:E_X_opt_general} predicts the standard deviation with average relative errors of 3.65\%, 2.13\%, and 4.01\% for $D=1,2,3$, respectively; when $m = 100$, the relative errors are 7.21\%, 1.29\%, and 2.92\% for $D=1,2,3$, respectively. 

\begin{figure}[ht!] 
\centering 
\includegraphics[width=\textwidth]{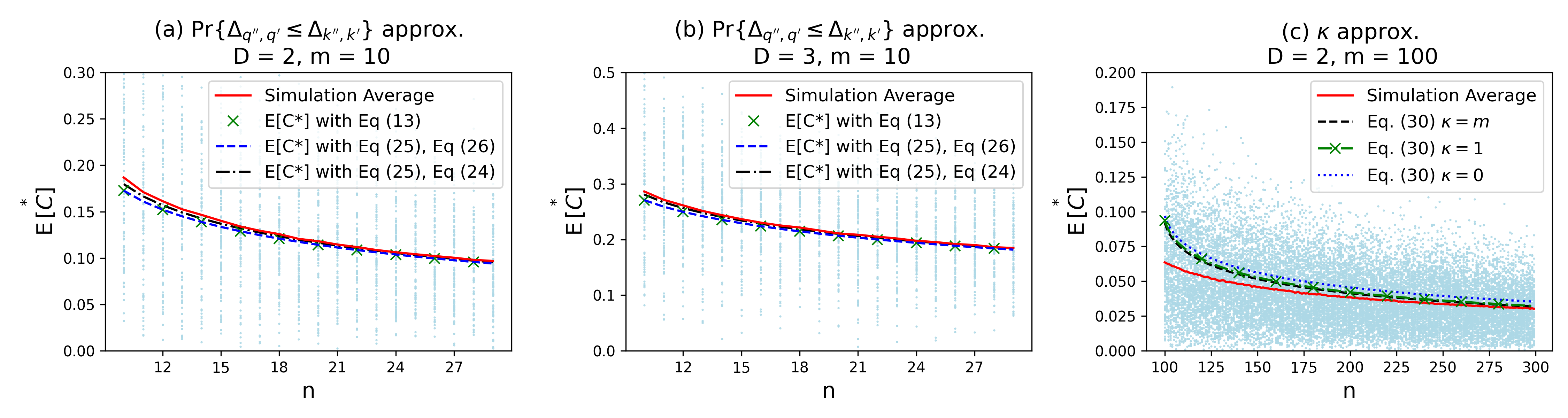}
\caption{
Impacts of RBMP-I approximations.}
\label{fig: RBMP-I_approx}
\end{figure}

Finally, Figure \ref{fig: RBMP-I_approx} verifies the impacts of several key approximations during the derivation, including the normal approximation used for  $\text{Pr}\{\Delta_{q'',q'}\leq \Delta_{k'',k'}\}$ in Equation \eqref{eq:P_normal} and the $\kappa$-approximation used for $\mathbb{E}[\bar{C}]$ in Equation \eqref{eq:E_X_greedy_kappa}. 

Recall from Section \ref{sec:P_approx} that, for computational simplicity, Equation \eqref{eq:E_X_opt_general} uses two approximations: assuming a normal distribution in computing $\text{Pr}\{\Delta_{q'',q'}\leq \Delta_{k'',k'}\}$, and applying the mean and variance computed for $D=1$ to all other dimensions (e.g., $D=2,3$), as given by Equations \eqref{eq:P_normal} and \eqref{eq:mu_sigma}.
Figures \ref{fig: RBMP-I_approx} (a)-(b) show formula predictions of $\mathbb{E}[\accentset{\ast}{C}]$ when 
$\text{Pr}\{\Delta_{q'',q'}\leq \Delta_{k'',k'}\}$ is computed using either (i) the exact CDF formula in Equation \eqref{eq:P_swap_one}, shown as green cross markers, or (ii) the normal approximation but with the exact mean and variance for the respective $D$ value, based on Equations \eqref{eq:P_normal} and \eqref{eq:E_Delta}, shown as blue dash lines, respectively. 
Since both formulas involve integrals that are computationally cumbersome for $D > 1$, we only showcase the results for smaller-scale parameter settings with $m=10$. 
It can be observed that the estimates from the two formulas closely align with the results obtained when directly applying Equations \eqref{eq:P_normal} and \eqref{eq:mu_sigma} (shown as the black dash-dotted lines).
This indicates that the normal approximation used in the proposed estimator, Equation \eqref{eq:E_X_opt_general}, performs well and is practically effective.

Figure \ref{fig: RBMP-I_approx}(c) shows the $\kappa$-approximation estimation for the greedy matching cost $\mathbb{E}[\bar{C}]$ from Equation \eqref{eq:E_X_greedy_kappa}. As discussed in Section \ref{sec:greedy_cost}, parameter $\kappa$ indicates the number of dominating gamma function terms that are not computed by Stirling's approximation.
Here we use blue dotted lines, green dash-dotted lines with markers, and black dashed lines to represent $\kappa=0, 1, m$, respectively. It is observed that the difference between those from $\kappa = 1$ and $\kappa = m$ is minimal, and they both align closely with the simulated results as $n \gg m$. 
In addition, setting $\kappa=0$ here (which is computationally very appealing) does not sacrifice too much accuracy as compared to $\kappa=1$.
This suggests that choosing $\kappa = 0$ or $1$ can already provide a sufficiently accurate estimation, as compared to the computationally more challenging choice of $\kappa = m$. This observation is particularly useful for applying the approximate formula to RBMP-B later in Section \ref{sec: application}.

\subsection{Verification of RBMP-S}
Next, we examine the accuracy of the proposed formulas for $\mathbb{E}[\accentset{\ast}{C}_{\text{S}}(m,n,D)]$ in RBMP-S: (i) Equation \eqref{eq:E_Y_swap} with the spatial correlation correction term, which approximates the optimal distance for general $m$ and $n$ values; and (ii) Equations \eqref{eq:E_Y_greedy} and \eqref{eq:convergence_S}, which do not include the correction term but are applicable when $n \gg m$.
Here we set $D \in \{1,2\}$, as most transportation problems in the real world are related to either $D=1$ (such as a ring road) or $D=2$ (such as surface of the globe). Other parameter values remain the same; e.g., $m\in \{10, 100\}$, and $n$ ranges from $m$ to $3m$.  

\begin{figure}[ht!] 
\centering 
\includegraphics[width=0.8\textwidth]{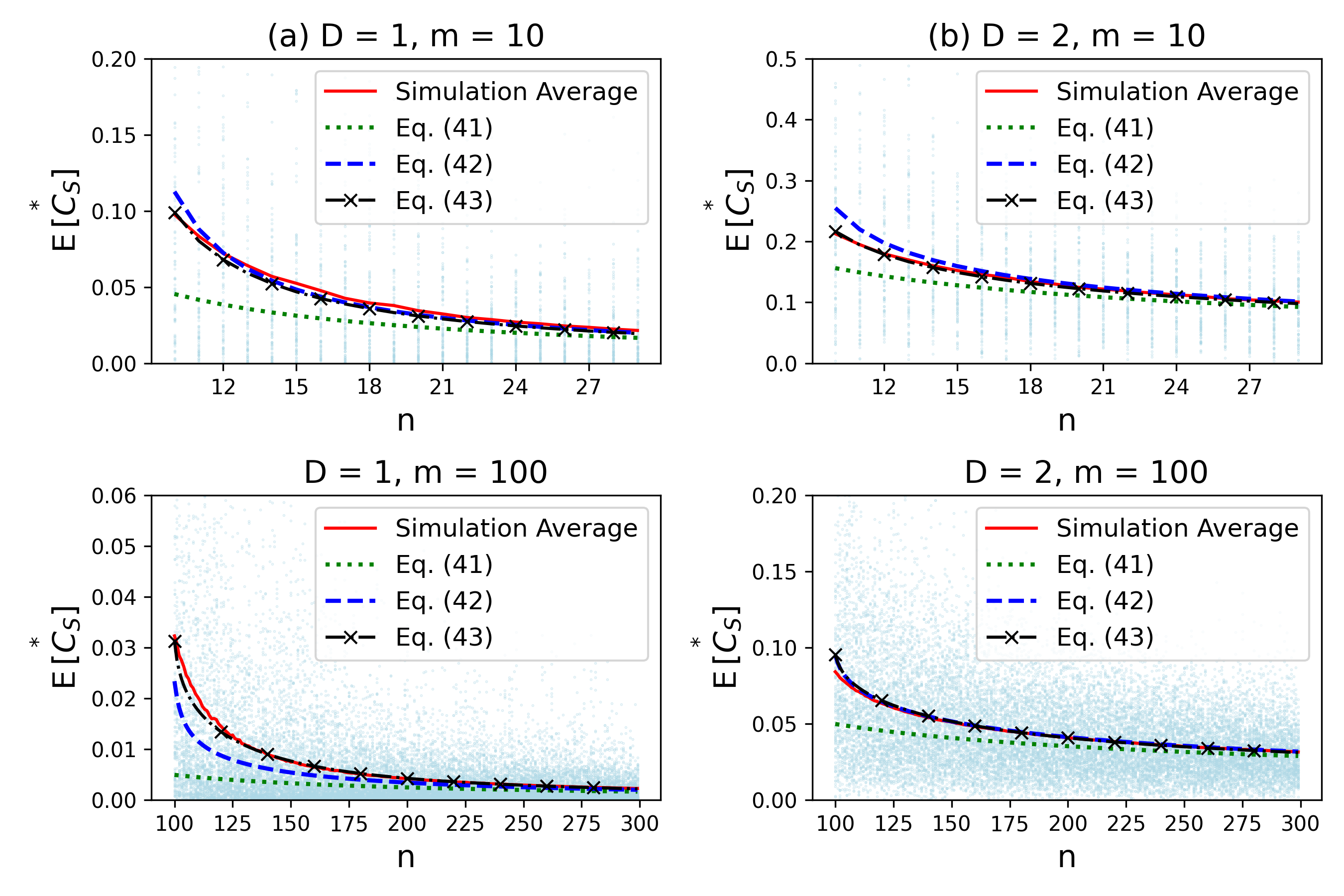}
\caption{
RBMP-S estimators vs. simulations under different dimensions.}
\label{fig: RBMP-S}
\end{figure}

Figure \ref{fig: RBMP-S} compares the Monte-Carlo simulation results with the formula estimates. 
The optimal matching cost for each demand vertex from each
simulated RBMP-S instance is represented by a light-blue dot. The corresponding sample mean across all vertices is represented by the red solid line. 
The estimates from Equations \eqref{eq:convergence_S}, \eqref{eq:E_Y_greedy},  and \eqref{eq:E_Y_swap} are represented by the blue dashed line, green dotted line, and black dashed line with cross markers, respectively. 

It can be observed that the estimations from Equation \eqref{eq:E_Y_swap} closely match with 
the sample means across all parameter combinations $(m, n, D)$. 
In particular, when $m=10$, Equation \eqref{eq:E_Y_swap} has average relative errors 8.83\% and 2.14\% for $D=1,2$, across all $n$ values, respectively; when $m=100$, the average relative errors are 2.69\% and 1.39\% for $D=1,2$, respectively. 
This indicates that Equation \eqref{eq:E_Y_swap} can serve as a reasonably good 
estimator for RBMP-S under general parameter settings.
Moreover, when $n \gg m$, similar to RBMP-I, the estimates from approximate formulas \eqref{eq:E_Y_greedy} and \eqref{eq:convergence_S} also converge well to the expected optimal distance.
This suggests that spatial correlation becomes negligible in such settings, and hence these simpler formulas can also be used effectively. 

However, while Equation \eqref{eq:convergence_S} remains a valid lower bound, we must note that Equation \eqref{eq:E_Y_greedy} no longer serves strictly as an upper bound; e.g., the estimates fall slightly below the corresponding simulated values for some cases. This indicates the importance of including the correction term to address the impact of spatial correlation. 
In particular, when the problem becomes more balanced (i.e., $n \approx m$), when $D$ is smaller, and when $m$ is larger, the effect of spatial correlation becomes more significant. 
This is expected, since the spatial correlation primarily arises from ``unavoidability'' of competitions among demand vertices --- the effect is stronger with fewer supply vertices, in lower dimensions (especially when $D=1$), and with more demand vertices competitors.

Another interesting observation is that, under certain parameter settings (e.g., when $D$ and $m$ is larger), Equation \eqref{eq:E_Y_greedy} can serve as an effective estimator not only for highly unbalanced cases ($n \gg m$) but also for nearly balanced cases $(n \approx m$).
For example, when $m=100$ and $D=2$, Equation \eqref{eq:E_Y_greedy} achieves an average relative error of only 1.33\% across all $n$ values.
A similar pattern will also be observed for RBMP-B in the next section. These findings suggest that the simpler greedy formula, even without correction, nevertheless 
provides suitable approximate estimates of the expected optimal distance for practical use.

\subsection{Verification of RBMP-B}
\label{sec:experiment_B}
\begin{figure}[ht!] 
\centering 
\includegraphics[width=\textwidth]{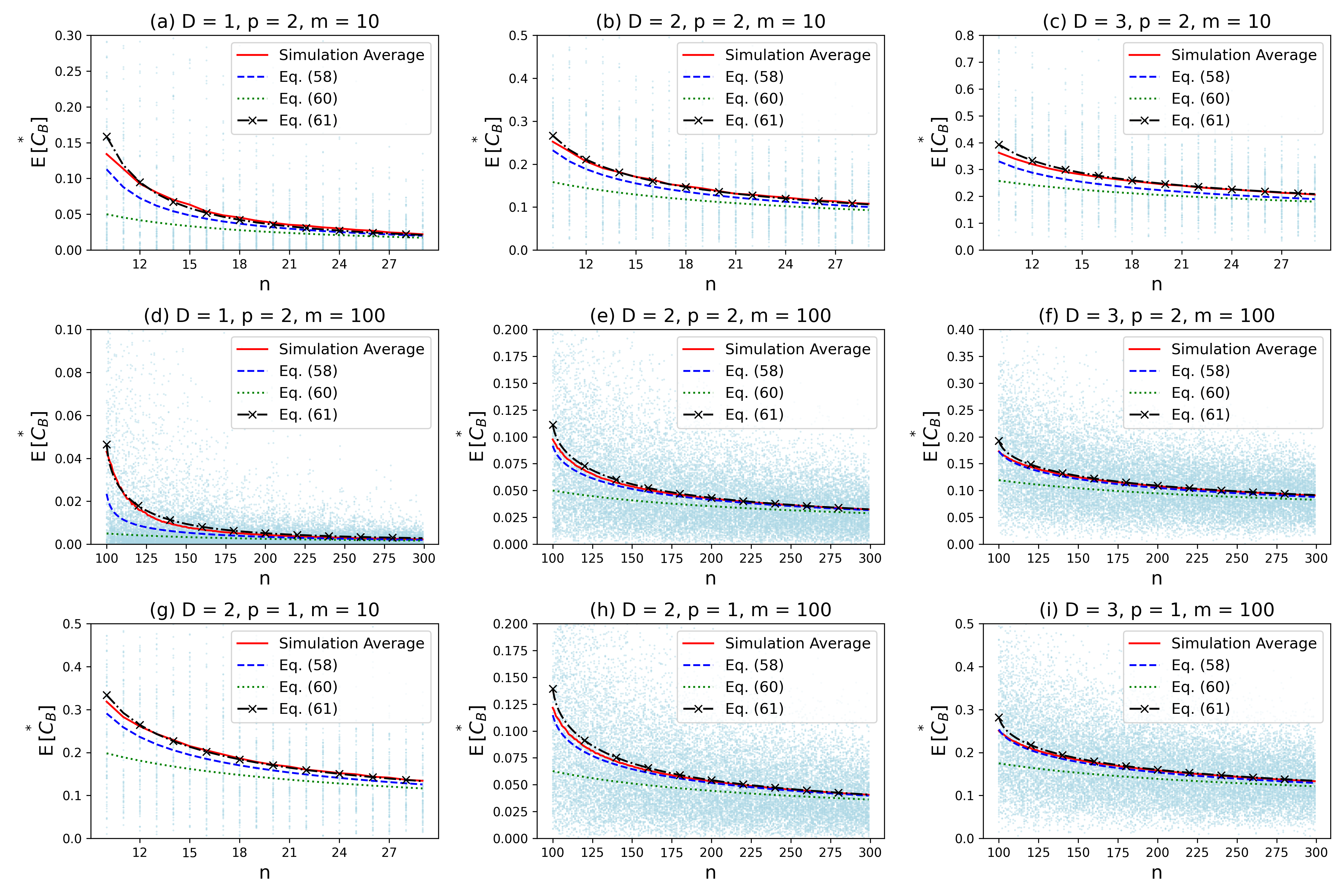}
\caption{
RBMP-B estimators vs. simulations under various $D,p$ values.}
\label{fig: RBMP-B}
\end{figure}

Finally, we examine the accuracy of the proposed formulas for $\mathbb{E}[\accentset{\ast}{C}_{\text{B}}(m,n,D,p)]$ in RBMP-B: (i) Equation \eqref{eq:E_Z_swap} that works for general $m$ and $n$ values; and (ii) Equations \eqref{eq:E_Z_greedy} and \eqref{eq:E_Z_nearest} that are applicable when $n \gg m$.
Here we set $D\in \{1,2,3\}$ and $p\in \{1,2\}$ to address the most commonly encountered real-world setting: one to three-dimensional spaces with Manhattan or Euclidean distance metric. All other parameter choices remain the same. 
Figure \ref{fig: RBMP-B} summarizes the results, where the estimations by Equations \eqref{eq:E_Z_greedy}, \eqref{eq:E_Z_nearest}, and \eqref{eq:E_Z_swap} are represented by the blue dashed lines, green dotted lines, and black dashed lines with cross markers, respectively. 

The overall patterns are similar to what we found in RBMP-S. In general, the estimations from Equation \eqref{eq:E_Z_swap} closely match with the simulation averages across all tested combinations of $(m,n,D,p)$.
For $p=2$, when $m=10$, Equation \eqref{eq:E_Z_swap} has average relative errors of 6.90\%, 1.8\%, and 1.94\% for $D=1,2,3$, respectively; 
when $m=100$, the average relative errors are 6.18\%, 1.64\%, and 0.81\%, respectively.
For $p=1$, when $m=10$, Equation \eqref{eq:E_Z_swap} has average relative errors of 1.31\% and 1.92\% for $D=2,3$, respectively; 
when $m=100$, the average relative errors are 2.86\% and 1.75\%, respectively. 
Recall that in Section \ref{sec:dist_k_RBMP-B}, when deriving the volume of the hyper-spherical cap: $\mathcal{V}(R, h, D, p)$ as given by Equation \eqref{eq:V_general}, we approximate the result for general $p$ values using the formula for $p=2$. 
These results verify the effectiveness of such an approximation. 
When $n \gg m$, Equations \eqref{eq:E_Z_greedy} and \eqref{eq:E_Z_nearest} converge to the expected optimal distance as the correlation becomes negligible.  
When $n\approx m$, Equation \eqref{eq:E_Z_greedy} no longer serve as a valid upper bound due to the absence of the correction term. The correlation's effects, similar to those observed in RBMP-S, become more significant when $D$ is smaller and $m$ is larger. 
Again, at certain parameter combinations, such as when $p=2$ and $m = 100$, Equation \eqref{eq:E_Z_greedy} provides very accurate estimations; e.g., the average relative errors are 3.59\% and 2.82\% over all $n$ values for $D=2, 3$, respectively. 

These results suggest that Equation \eqref{eq:E_Z_swap} serves as an effective distance estimator for RBMP-B under general parameter settings of $(m,n,D,p)$. At the special case when $n\gg m$, the simpler formulas \eqref{eq:E_Z_greedy} and \eqref{eq:E_Z_nearest} are also effective. Moreover, under certain conditions such as $D=2$ and $D=3$, which corresponds to many real-world transportation problems (e.g., ride-hailing or drone delivery), Equation \eqref{eq:E_Z_greedy} can be used as a practically good estimator. To further simplify the computation for Equation \eqref{eq:E_Z_greedy}, the $\kappa$-approximations as given by Equation \eqref{eq:E_Z_greedy_kappa} can be applied by setting $\kappa=0$ or $1$. Their effectiveness has already been verified in RBMP-I (see Figure \ref{fig: RBMP-I_approx}).


\section{Insights and Applications to Mobility Systems} \label{sec: application}

In this section, we show an application of our proposed formulas to estimate the expected optimal matching distance for RBMP-B in the context of 2-dimensional mobility problems (i.e., $D=2$) under a general L$^p$ distance metric. Using the taxi systems as an example, we first show that the convergence property of $\mathbb{E}[\accentset{\ast}{C}_{\text{B}}(m,n,2,p)]$, as revealed in Proposition \ref{prop:X_greedy_opt}, can be used to provide a theoretical explanation (or justification) on the conditions under which the empirically assumed taxi meeting function in the Cobb-Douglas form \citep{yang_equilibria_2010} are suitable, and how some of its parameter values should be set. Then, we demonstrate how our closed-form formulas can be integrated into simple optimization models to design the best operational strategies (e.g., demand pooling interval) for e-hailing taxi systems. The effectiveness of these optimized strategies is verified through two simulation programs under a variety of system settings. 

\subsection{Theoretical Explanation on the Cobb-Douglas Meeting Function} 
\label{sec: cobb}

The Cobb-Douglas meeting function proposed by \citet{yang_equilibria_2010} is widely applied to the analysis of customer-vehicle matching process in taxi systems. It was originally developed for the street-hailing service, where customers must be picked up by vehicles at predetermined locations, but was later extended to the e-hailing services \citep{zha_economic_2016}. 
The model says that the steady-state meeting (or ``pickup") rate between $m$ randomly located waiting customers and $n$ randomly located idle vehicles, denoted $\mathcal{M}(m,n)$ [tu$^{-1}$], 
is given by the following matching function in the Cobb-Douglas form: 
\begin{align*}
\mathcal{M}(m,n) = \alpha_0 \cdot m^{\alpha_1} \cdot n^{\alpha_2},
\end{align*}
where $\alpha_0$ [tu$^{-1}$] is a scaling parameter; $\alpha_1 \in (0,1]$ and $\alpha_2 \in (0,1]$ are unitless elasticities with respect to $m$ and $n$, respectively.
From the Little's formula, one can derive the expected customer waiting time as follows:
\begin{align} \label{eq:cobb_t}
\frac{m}{\mathcal{M}(m,n)} = \alpha_0^{-1} \cdot m^{1-\alpha_1} \cdot n^{-\alpha_2}.
\end{align}
Multiplying this waiting time by the vehicle speed gives the expected matching distance $\mathbb{E}[\accentset{\ast}{C}_{\text{B}}(m,n,2,p)]$. If proper units du and tu are chosen, the value of the area size and vehicle speed of any given service region can both have a value of $1$, and hence Equation \eqref{eq:cobb_t} directly gives the value of $\mathbb{E}[\accentset{\ast}{C}_{\text{B}}(m,n,2,p)]$.

A series of studies have tried to explore the values of parameters $\alpha_0, \alpha_1, \alpha_2$, either theoretical or empirical, under different settings. For instance, \citet{yang_equilibrium_2011} and \citet{zha_economic_2016} argued qualitative that, based on microeconomics theory, 
the value of $\alpha_1+\alpha_2$ should satisfy $>1, =1$ or $<1$ based on whether the market exhibits an increasing, constant or decreasing return to scale, 
while $\alpha_0$ should depend on the service region's characteristics (e.g., area size and vehicle speed). 
In their studies, they simply tested a range of $\alpha_0$ values, and assumed that $\alpha_1 = \alpha_2$ (``for symmetry") with varying values from 0.25 to 0.8. Table \ref{tab_cobb} gives a summary of their assumed parameter values.
These values were chosen purely to demonstrate, in theory, how the market's economic properties (e.g.,  service quality and market profitability) could vary under different parameter combinations; no empirical analyses were conducted to validate these choices.

Other studies have tried to calibrate these parameters based on simulated or empirically observed data, normally via linear regression. 
For example, \citet{yang_returns_2014} was the first to use empirical data, collected from Hong Kong's street-hailing service statistics in years 1986-2009, to estimate the parameters. 
\citet{zhang_efficiency_2019} later used e-hailing service data collected from Shenzhen in year 2016, and parameter values were estimated with or without the assumption that $\alpha_1 = \alpha_2$ --- Significant differences were found in the estimated parameter values (as shown in Table \ref{tab_cobb}), but this ``symmetry'' assumption did not notably affect the goodness-of-fit. 
Furthermore, \citet{wang_-demand_2022} calibrated the model by simulating customer and vehicle data in two hypothetical square regions (with $p=1, 2$) and one real-world region in New York City.
The parameter values obtained from the two hypothetical square regions are similar to each other (i.e., the value of $p$ did not make any major difference), but those from the New York City case differ notably. 

\begin{table}[ht] 
\caption{Parameter values in the literature.}
\centering
\label{tab_cobb}
\resizebox{\textwidth}{!}{%
\begin{tabular}{|cc|c|c|c|c|c|}
\hline
\multicolumn{2}{|c|}{Model   Type and Study} & Region & $m, n$ & $\alpha_0$ & $\alpha_1$ & $\alpha_2$ \\ \hline
\multicolumn{1}{|c|}{\multirow{6}{*}{Assumed}} & \multirow{3}{*}{\citet{yang_equilibrium_2011}} & \multirow{3}{*}{-} & \multirow{3}{*}{-} & $0.2$ & $0.75$ & $0.75$ \\ \cline{5-7} 
\multicolumn{1}{|c|}{} &  &  &  & $10$ & $0.5$ & $0.5$ \\ \cline{5-7} 
\multicolumn{1}{|c|}{} &  &  &  & $200$ & $0.25$ & $0.25$ \\ \cline{2-7} 
\multicolumn{1}{|c|}{} & \multirow{3}{*}{\citet{zha_economic_2016}} & \multirow{3}{*}{-} & \multirow{3}{*}{-} & \multirow{3}{*}{$(0,3]$} & $0.6$ & $0.6$ \\ \cline{6-7} 
\multicolumn{1}{|c|}{} &  &  &  &  & $0.7$ & $0.7$ \\ \cline{6-7} 
\multicolumn{1}{|c|}{} &  &  &  &  & $0.8$ & $0.8$ \\ \hline
\multicolumn{1}{|c|}{\multirow{3}{*}{Regression   (empirical data)}} & \citet{yang_returns_2014} & Hong Kong & - & 0.0177 & $0.4080$ & $0.7550$ \\ \cline{2-7} 
\multicolumn{1}{|c|}{} & \multirow{2}{*}{\citet{zhang_efficiency_2019}} & \multirow{2}{*}{Shenzhen} & \multirow{2}{*}{-} & $13.3875$ & $0.4569$ & $0.4569$ \\ \cline{5-7} 
\multicolumn{1}{|c|}{} &  &  &  & $1.2011$ & $0.2052$ & $0.6760$ \\ \hline
\multicolumn{1}{|c|}{\multirow{3}{*}{Regression   (simulated data)}} & \multirow{3}{*}{\citet{wang_-demand_2022}} & Square, $p=1$ & \multirow{3}{*}{-} & $4.2083$ & $0.5274$ & $0.5270$ \\ \cline{3-3} \cline{5-7} 
\multicolumn{1}{|c|}{} &  & Square, $p=2$ &  & $4.1950$ & $0.5246$ & $0.5260$ \\ \cline{3-3} \cline{5-7} 
\multicolumn{1}{|c|}{} &  & New York City &  & $7.2849$ & $0.3877$ & $0.4077$ \\ \hline
\multicolumn{1}{|c|}{Conjecture} & \citet{zha_economic_2016} & - & $m \ll n$ & - & $1$ & - \\ \hline
\multicolumn{1}{|c|}{\multirow{5}{*}{Theoretical}} & \multirow{2}{*}{\citet{daganzo_approximate_1978}} & Quasi-circle, $p=1$ & \multirow{2}{*}{$1=m \ll n$} & $2\sqrt{2/\pi}$ & $1$ & $0.5$ \\ \cline{3-3} \cline{5-7} 
\multicolumn{1}{|c|}{} &  & Quasi-circle, $p=2$ &  & $2$ & $1$ & $0.5$ \\ \cline{2-7} 
\multicolumn{1}{|c|}{} & Equation \eqref{eq:E_Z_nearest} & Hyper-ball & $1=m \ll n$ & $1/R_{\text{B}}(2,p) \Gamma(\frac{3}{2})$ & $1$ & $0.5$ \\ \cline{2-7} 
\multicolumn{1}{|c|}{} & Equation \eqref{eq:E_Z_greedy} & Hyper-ball & $1< m \ll n$ & $1/R_{\text{B}}(2,p) \Gamma(\frac{3}{2})$ & $1$ & $0.5$ \\ \cline{2-7} 
\multicolumn{1}{|c|}{} & Equation \eqref{eq:E_Z_swap} & Hyper-ball & $\forall m, n$  & - & - & - \\ \hline
\end{tabular}%
}
\end{table}

These studies presented very interesting findings (e.g., showing that the parameter values should depend on the region's characteristics), but did not reveal any clear quantitative relationships. 
More importantly, it is not clear whether the Cobb-Douglas form is universally applicable to all settings (e.g., region shape, size, distance metric, dimensions, and $m$ versus $n$ values). Regarding the latter, \citet{wei_calibration_2022} conducted a large set of regression analyses to compare the model predictions with simulation outcomes under 420 different settings. It was found that the parameters in the Cobb-Douglas function vary notably across these settings, and the Cobb-Douglas function might fit quite poorly in certain settings (e.g., when $n$ is small). 

\citet{zha_economic_2016} also argued (without theoretical proof) that one may simply set $\alpha_1 = 1$ when the number of idle vehicles far exceeds the number of customers (i.e., $m \ll n$). 
This might be inspired by the intuition that when $m$ is small, the customer waiting time should be solely dependent on the number of vehicles $n$. The only known theoretical result, derived first in \citet{daganzo_approximate_1978} and later in \citet{arnott_taxi_1996}, states that for the special case of $1=m\ll n$, the customer waiting time can be expressed in the form of Equation \eqref{eq:cobb_t} with $\alpha_0 = 2\sqrt{2/\pi}$ (when $p=1$) or $2$ (when $p=2$), $\alpha_1 = 1$ and $\alpha_2 = 0.5$. 
However, no theoretical proof was provided for the more general case when $m \neq 1$.

Next we show how the convergence property in Proposition \ref{prop:X_greedy_opt} and the distance estimators developed for RBMP-B in Section \ref{sec:RBMP-B} can be used to reveal (i) conditions under which the Cobb-Douglas meeting function Equation \eqref{eq:cobb_t} can be applied, and under such conditions, (ii) the closed-form formulas for the values of parameters $\alpha_0$, $\alpha_1$ and $\alpha_2$.
According to Equation \eqref{eq:E_Z_nearest},
when $1 \le m \ll n$, 
\begin{align}
\mathbb{E}[\accentset{\ast}{C}_{\text{B}}(m,n,2,p)] 
\xrightarrow{n\gg m} \mathbb{E}[\bar{C}_{\text{B}}(1,n,2,p)] 
= R_{\text{B}}(2,p) \cdot \Gamma(\frac{3}{2}) \cdot n^{-\frac{1}{2}} 
= \frac{\Gamma(\frac{3}{2})\left(\Gamma(\frac{2}{p}+1)\right)^{\frac{1}{2}}}{2\Gamma(\frac{1}{p}+1)}\cdot n^{-\frac{1}{2}},
\end{align}
where $R_{\text{B}}(2,p)$ is given by Equation \eqref{eq:R_B}.
Note that this provides us with a Cobb-Douglas form formula. By comparing it with Equation \eqref{eq:cobb_t}, we identify the corresponding parameter values $\alpha_0$, $\alpha_1$, and $\alpha_2$ as follows.
\begin{align*}
    \alpha_0 = \frac{1}{R_{\text{B}}(2,p) \cdot \Gamma(\frac{3}{2})} = \frac{2\Gamma(\frac{1}{p}+1)}{\Gamma(\frac{3}{2})\left(\Gamma(\frac{2}{p}+1)\right)^{\frac{1}{2}}}, \quad
    \alpha_1 = 1, \quad \alpha_2=0.5.
\end{align*} 
For the special case when $m=1$, such a result is clearly consistent with, but generalizes, the result in \citet{daganzo_approximate_1978} with any L$^p$ metric. It also generalizes the result to $m> 1$, as long as $n\gg m$. 
However, it is important to note that for a more accurate estimation, especially in cases with more general combinations of $n$ and $m$, we should refrain from using the simple Cobb-Douglas function. Instead, the formula Equation \eqref{eq:E_Z_swap} shall be used for these cases. These are summarized in the last few rows of Table \ref{tab_cobb}.

\subsection{Optimal Demand Pooling for E-Hailing Services} \label{sec: pooling}

Next, we show how our formulas can be used to optimize e-hailing service strategies.
Consider a 2-dimensional service region (i.e., $D=2$) 
with an area size of 1 [du$^2$], served by 
vehicles traveling at an average speed of 1 [du/tu].
Customer trips are generated from a spatiotemporally homogeneous Poisson process, with rate $\lambda$ [trips/du$^2$-tu], such that each trip's origin and destination are uniformly distributed in the region. 
The vehicles in the system may transition among three states: idle, assigned, and in-service. An idle vehicle, upon assignment to a customer, will start a deadheading trip from its current location to the customer's origin, while at the same time, the customer starts to wait for pickup. After reaching the customer's origin, the vehicle becomes in-service and starts to move towards the customer's destination. When the in-service vehicle drops off the customer, it becomes idle again. 

The service platform periodically matches and assigns idle vehicles to the customers at a set of discrete decision epochs. Between any two consecutive epochs, 
the newly idle vehicles and newly arriving customers are accumulated into their respective pools. At the latter decision epoch, an RBMP-B instance is solved to match all those customers and vehicles in these pools. The average optimal matching distance per customer-vehicle pair, with an expectation of $\mathbb{E}[\accentset{\ast}{C}_{\text{B}}(m,n,2,p)] $, is the expected deadheading distance to pick up a customer. 

While the system might be time-varying, we focus on a specific decision epoch, and denote the length of the corresponding pooling interval as $\tau$. The value of $\tau$ directly affects the numbers of accumulated customers and vehicles and in turn affects $\mathbb{E}[\accentset{\ast}{C}_{\text{B}}(m,n,2,p)] $. Yet, it is also directly experienced as extra ``waiting" or idle time of both the customers and the vehicles.
Therefore, the ride sharing platform needs to identify the optimal $\tau$ value that minimizes the overall system ``waiting" time per customer, including that for vehicle deadheading, $\mathbb{E}[\accentset{\ast}{C}_{\text{B}}(m,n,2,p)] $, and the average wait for pooling, $\frac{\tau}{2}$: 
\begin{align}
\mathbb{E}[\accentset{\ast}{C}_{\text{B}}(m,n,2,p)] + \gamma \cdot \frac{\tau}{2}. \label{opt_obj}
\end{align}
Here $\gamma$ is a relative weight that should normally be larger than 1; e.g., to account for the customer anxiety before a vehicle is assigned, and the fact that we assume there are more idle vehicles than customers at each matching epoch. 

Next, two types of systems with different vehicle arrival dynamics are analyzed: one is an ``open-loop" system without vehicle conservation, possibly representing services with freelance drivers; the other is a ``closed-loop" system with conservation of a fixed fleet of vehicles, possibly representing services with full-time drivers or robo-taxis.

\subsubsection{Open-loop system} \label{sec: open_sys}
In the open-loop system that we consider, new idle vehicles arrive from outside of the system according to an independent process, and their initial locations are uniformly distributed in the region. This mimics a service scenario discussed in \citet{yang_optimizing_2020}, where the total fleet size in the system, as well as the numbers of vehicles in the three states, changes over time. 
Let $n_i$ denote the number of idle vehicles in the system right after a previous decision epoch, and let $\lambda'$ [trips/du$^2$-tu] denote the expected arrival rate of idle vehicles (including those transitioned from the in-service state, and those new arrivals from outside). 
Right before the next decision epoch, the expected number of waiting customers in the pool is $\lambda \tau$, and the expected number of idle vehicles is $n_i + \lambda' \tau$. Hence, from formula \eqref{eq:E_Z_swap}, we know that $\mathbb{E}[\accentset{\ast}{C}_{\text{B}}(m,n,2,p)]$ is solely dependent on $\tau$, as
\begin{align}
m = \lambda \tau, \quad n = n_i + \lambda' \tau \label{opt_m_n}.
\end{align}
Then,  we can write a simple optimization model to minimize \eqref{opt_obj}, as follows. 
    \begin{align}
        \min_{\tau} \quad & \eqref{opt_obj}  \nonumber \\
        \text{s.t.} \quad
        & \eqref{eq:E_Z_swap} \text{ and } \eqref{opt_m_n}, \nonumber\\
        \quad & \tau \ge \frac{1}{\lambda} \label{con_tau}. 
    \end{align}
Here Equation \eqref{eq:E_Z_swap} can also be replaced by Equation \eqref{eq:E_Z_greedy_kappa} 
for computational efficiency, as revealed in Section \ref{sec:experiment_B}.
Constraint \eqref{con_tau} specifies a boundary condition that at least one customer arrives during the pooling interval. 

This optimization problem takes $\lambda$, $\lambda'$ and $n_i$ as input and involves only one decision variable $\tau$. Thus, it can be easily solved numerically (without simulations). 
Furthermore, we can obtain some qualitative insights from the analytical formulation. Note the following: (i) from Equation \eqref{opt_m_n}, $m$ and $n$ both monotonically increase with $\tau$; (ii) intuitively, $\mathbb{E}[\accentset{\ast}{C}_{\text{B}}(m,n,2,p)]$ decreases monotonically with $n$ and increases monotonically with $m$. Hence, it is possible for $\mathbb{E}[\accentset{\ast}{C}_{\text{B}}(m,n,2,p)]$, as well as Equation \eqref{opt_obj}, to either increase or decrease with $\tau$. 
Two representative scenarios may arise here: (i) the objective function \eqref{opt_obj} monotonically increases with $\tau$ and the optimal solution is simply $\tau^*=\frac{1}{\lambda}$
, which implies that 
instantly matching each customer upon its arrival (i.e., ``instant matching'' as in \citet{daganzo_general_2019, ouyang_measurement_2023, shen_dynamic_2023}) might be a (near-) optimal strategy; 
or (ii) the objective function \eqref{opt_obj} 
is a non-monotonic function of $\tau$, such that the optimal solution $\tau^* > \frac{1}{\lambda}$, which indicates that pooling multiple customers into a batch for matching (i.e., ``batch matching" as in \citet{yang_optimizing_2020}) would be more favorable.

\subsubsection{Closed-loop system} \label{sec: close_sys}
Now, we consider a closed-loop system which operates with a fixed fleet of $S$ identical vehicles (i.e. no idle vehicles arrive from outside). 
Following the steady-state aspatial queuing network model in \citet{daganzo_public_2010}, we know that 
(i) the idle vehicle arrival rate equals the customer arrival rate; i.e., 
\begin{align} \label{opt_lambda}
\lambda' = \lambda, 
\end{align}
and (ii) the numbers of vehicles in the three states will reach a certain equilibrium (or jump between two equilibria). 
We denote the numbers of assigned and in-service vehicles right after the previous decision epoch as $n_a$ and $n_s$, respectively. They must satisfy vehicle conservation; i.e.,
\begin{align}
    S =n_i+n_s+n_a \label{veh_con}.
\end{align} 
Let us further denote $l$ [du] as the expected travel distance for an in-service vehicle to deliver a customer inside the unit-volume hyper-ball. Its value depends solely on the distance metric. Per Little's formula, in the steady state, $n_s$ and $n_a$ must also satisfy the following:
\begin{align} 
    n_s =  \lambda l, \quad  n_a = \lambda \left( \mathbb{E}[\accentset{\ast}{C}_{\text{B}}(m,n,2,p)] + \frac{\tau}{2} \right) \label{little}.
\end{align}
Then, the optimization model 
in Section \ref{sec: open_sys} can be adapted 
into the following:
    \begin{align}
        \min_{\tau} \quad 
            & \eqref{opt_obj}, \nonumber\\
            \text{s.t.} \quad
            & \eqref{eq:E_Z_swap} \text{ and } \eqref{opt_m_n}-\eqref{little}, \nonumber\\
            & 
            0 \le n_i \le S - \lambda \tau, \quad \lambda \tau \le n_a \le S, \quad 0 \le n_s \le S. \label{constraints-opt}
    \end{align}
Constraints \eqref{constraints-opt} are boundary conditions that ensure the number of idle, assigned, or in-service vehicles, either before or after each decision epoch, cannot be negative nor exceed the fleet size. 

This new optimization problem takes $\lambda$ and $S$ as input, 
but $n_i$ now must be solved out of Equations \eqref{opt_m_n}, \eqref{opt_lambda}-\eqref{little} for any $\tau$. 
Two roots of $n_i$ may exist in correspondence to two equilibrium states of the system 
\citep{ouyang_measurement_2023, shen_dynamic_2023}: the one with a larger $n_i$ value is referred to as the efficient equilibrium; and the other with a smaller $n_i$ value is referred to as the inefficient equilibrium (i.e., the so-called ``wild goose chase (WGC)" phenomenon). 
The optimal value of $\tau$ shall depend on the relative values of $m$ and $n$, similar to the discussed cases in Section \ref{sec: open_sys}, around either equilibrium vertex. 



\subsection{Verification of the Effectiveness of Demand Pooling}

Two simulation programs are developed to verify the effectiveness of the demand pooling strategy for the open- and closed-loop e-hailing taxi systems in Section \ref{sec: pooling}. 
A series of numerical experiments are conducted to compare simulation measurements with model results. In all experiments, we set $p = 2$ and $\gamma = 1$. For computational efficiency, $\mathbb{E}[\accentset{\ast}{C}_{\text{B}}(m,n,2,2)]$ is estimated from the approximate formula \eqref{eq:E_Z_greedy_kappa} with $\kappa=0$. The effectiveness of such approximation has been discussed in Section \ref{sec:experiment_B}. 

The program for the open-loop system simulates the arrivals of customers and vehicles in a specific pooling interval before a decision epoch. 
At the beginning of each simulation (i.e., time 0), $n_i$ idle vehicles are distributed uniformly within a unit-area circle. 
Then, customer trips and new idle vehicles are generated from two independent homogeneous spatio-temporal Poisson processes with rates $\lambda$ and $\lambda'$ [trips/du$^2$-tu], respectively. We select a sample of $\tau$ [tu] values. For each $\tau$ value, we use the modified Jonker-Volgenant algorithm \citep{modified_JV} to solve an RBMP-B instance with the accumulated customers and idle vehicles by time $\tau$, and each matched customer's waiting times for pooling and for pickup are recorded. 

\begin{figure}[ht!] 
\centering 
\includegraphics[scale=0.535]{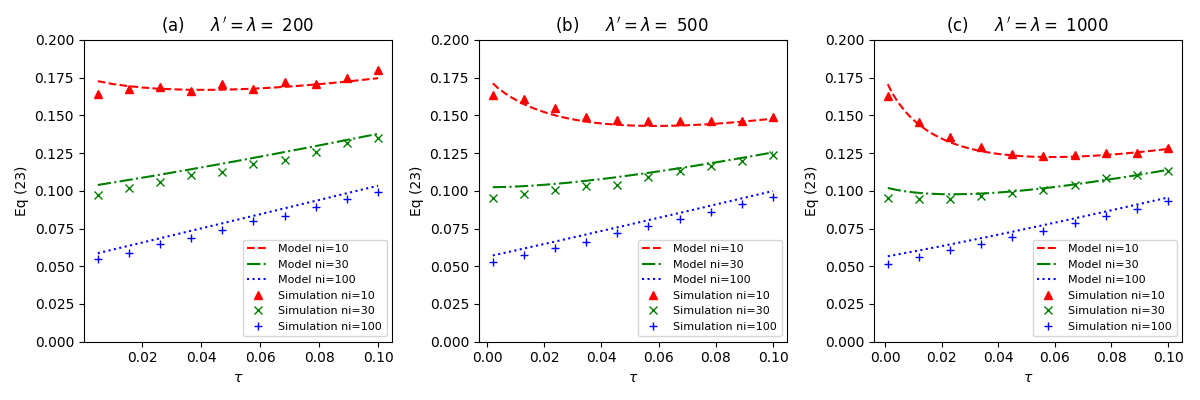}
\caption{Estimated vs. simulated Eq \eqref{opt_obj} of an open-loop system under different $\lambda$ and $n_i$.}
\label{fig: open_loop}
\end{figure}

Experiments are conducted under a set of parameter combinations: 
$\lambda=\lambda' \in \{200, 500, 1000\}$, and $n_i \in \{10, 30, 100\}$. 
We select 10 discrete $\tau$ values that are evenly distributed between $\frac{1}{\lambda}$ and $0.1$ [tu]. For each parameter combination and $\tau$ value, we run $30/\tau$ simulations and take the average. 
The results are shown in Figure \ref{fig: open_loop}. 
Model estimations of Equation \eqref{opt_obj} for $n_i \in \{10, 30, 100\}$ are represented by the red dashed curve, green dash-dot curve, and blue dotted curve, respectively. Meanwhile, the simulation measurements are represented by the red triangles, green cross markers, and blue plus markers, respectively. 
It can be observed that the objective value estimations, for all $\tau$, match quite well with their corresponding simulation averages. 
When demand is relatively lower, $\lambda=\lambda' \in \{200, 500\}$ and there are abundant idle vehicles, $n_i = 30$ and $100$, as we can see from Figures \ref{fig: open_loop}(a)-(b), the value of Equation \eqref{opt_obj} monotonically increases with $\tau$, and hence $\tau^* = \frac{1}{\lambda}$. This corresponds to the case of scenario (i) as discussed in Section \ref{sec: open_sys}, where ``instant matching" is more favorable. However, when the number of idle vehicles is smaller, $n_i = 10$, Equation \eqref{opt_obj} first decreases and then increases with $\tau$, such that notably $\tau^* > \frac{1}{\lambda}$. This corresponds to the case of scenario (ii) where ``batch matching" is more favorable. 
Also, when demand increases to $\lambda=\lambda' = 1000$, as shown in Figure \ref{fig: open_loop}(c), it is more likely for batch matching to be favorable for $n_i = 10$ and $30$. 
These observations are consistent with our analytical insights in Section \ref{sec: pooling}, and suggest that: the simpler instant matching strategy might be very suitable when the taxi system has a sufficient number of idle vehicles and it does not expect a large number of new idle vehicles arriving from outside of the system; meanwhile, batch matching is likely to be beneficial when the system is about to run out of idle vehicles, or if it expects a large number of new idle vehicles to arrive from outside. 

For the closed-loop system, an agent-based simulation program is developed. It tracks each customer's and vehicle's entire travel experience, including arrival, assignment, pickup and drop-off. 
The program generates customer trips from a Poisson process with rate $\lambda$ [trips/du$^2$-tu], similar to that for the open-loop system, but now a fixed fleet of $S$ vehicles are used to serve these trips. 
The customers and vehicles evolve following the processes described in Section \ref{sec: pooling}, and bipartite matching between waiting customers and idle vehicles is conducted at times $\tau$, $2\tau, \ldots$ [tu]. Any unmatched customers (e.g., due to random demand surges) are removed from the system and considered as lost. 

Experiments are conducted under a set of parameter combinations. 
We use the same $\lambda$ and $\tau$ values as those for the open-loop system. 
We can compute the minimum required fleet size to serve all customers in a steady state \citep{daganzo_public_2010}, which are approximately 143, 330, and 630 for $\lambda \in \{200, 500, 1000\}$, respectively. Based on this information, we try two fleet sizes to represent smaller and large fleet scenarios: the values of $S \in \{155, 220\}$, $\{360, 500\}$ and $\{680, 1000\}$, respectively. 
For each parameter combination, one simulation run with a duration of 30 [tu] is performed. The first 6 [tu] of each simulation run is considered as the warm-up period, and we only record data from the later part. The corresponding model estimations are computed based on Equations \eqref{opt_obj}-\eqref{opt_m_n} and \eqref{opt_lambda}-\eqref{little}. 

\begin{figure}[ht!] 
\centering 
\includegraphics[scale=0.535]{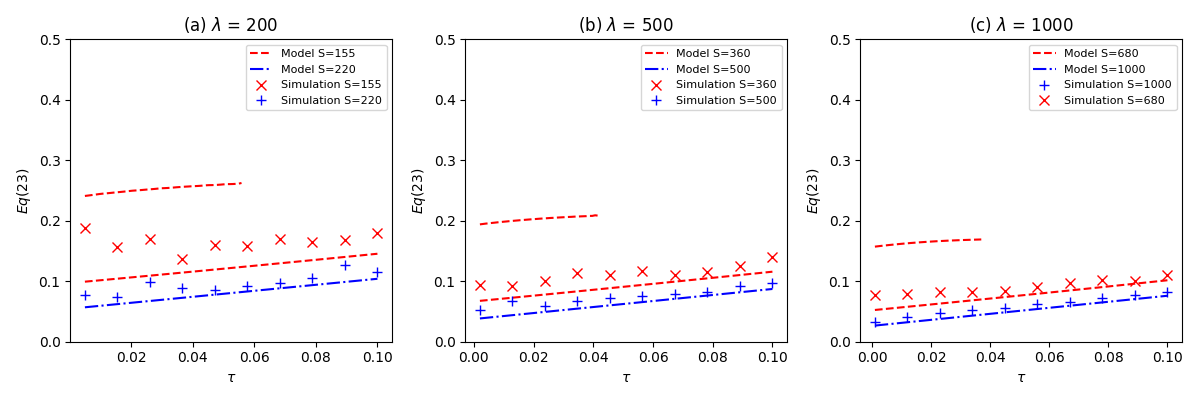}
\caption{Estimated vs. simulated Eq \eqref{opt_obj} of a closed-loop system under different $\lambda$ and $S$.
}
\label{fig: close_loop}
\end{figure}

The results are shown in Figure \ref{fig: close_loop}. 
Model estimations are represented by the red dashed curve and blue dash-dot curve, respectively, for the small and large $S$ values. The corresponding simulation averages are represented by the red cross markers and blue plus markers, respectively. 
When $S$ is sufficiently large, only one blue dash-dot curve appears in each sub-figure, indicating the system has only one (efficient) equilibrium state \citep{ouyang_measurement_2023}. In this case, Equation \eqref{opt_obj} monotonically increases with $\tau$, which falls into the case of scenario (i), 
and the simulation averages across all demand levels match well with the model estimations (as shown in Figure \ref{fig: close_loop}(a)-(c)). 
When $S$ is relatively small, two red dashed curves appear in each sub-figure, indicating the system can potentially have both inefficient and efficient equilibrium states. It can also be seen that the upper red dashed curve is truncated as $\tau$ increases beyond a certain threshold, which indicates that the inefficient equilibrium (WGC) may completely disappear if demand is pooled for a sufficiently long time period --- a significant potential benefit of demand pooling. 
Yet, if we only focus on only the efficient equilibrium, we shall still do instant matching because Equation \eqref{opt_obj} for the efficient equilibrium monotonically increases with $\tau$. 

In addition, we shall note that when two equilibrium states exist, the simulation measurements may be the average of the two equilibrium vertices (weighted by the duration the system spends in either equilibrium); see more discussion in \citep{ouyang_measurement_2023}. 
Yet, most of these simulation measurements still align well with the model estimations (especially when the fleet size is large) for the efficient equilibrium --- implying that the system predominantly stay in the efficient equilibrium state under those conditions. 
All these observations imply that, in a closed-loop system, instant matching might again be an effective strategy overall. Batch matching may only be beneficial for counteracting the WGC phenomenon 
when the system with a small fleet size is currently stuck in an inefficient equilibrium state. 


\section{Conclusion} \label{sec: conclusion}

This paper proposes a comprehensive modeling framework that yields
analytical formulas to estimate the expected optimal matching cost (or distance) from a random bipartite matching problem (RBMP). We study three interrelated but increasingly complex versions of RBMPs: 
(i) RBMP-I with i.i.d. edge costs; 
(ii) RBMP-S with uniformly distributed vertices on the surface of a hyper-sphere in a $D$-dimensional Euclidean space; and 
(iii) RBMP-B with uniformly distributed vertices inside a hyper-ball in a $D$-dimensional L$^p$ space. 
The model first derives closed-form formulas for the conditional expected distance for a demand vertex matched to its $k$-th nearest neighbor in RBMP-I, followed by that for the probability of such a match (via a two-step heuristic process). Using these RBMP-I results as a first approximation, we then extend the model to account for spatial correlation in RBMP-S and boundary effects in RBMP-B. 
Through a series of Monte-Carlo simulation experiments, the proposed formulas are shown to be able to provide accurate estimations under a variety of conditions (e.g., regarding the numbers of bipartite vertices, spatial dimensions, and distance metrics). 

Using mobility services as an example, this paper shows how the proposed formulas can be used to provide managerial insights on the expected matches between randomly distributed customers and service vehicles. For instance, 
the formulas provide a theoretical explanation on when and why the empirically assumed Cobb-Douglas matching functions are suitable for mobility systems. 
The proposed distance formulas are also particularly useful for strategic planning or allocation of resources in mobility services. For example, this paper shows how they can be integrated into optimization models to plan e-hailing taxi operational strategies (e.g., choosing demand pooling intervals) under two system settings. The predicted performance of these strategies are verified by agent-based simulations. 
The results not only show that the model estimates match quite well with the simulation measurements under all tested service settings, but also provide insights on suitability of instant matching vs. batch matching strategies: instant matching is more favorable when the system has a sufficiently large fleet, while batch matching can be beneficial when the system with a small fleet (e.g., to avoid the inefficient equilibrium).

The current modeling approach builds upon several simplifying assumptions or hypotheses that could be relaxed in the future. 
First, in deriving the matching probabilities in Section \ref{sec:P_k}, we employ a two-step heuristic method and use the refined matching probability $\hat{\mathbb{P}}(k)$ as an approximation of the optimal matching probability $\accentset{\ast}{\mathbb{P}}(k)$. Although the resulting formula already yields highly accurate estimations in general, it will be ideal to derive 
exact formulas based on strict optimality conditions. Moreover, the current formula for $\hat{\mathbb{P}}(k)$, while already involving certain approximations (such as for $\Pr\{\Delta_{q'',q'}\le\Delta_{k'',k'}\}$), is still computationally expensive. A useful direction would be to simplify the formulas to reduce computational cost while maintaining a similar level of accuracy. 
Third, in deriving the correction term that accounts for spatial correlation among vertex distributions in $D$-dimensional spaces, the proposed functional form of Equation \eqref{eq:delta} is motivated by insights from the previous literature, and the value of the correction coefficient $d$, which captures the impact of imbalance between the two vertex sets, is obtained through data regression. A promising direction for future work would be to establish more theoretically grounded correlation terms, at least for special cases such as $D=1$ and $D=2$. For higher dimensions ($D \geq 3$), the correlations become negligible and are thus less critical. 

The proposed model and formulas could also be further extended to other problem settings or application contexts. 
First, more accurate asymptotic approximations for certain special cases (such as when $m \approx n \gg 1$) hold the promise for providing more theoretical and practical connections to those in the field of statistical physics (e.g., those in \citet{caracciolo_scaling_2014}) 
for the balanced case (i.e., $m=n$). It would be interesting to see if insights and theoretical justifications (similar to those in Section \ref{sec: cobb}) could be found for those results.
Second, it will be interesting to extend the model as building blocks for other problems. 
For example, the current RBMP is defined as a static problem, and the matching decisions are made with full knowledge of all vertex locations. 
It can be generalized into a dynamic matching problem where matching decisions must be made based on incomplete information on future vertex arrivals. 
As discussed in Section \ref{sec: literature}, \citet{kanoria_dynamic_2022} showed, by adopting a hierarchical greedy algorithm, how the bounds to the expected matching distance for a static problem can be adapted into bounds of minimum achievable matching distance in a dynamic problem. 
Now that we have formulas to estimate the exact matching distance rather than the bounds, it would be interesting to generalize our model to directly predict the expected matching distance in such a dynamic setting. 



\newpage
\bibliographystyle{apalike}
\bibliography{Batch}

\begin{thebibliography}{}

\bibitem[Ahlgren, 2014]{ahlgren_probability_2014}
Ahlgren, J. (2014).
\newblock The {Probability} {Distribution} for {Draws} {Until} {First} {Success} {Without} {Replacement}.
\newblock arXiv:1404.1161 [math] version: 1.

\bibitem[Ajtai et~al., 1984]{ajtai_optimal_1984}
Ajtai, M., Komlós, J., and Tusnády, G. (1984).
\newblock On optimal matchings.
\newblock {\em Combinatorica}, 4(4):259--264.

\bibitem[Alm and Sorkin, 2002]{alm_exact_2002}
Alm, S.~E. and Sorkin, G.~B. (2002).
\newblock Exact {Expectations} and {Distributions} for the {Random} {Assignment} {Problem}.
\newblock {\em Combinatorics, Probability and Computing}, 11(3):217--248.

\bibitem[Arnold et~al., 2008]{arnold_first_2008}
Arnold, B.~C., Balakrishnan, N., and Nagaraja, H.~N. (2008).
\newblock {\em A {First} {Course} in {Order} {Statistics} ({Classics} in {Applied} {Mathematics})}.
\newblock Society for Industrial and Applied Mathematics, USA.

\bibitem[Arnott, 1996]{arnott_taxi_1996}
Arnott, R. (1996).
\newblock Taxi {Travel} {Should} {Be} {Subsidized}.
\newblock {\em Journal of Urban Economics}, 40(3):316--333.

\bibitem[Boniolo et~al., 2014]{boniolo_correlation_2014}
Boniolo, E., Caracciolo, S., and Sportiello, A. (2014).
\newblock Correlation function for the {Grid}-{Poisson} {Euclidean} matching on a line and on a circle.
\newblock {\em Journal of Statistical Mechanics: Theory and Experiment}, 2014(11):P11023.

\bibitem[Brunetti et~al., 1991]{brunetti_extensive_1991}
Brunetti, R., Krauth, W., Mézard, M., and Parisi, G. (1991).
\newblock Extensive {Numerical} {Simulations} of {Weighted} {Matchings}: {Total} {Length} and {Distribution} of {Links} in the {Optimal} {Solution}.
\newblock {\em Europhysics Letters (EPL)}, 14(4):295--301.

\bibitem[Cai et~al., 2013]{cai_distributions_2013}
Cai, T., Fan, J., and Jiang, T. (2013).
\newblock Distributions of {Angles} in {Random} {Packing} on {Spheres}.
\newblock arXiv:1306.0256 [math].

\bibitem[Caracciolo et~al., 2019]{caracciolo_selberg_2019}
Caracciolo, S., Di~Gioacchino, A., Malatesta, E.~M., and Molinari, L.~G. (2019).
\newblock Selberg integrals in {1D} random {Euclidean} optimization problems.
\newblock {\em Journal of Statistical Mechanics: Theory and Experiment}, 2019(6):063401.
\newblock Publisher: IOP Publishing and SISSA.

\bibitem[Caracciolo et~al., 2014]{caracciolo_scaling_2014}
Caracciolo, S., Lucibello, C., Parisi, G., and Sicuro, G. (2014).
\newblock Scaling hypothesis for the {Euclidean} bipartite matching problem.
\newblock {\em Physical Review E}, 90(1):012118.

\bibitem[Caracciolo and Sicuro, 2014]{caracciolo_one-dimensional_2014}
Caracciolo, S. and Sicuro, G. (2014).
\newblock One-dimensional {Euclidean} matching problem: {Exact} solutions, correlation functions, and universality.
\newblock {\em Physical Review E}, 90(4):042112.

\bibitem[Caracciolo and Sicuro, 2015a]{caracciolo_quadratic_2015}
Caracciolo, S. and Sicuro, G. (2015a).
\newblock Quadratic stochastic {Euclidean} bipartite matching problem.
\newblock {\em Physical Review Letters}, 115(23):230601.

\bibitem[Caracciolo and Sicuro, 2015b]{caracciolo_scaling_2015}
Caracciolo, S. and Sicuro, G. (2015b).
\newblock Scaling hypothesis for the {Euclidean} bipartite matching problem. {II}. {Correlation} functions.
\newblock {\em Physical Review E}, 91(6):062125.

\bibitem[Coppersmith and Sorkin, 1998]{coppersmith_constructive_1998}
Coppersmith, D. and Sorkin, G.~B. (1998).
\newblock Constructive {Bounds} and {Exact} {Expectations} for the {Random} {Assignment} {Problem}.
\newblock In Luby, M., Rolim, J. D.~P., and Serna, M., editors, {\em Randomization and {Approximation} {Techniques} in {Computer} {Science}}, pages 319--330, Berlin, Heidelberg. Springer.

\bibitem[Crouse, 2016]{modified_JV}
Crouse, D.~F. (2016).
\newblock On implementing 2d rectangular assignment algorithms.
\newblock {\em IEEE Transactions on Aerospace and Electronic Systems}, 52(4):1679--1696.

\bibitem[Daganzo, 1978]{daganzo_approximate_1978}
Daganzo, C.~F. (1978).
\newblock An approximate analytic model of many-to-many demand responsive transportation systems.
\newblock {\em Transportation Research}, 12(5):325--333.

\bibitem[Daganzo, 2010]{daganzo_public_2010}
Daganzo, C.~F. (2010).
\newblock Public {Transportation} {Systems}:{Basic} {Principles} of {System} {Design},{Operations} {Planning} and {Real}-{TimeControl}.

\bibitem[Daganzo and Ouyang, 2019]{daganzo_general_2019}
Daganzo, C.~F. and Ouyang, Y. (2019).
\newblock A general model of demand-responsive transportation services: {From} taxi to ridesharing to dial-a-ride.
\newblock {\em Transportation Research Part B: Methodological}, 126:213--224.

\bibitem[Daganzo et~al., 2020]{daganzo_analysis_2020}
Daganzo, C.~F., Ouyang, Y., and Yang, H. (2020).
\newblock Analysis of ride-sharing with service time and detour guarantees.
\newblock {\em Transportation Research Part B: Methodological}, 140:130--150.

\bibitem[Daganzo and Smilowitz, 2004]{daganzo_bounds_2004}
Daganzo, C.~F. and Smilowitz, K.~R. (2004).
\newblock Bounds and {Approximations} for the {Transportation} {Problem} of {Linear} {Programming} and {Other} {Scalable} {Network} {Problems}.
\newblock {\em Transportation Science}, 38(3):343--356.

\bibitem[Houdayer et~al., 1998]{houdayer_comparing_1998}
Houdayer, J., Boutet De~Monvel, J., and Martin, O. (1998).
\newblock Comparing mean field and {Euclidean} matching problems.
\newblock {\em The European Physical Journal B}, 6(3):383--393.

\bibitem[Jonker and Volgenant, 1987]{jonker_shortest_1987}
Jonker, R. and Volgenant, A. (1987).
\newblock A shortest augmenting path algorithm for dense and sparse linear assignment problems.
\newblock {\em Computing}, 38(4):325--340.

\bibitem[Kanoria, 2022]{kanoria_dynamic_2022}
Kanoria, Y. (2022).
\newblock Dynamic {Spatial} {Matching}.
\newblock In {\em Proceedings of the 23rd {ACM} {Conference} on {Economics} and {Computation}}, pages 63--64, Boulder CO USA. ACM.

\bibitem[Kershaw, 1983]{kershaw_extensions_1983}
Kershaw, D. (1983).
\newblock Some {Extensions} of {W}. {Gautschi}'s {Inequalities} for the {Gamma} {Function}.
\newblock {\em Mathematics of Computation}, 41(164):607--611.
\newblock Publisher: American Mathematical Society.

\bibitem[Kuhn, 1955]{kuhn_hungarian_1955}
Kuhn, H.~W. (1955).
\newblock The {Hungarian} method for the assignment problem.
\newblock {\em Naval Research Logistics Quarterly}, 2(1-2):83--97.

\bibitem[Lei and Ouyang, 2024]{lei_average_2024}
Lei, C. and Ouyang, Y. (2024).
\newblock Average minimum distance to visit a subset of random points in a compact region.
\newblock {\em Transportation Research Part B: Methodological}, 181:102904.

\bibitem[Li, 2010]{li_concise_2010}
Li, S. (2010).
\newblock Concise {Formulas} for the {Area} and {Volume} of a {Hyperspherical} {Cap}.
\newblock {\em Asian Journal of Mathematics \& Statistics}, 4(1):66--70.

\bibitem[Linusson and Wästlund, 2004]{linusson_proof_2004}
Linusson, S. and Wästlund, J. (2004).
\newblock A proof of {Parisi}’s conjecture on the random assignment problem.
\newblock {\em Probability Theory and Related Fields}, 128(3):419--440.

\bibitem[Liu and Ouyang, 2021]{liu_mobility_2021}
Liu, Y. and Ouyang, Y. (2021).
\newblock Mobility service design via joint optimization of transit networks and demand-responsive services.
\newblock {\em Transportation Research Part B: Methodological}, 151:22--41.

\bibitem[Liu and Ouyang, 2023]{liu_planning_2023}
Liu, Y. and Ouyang, Y. (2023).
\newblock Planning ride-pooling services with detour restrictions for spatially heterogeneous demand: {A} multi-zone queuing network approach.
\newblock {\em Transportation Research Part B: Methodological}, 174:102779.

\bibitem[Mézard and Parisi, 1985]{mezard_replicas_1985}
Mézard, M. and Parisi, G. (1985).
\newblock Replicas and optimization.
\newblock {\em Journal de Physique Lettres}, 46(17):771--778.

\bibitem[Mézard and Parisi, 1988]{mezard_euclidean_1988}
Mézard, M. and Parisi, G. (1988).
\newblock The {Euclidean} matching problem.
\newblock {\em Journal de Physique}, 49(12):2019--2025.

\bibitem[Nair et~al., 2005]{nair_proofs_2005}
Nair, C., Prabhakar, B., and Sharma, M. (2005).
\newblock Proofs of the {Parisi} and {Coppersmith}‐{Sorkin} random assignment conjectures.
\newblock {\em Random Structures \& Algorithms}, 27(4):413--444.

\bibitem[Olver et~al., 2010]{olver_nist_2010}
Olver, F.~W., Lozier, D.~W., Boisvert, R.~F., and Clark, C.~W. (2010).
\newblock {\em {NIST} {Handbook} of {Mathematical} {Functions}}.
\newblock Cambridge University Press, USA, 1st edition.

\bibitem[Ouyang and Yang, 2023]{ouyang_measurement_2023}
Ouyang, Y. and Yang, H. (2023).
\newblock Measurement and mitigation of the “wild goose chase” phenomenon in taxi services.
\newblock {\em Transportation Research Part B: Methodological}, 167:217--234.

\bibitem[Ouyang et~al., 2021]{ouyang_performance_2021}
Ouyang, Y., Yang, H., and Daganzo, C.~F. (2021).
\newblock Performance of reservation-based carpooling services under detour and waiting time restrictions.
\newblock {\em Transportation Research Part B: Methodological}, 150:370--385.

\bibitem[Parisi, 1998]{parisi_conjecture_1998}
Parisi, G. (1998).
\newblock A {Conjecture} on random bipartite matching.
\newblock arXiv:cond-mat/9801176.

\bibitem[Shen and Ouyang, 2023]{shen_dynamic_2023}
Shen, S. and Ouyang, Y. (2023).
\newblock Dynamic and {Pareto}-improving swapping of vehicles to enhance integrated and modular mobility services.
\newblock {\em Transportation Research Part C: Emerging Technologies}, 157:104366.

\bibitem[Simmons et~al., 2019]{simmons_motifs_2019}
Simmons, B.~I., Cirtwill, A.~R., Baker, N.~J., Wauchope, H.~S., Dicks, L.~V., Stouffer, D.~B., and Sutherland, W.~J. (2019).
\newblock Motifs in bipartite ecological networks: uncovering indirect interactions.
\newblock {\em Oikos}, 128(2):154--170.

\bibitem[Talagrand, 1992]{talagrand_matching_1992}
Talagrand, M. (1992).
\newblock Matching {Random} {Samples} in {Many} {Dimensions}.
\newblock {\em The Annals of Applied Probability}, 2(4):846--856.
\newblock Publisher: Institute of Mathematical Statistics.

\bibitem[Tanay et~al., 2004]{tanay_revealing_2004}
Tanay, A., Sharan, R., Kupiec, M., and Shamir, R. (2004).
\newblock Revealing modularity and organization in the yeast molecular network by integrated analysis of highly heterogeneous genomewide data.
\newblock {\em Proceedings of the National Academy of Sciences}, 101(9):2981--2986.

\bibitem[Wang et~al., 2022]{wang_-demand_2022}
Wang, G., Zhang, H., and Zhang, J. (2022).
\newblock On-{Demand} {Ride}-{Matching} in a {Spatial} {Model} with {Abandonment} and {Cancellation}.
\newblock {\em Operations Research}.

\bibitem[Wei et~al., 2022]{wei_calibration_2022}
Wei, S., Feng, S., Ke, J., and Yang, H. (2022).
\newblock Calibration and validation of matching functions for ride-sourcing markets.
\newblock {\em Communications in Transportation Research}, 2:100058.

\bibitem[Wise, 1960]{wise_normalizing_1960}
Wise, M.~E. (1960).
\newblock On normalizing the incomplete beta-function for fitting to dose-response curves.
\newblock {\em Biometrika}, 47(1-2):173--175.

\bibitem[Wu et~al., 2022]{wu_graph_2022}
Wu, S., Sun, F., Zhang, W., Xie, X., and Cui, B. (2022).
\newblock Graph {Neural} {Networks} in {Recommender} {Systems}: {A} {Survey}.
\newblock {\em ACM Computing Surveys}, 55(5):97:1--97:37.

\bibitem[Yang et~al., 2010]{yang_equilibria_2010}
Yang, H., Leung, C.~W., Wong, S., and Bell, M.~G. (2010).
\newblock Equilibria of bilateral taxi–customer searching and meeting on networks.
\newblock {\em Transportation Research Part B: Methodological}, 44(8-9):1067--1083.

\bibitem[Yang et~al., 2020]{yang_optimizing_2020}
Yang, H., Qin, X., Ke, J., and Ye, J. (2020).
\newblock Optimizing matching time interval and matching radius in on-demand ride-sourcing markets.
\newblock {\em Transportation Research Part B: Methodological}, 131:84--105.

\bibitem[Yang and Yang, 2011]{yang_equilibrium_2011}
Yang, H. and Yang, T. (2011).
\newblock Equilibrium properties of taxi markets with search frictions.
\newblock {\em Transportation Research Part B: Methodological}, 45(4):696--713.

\bibitem[Yang et~al., 2014]{yang_returns_2014}
Yang, T., Yang, H., Wong, S.~C., and Sze, N.~N. (2014).
\newblock Returns to scale in the production of taxi services: an empirical analysis.
\newblock {\em Transportmetrica A: Transport Science}, 10(9):775--790.

\bibitem[Yu et~al., 2020]{yu_learning_2020}
Yu, H., Ye, W., Feng, Y., Bao, H., and Zhang, G. (2020).
\newblock Learning {Bipartite} {Graph} {Matching} for {Robust} {Visual} {Localization}.
\newblock In {\em 2020 {IEEE} {International} {Symposium} on {Mixed} and {Augmented} {Reality} ({ISMAR})}, pages 146--155, Porto de Galinhas, Brazil. IEEE.

\bibitem[Zha et~al., 2016]{zha_economic_2016}
Zha, L., Yin, Y., and Yang, H. (2016).
\newblock Economic analysis of ride-sourcing markets.
\newblock {\em Transportation Research Part C: Emerging Technologies}, 71:249--266.

\bibitem[Zhai et~al., 2024]{zhai_average_2024}
Zhai, Y., Shen, S., and Ouyang, Y. (2024).
\newblock Average {Distance} of {Random} {Bipartite} {Matching} in {Discrete} {Networks}.
\newblock arXiv:2409.18292 [math].

\bibitem[Zhang et~al., 2023]{zhang_survey_2023}
Zhang, J., Liu, C., Li, X., Zhen, H.-L., Yuan, M., Li, Y., and Yan, J. (2023).
\newblock A survey for solving mixed integer programming via machine learning.
\newblock {\em Neurocomputing}, 519:205--217.

\bibitem[Zhang et~al., 2019]{zhang_efficiency_2019}
Zhang, K., Chen, H., Yao, S., Xu, L., Ge, J., Liu, X., and Nie, M. (2019).
\newblock An {Efficiency} {Paradox} of {Uberization}.
\newblock {\em SSRN Electronic Journal}.

\bibitem[Zhou et~al., 2007]{zhou_bipartite_2007}
Zhou, T., Ren, J., Medo, M., and Zhang, Y. (2007).
\newblock Bipartite network projection and personal recommendation.
\newblock {\em Physical Review E}, 76(4):046115.

\end{thebibliography}

\newpage
\appendix

\section{Error Bounds of Using the Gamma Approximations}
\label{app:kappa_approx}

Based on Equation \eqref{eq:gamma_ratio}, the gamma function ratio outside the summation in Equation \eqref{eq:E_X_greedy_geo}, $\frac{\Gamma(n+1)}{\Gamma(n+\frac{1}{D}+1)}$, can be well approximated by $(n+1)^{-\frac{1}{D}} \approx n^{-\frac{1}{D}}$ as $n \gg 1$.
Similarly, the gamma function ratios inside the summation can be well approximated as $ \frac{\Gamma(k+\frac{1}{D})}{\Gamma(k)}\approx k^{-\frac{1}{D}}$ and $\frac{\Gamma(i+\frac{1}{D})}{\Gamma(i)}\approx i^{-\frac{1}{D}}$ for $i \ge k \gg 1$.
The approximation error arises only when either $k$ or $i$ is very small. The following lemma presents the error bounds of using Equation \eqref{eq:gamma_ratio} to approximate the gamma functions.
\begin{lemma} \label{lemma:error_gamma}
For $z \in \mathbb{Z}^+, D\in \mathbb{Z}^+$, 
\begin{align} \label{eq:error_gamma} 
    0 \le z^{\frac{1}{D}} - \frac{\Gamma(z+\frac{1}{D})}{\Gamma(z)} < 
    1-\frac{\sqrt{3}}{2}.
\end{align}
\end{lemma}

\noindent{\bf Proof}.
For the trivial case when $D = 1$, $\frac{\Gamma(z+1)}{\Gamma(z)} = z$ and the inequalities clearly hold. 
We now consider the case when $D \ge 2$.
\citet{kershaw_extensions_1983} showed that for any $x > 0$ and $0 < s < 1$, the following must hold:
\begin{align*}
\left( x + \frac{s}{2} \right)^{1-s} < \frac{\Gamma(x+1)}{\Gamma(x+s)} < \left[ x - \frac{1}{2} +\left( s + \frac{1}{4}\right)^{\frac{1}{2}} \right]^{1-s}.
\end{align*}
Since $z \in \mathbb{Z}^+, D\in \mathbb{Z}^+$, we can let $x=z - s > 0 $ and $0 < s=1 - \frac{1}{D} < 1$, and as such: 
\begin{align} \label{q_gamma}
\left( z - \frac{1-\frac{1}{D}}{2} \right)^{\frac{1}{D}} < \frac{\Gamma(z+\frac{1}{D})}{\Gamma(z)} < \left[ z  - \left( \frac{3}{2} -\frac{1}{D} \right) +\left( \frac{5}{4} - \frac{1}{D}\right)^{\frac{1}{2}} \right]^{\frac{1}{D}}.
\end{align}
When $D \ge 2, z\in \mathbb{Z}^+$, 
simple algebra shows that the right hand side of \eqref{q_gamma} is strictly less than $ z^{\frac{1}{D}} $, and therefore $ z^{\frac{1}{D}} - \frac{\Gamma(z+\frac{1}{D})}{\Gamma(z)} > 0$. 

Meanwhile, 
the left hand side of \eqref{q_gamma} leads to $z^{\frac{1}{D}} - \frac{\Gamma(z+\frac{1}{D})}{\Gamma(z)} < z^{\frac{1}{D}} - \left( z - \frac{1-\frac{1}{D}}{2} \right)^{\frac{1}{D}}$.   
Let $\delta(z,D) = z^{\frac{1}{D}} - \left( z - \frac{1-\frac{1}{D}}{2} \right)^{\frac{1}{D}}$. It is easy to show that $\frac{\partial \delta(z,D)}{\partial z} < 0$ for all $z \in \mathbb{Z}^+, D\in \mathbb{Z}^+$, and hence $\delta(z,D)$ should monotonically decrease with $z$ for any $D$; i.e., $ \delta(z,D) \le \delta(1,D)=1- \left( \frac{1+\frac{1}{D}}{2} \right)^{\frac{1}{D}}$. 
For $D\in \mathbb{Z}^+$, this upper bound 
should take its maximum value $1-\frac{\sqrt{3}}{2}$ at $D = 2$. 
Therefore, 
$z^{\frac{1}{D}} - \frac{\Gamma(z+\frac{1}{D})}{\Gamma(z)} <
    1-\frac{\sqrt{3}}{2}.$
This completes the proof.
\hfill \qedsymbol\\

Note from Lemma \ref{lemma:error_gamma} that the approximation error to the gamma functions is always non-negative, and hence, for any $\kappa_1 \le \kappa_2 \in \{0,\ldots,m\}$, the following must hold: 
\begin{align} \label{eq:inequality_kappa}
\mathbb{E}_{\kappa_1} [\bar{C}] \ge \mathbb{E}_{\kappa_2} [\bar{C}] \ge \mathbb{E} [\bar{C}].
\end{align}
This leads to the inequalities as shown in Equation \eqref{eq:E_X_greedy_kappa}. The error bounds of using these approximations are shown in the following lemma.

\begin{lemma} \label{lemma:error_kappa}
When $n \to \infty$, $\kappa \in \{0,\ldots,m\}$,
\begin{align*}
0 \le \mathbb{E}_\kappa [\bar{C}] - \mathbb{E} [\bar{C}] \le 
\left( 1-\frac{\sqrt{3}}{2} \right)R\cdot n^{-\frac{1}{D}},
\end{align*}
and 
$
\lim\limits_{n \rightarrow +\infty} \left( \mathbb{E}_\kappa [\bar{C}] - \mathbb{E} [\bar{C}] \right) \rightarrow 0
$ for $D < \infty$.
\end{lemma}

\noindent{\bf Proof}.
To show $\mathbb{E}_\kappa [\bar{C}] - \mathbb{E} [\bar{C}] \le \left( 1-\frac{\sqrt{3}}{2} \right)
\frac{R}{n^{\frac{1}{D}}}$ for all $\kappa \in \{0,\ldots,m\}$,
it suffices to show it for the extreme case when $\kappa = 0$: 
\begin{align*} 
\begin{split}
    \mathbb{E}_0 [\bar{C}] - \mathbb{E} [\bar{C}]
        & = \frac{R}{mn^{\frac{1}{D}}} \sum_{i=1 }^{m} 
        \sum_{k=1}^{i} \left[ \left( \frac{i-1}{n} \right)^{k-1} \left( 1- \frac{i-1}{n} \right) \left( k^{\frac{1}{D}} - \frac{\Gamma(k+\frac{1}{D}) }{\Gamma(k)}\right) \right. \\
        & +  \left. \left( \frac{i-1}{n} \right)^{i} \left(  i^{\frac{1}{D}} - \frac{\Gamma(i+\frac{1}{D}) }{\Gamma(i)}\right) \right]< \left( 1-\frac{\sqrt{3}}{2} \right)
\frac{R}{n^{\frac{1}{D}}}.
\end{split}
\end{align*}
Note the inequality holds because $
\frac{1}{m}\sum_{i=1 }^{m} 
        \sum_{k=1}^{i} \left[ \left( \frac{i-1}{n} \right)^{k-1} \left( 1- \frac{i-1}{n} \right) \right. +  \left. \left( \frac{i-1}{n} \right)^{i} \right]
=1$, 
and both  $\left(  k^{\frac{1}{D}} - \frac{\Gamma(k+\frac{1}{D}) }{\Gamma(k)}\right)$ and $ \left( i^{\frac{1}{D}} - \frac{\Gamma(i+\frac{1}{D}) }{\Gamma(i)}\right)$ are less than $1-\frac{\sqrt{3}}{2}$ from Equation \eqref{eq:error_gamma}. 

Obviously, for any finite value of $D$, the right hand side of the above inequality approaches 0 as $n$ goes to infinity, and as such $\lim_{n \rightarrow +\infty} \left( \mathbb{E}_0 [\bar{C}] - \mathbb{E} [\bar{C}] \right) \rightarrow 0$.
\hfill \qedsymbol\\

\section{Proof for Lemma \ref{lemma:error_upper_lower}} 
\label{app:error_upper_lower}
We begin by deriving a further relaxed upper bound for $\mathbb{E}_0[\bar{C}]$, as stated in the following lemma.
\begin{lemma} \label{lemma:E_0_ub}
Given $m\in \mathbb{Z}^+$, $n\to \infty$, 
\begin{align} \label{eq:E_0_ub}
\begin{split}
    \mathbb{E}_0 [\bar{C}]
    &\le \frac{R}{n^{\frac{1}{D}}} \left( -\frac{n}{m} \right) \ln \left(1-\frac{m}{n}\right).
\end{split}
\end{align}
\end{lemma}
\noindent{\bf Proof}.
By setting $\kappa=0$ in Equation \eqref{eq:E_X_greedy_kappa}, we have:
\begin{align} \label{app_gamma2}
\begin{split}
    \mathbb{E}_0 [\bar{C}] 
        & = \frac{R}{mn^{\frac{1}{D}}} \sum_{i=1 }^{m} 
        \left[ \sum_{k=1}^{i}  \left( \frac{i-1}{n} \right)^{k-1} \left( 1- \frac{i-1}{n} \right) k^{\frac{1}{D}}  +  \left( \frac{i-1}{n} \right)^{i} i^{\frac{1}{D}}\right].
\end{split}
\end{align}
The first term within the square brackets can be rewritten into the difference between two summations from $k=1$ to $\infty$,
as follows:
\begin{align*} 
\begin{split}
    \sum_{k=1}^{i}  \left( \frac{i-1}{n} \right)^{k-1} \left( 1-  \frac{i-1}{n} \right) k^{\frac{1}{D}}
    = \left( 1- \frac{i-1}{n} \right) \left[ \sum_{k=1}^{\infty} \left( \frac{i-1}{n} \right) ^{k-1} k^{\frac{1}{D}} 
    - \sum_{k=1}^{\infty} \left( \frac{i-1}{n} \right) ^{k-1+i} (k+i)^{\frac{1}{D}} \right].
\end{split}
\end{align*}
Noting $\left( 1-  \frac{i-1}{n}  \right) \sum_{k=1}^{\infty} \left( \frac{i-1}{n} \right) ^{k-1} = 1$, the second term within the square brackets can also be written into a summation from $k=1$ to $\infty$, as follows:
\begin{align*} 
\left( \frac{i-1}{n} \right) ^{i} i^{\frac{1}{D}}
= \left( 1-  \frac{i-1}{n}  \right) \sum_{k=1}^{\infty} \left( \frac{i-1}{n} \right) ^{k-1}\left( \frac{i-1}{n} \right)^{i} i^{\frac{1}{D}}.
\end{align*}
Equation \eqref{app_gamma2} then becomes:
\begin{align} \label{app_uu}
\begin{split}
\mathbb{E}_0[\bar{C}]
        &= \frac{R}{mn^{\frac{1}{D}}} \sum_{i=1}^{m} \left[  
        \left( 1- \frac{i-1}{n} \right) \sum_{k=1}^{\infty} \left( \frac{i-1}{n} \right) ^{k-1} k^{\frac{1}{D}} \right.\\
        & \left. - \left( 1- \frac{i-1}{n} \right) \left( \frac{i-1}{n} \right) ^i \sum_{k=1}^{\infty} \left( \frac{i-1}{n} \right) ^{k-1} \left( (k+i)^{\frac{1}{D}} - i^{\frac{1}{D}} \right)
        \right].
\end{split}
\end{align}
By omitting the second term inside the square brackets in Equation \eqref{app_uu}, and using the poly-logarithm function $\text{Li}_{-\frac{1}{D}}\left( \frac{i-1}{n} \right) = \sum_{k=1}^{\infty} \left( \frac{i-1}{n} \right)^{k} k^{\frac{1}{D}}$ to substitute the first term, we obtain the following formulation:
\begin{align} \label{eq:E_0}
     \bar{\mathbb{E}}_0[\bar{C}] 
     = \frac{R}{mn^{\frac{1}{D}}}  \left[ 1 + \sum_{i=2}^{m} \left( \frac{n}{i-1} -1 \right) \text{Li}_{-\frac{1}{D}}\left( \frac{i-1}{n} \right) \right].
\end{align}

Note here for $0 \le \frac{1}{D} \le 1, k\ge 1, i\ge 1$, it is well-known that $0 < (k+i)^{\frac{1}{D}} \le k^\frac{1}{D} + i^\frac{1}{D}$, where the right hand side equality occurs only at $D = 1$. As such, 
$$\mathbb{E}_0[\bar{C}] \le \bar{\mathbb{E}}_0[\bar{C}].$$
For $D\in \mathbb{Z}^+$, and $u \in (0,1)$, it is well-known that $\frac{u}{(1-u)} = \text{Li}_{0}(u) \le \text{Li}_{-\frac{1}{D}}(u) \le \text{Li}_{-1}(u) = \frac{u}{(1-u)^2}$. 
Let $u = \frac{i-1}{n}$, 
we can get the two inequalities below:
\begin{align} \label{li_bound_1}
m \le 
1 + \sum_{i=2}^{m} \left(\frac{n}{i-1} -1 \right) \text{Li}_{-\frac{1}{D}}\left( \frac{i-1}{n} \right)
\le n\sum_{i=1}^{m} \frac{1}{n-i+1}\le n\ln\left(\frac{n}{n-m}\right).
\end{align}
The last inequality comes from the fact that $\frac{1}{n-i+1}$ monotonically increases with $i$, and hence $\sum_{i=1}^{m} \frac{1}{n-i+1} \le \int_1^{m+1} \frac{1}{n-i+1} \text{ d}i =\ln\left(\frac{n}{n-m}\right)$. 
This directly leads to the upper bound presented in Equation \eqref{eq:E_0_ub}. 
\hfill \qedsymbol\\


Now consider $\mathbb{E}_0[\bar{C}]$ corresponds to the case where Equation \eqref{eq:gamma_ratio} is used to approximate all gamma function ratios in Equation \eqref{eq:E_X_greedy_kappa}. If, instead, we apply Equation \eqref{eq:gamma_ratio} to all terms in Equation \eqref{eq:E_X_greedy_kappa} except the case $k = i = 1$, the resulting value
$$\mathbb{E}_0[\bar{C}] - \left(1 - \frac{\Gamma\left(1 + \frac{1}{D}\right)}{\Gamma(1)}\right)\frac{R}{n^{\frac{1}{D}}}$$
still serves as an upper bound on $\mathbb{E}[\bar{C}]$. Then, based on the relaxed upper bound of $\mathbb{E}_0[\bar{C}]$ given by Lemma \ref{lemma:E_0_ub}, we can derive the following upper bound for $\mathbb{E}[\bar{C}]$.
\begin{align} \label{eq:E_X_ub}
\begin{split}
    \mathbb{E}[\bar{C}] 
    &\le \mathbb{E}_0[\bar{C}]-\left(1-\frac{\Gamma(1+\frac{1}{D})}{\Gamma(1)}\right)\frac{R}{n^{\frac{1}{D}}} 
    \le \bar{\mathbb{E}}_0[\bar{C}]-\left(1-\frac{\Gamma(1+\frac{1}{D})}{\Gamma(1)}\right)\frac{R}{n^{\frac{1}{D}}}\\
    & \le \frac{R}{n^{\frac{1}{D}}} \left[ \left( -\frac{n}{m} \right) \ln \left(1-\frac{m}{n}\right) - 1+\Gamma(1+\frac{1}{D}) \right].
\end{split}
\end{align}
In addition, by setting $m=1$ in Equation \eqref{eq:E_X_greedy_kappa}, we have:
\begin{align} \label{eq:E_X_nearest}
    \mathbb{E}[\bar{C}(1,n)] = \frac{R}{n^{\frac{1}{D}}}\Gamma(1+\frac{1}{D}).
\end{align}

Equations \eqref{eq:E_X_ub} and \eqref{eq:E_X_nearest} jointly lead to the inequality presented in Equation \eqref{eq:error_upper_lower}.
This completes the proof.\\

\end{document}